\documentclass{amsart}

\usepackage{amsmath}
\usepackage{amsfonts}
\usepackage{amssymb}
\usepackage{a4wide}
\usepackage{amsthm}
\usepackage{mathrsfs}
\usepackage{epic}
\usepackage{hyperref}
\usepackage{array,booktabs}
\usepackage{graphicx}
\usepackage{tikz}

\usepackage{mathtools}
\definecolor{green}{rgb}{0.2, 0.8, 0.2}

\usepackage{caption}
\usepackage{pgfplots}
\usepackage{anyfontsize}
\pgfplotsset{compat=newest}
\usetikzlibrary{calc,intersections,datavisualization.formats.functions, patterns}
\usepgfplotslibrary{fillbetween, decorations.softclip} 
\newtheorem{Theorem}{Theorem}
\newtheorem{Proposition}{Proposition}
\newtheorem{Lemma}{Lemma}
\newtheorem{Corollary}{Corollary}

\extrafloats{100}
\theoremstyle{definition}
\newtheorem{Definition}{Definition}
\newtheorem{Remark}{Remark}
\newtheorem{Example}{Example}
\def\komma{\raise.5ex\hbox{,}}
\def\punt{\raise.5ex\hbox{.}}
\newcommand{\Mod}[1]{\ (\mathrm{mod}\ #1)}
\newcommand{\R}{\mathbb R}
\newcommand{\Z}{\mathbb Z}

\newcommand{\N}{\mathbb N}
\newcommand{\T}{T_{\alpha}}
\newcommand{\I}{I_{\alpha}}
\newcommand{\F}{\mathcal F}

\begin{document}
\date{\today}
\title[Orbits of $N$-expansions]{Orbits of $N$-expansions with a finite set of digits}
\author{Jaap de Jonge}
\address{Jaap de Jonge, University of Amsterdam, Korteweg - de Vries Institute for Mathematics, Science Park 105 - 107, 1098 XG Amsterdam, The Netherlands}
\email{c.j.dejonge@uva.nl}
\author{Cor Kraaikamp}
\address{Cor Kraaikamp, Delft University of Technology, department of Electrical Engineering, Mathematics and Computer Science, Mekelweg 4, 2628 CD  Delft, The Netherlands}
\email{C.Kraaikamp@tudelft.nl}
\author{Hitoshi Nakada}
\address{Hitoshi Nakada, Department of Mathematics, Keio University, Yokohama, 223-8522 Japan}
\email{nakada@math.keio.ac.jp}
\date{\today}

\begin{abstract}
For $N \in \N_{\geq 2}$ and $\alpha \in \R$ such that $0 < \alpha \leq \sqrt{N}-1$, we define $I_\alpha:=[\alpha,\alpha+1]$ and $I_\alpha^-:=[\alpha,\alpha+1)$ and investigate the continued fraction map $\T:\I \to \I^-$, which is defined as $\T(x):= N/x-d(x),$ where $d: \I \to \N$ is defined by $d(x):=\left \lfloor N/x -\alpha\right \rfloor$. For $N\in \N_{\geq 7}$, for certain values of $\alpha$, open intervals $(a,b) \subset I_{\alpha}$ exist such that for almost every $x \in I_{\alpha}$ there is an $n_0 \in \N$ for which $\T^n(x)\notin (a,b)$ for all $n\geq n_0$. These \emph{gaps} $(a,b)$ are investigated in the square $\Upsilon_\alpha:=\I \times \I^-$, where the \emph{orbits} $\T^k(x), k=0,1,2,\ldots$ of numbers $x \in \I$ are represented as cobwebs. The squares $\Upsilon_\alpha$ are the union of \emph{fundamental regions}, which are related to the cylinder sets of the map $\T$, according to the finitely many values of $d$ in $\T$. In this paper some clear conditions are found under which $\I$ is gapless. When $\I$ consists of at least five cylinder sets, it is always gapless. In the case of four cylinder sets there are usually no gaps, except for the rare cases that there is one, very wide gap. Gaplessness in the case of two or three cylinder sets depends on the position of the endpoints of $\I$ with regard to the fixed points of $\I$ under $T$.
\end{abstract}

\maketitle

\section{Introduction}
\label{Introduction}

In 2008, Edward Burger and his co-authors introduced in \cite{Bu} new continued fraction expansions, the so-called $N$-expansions, which are nice variations of the \emph{regular continued fraction} (RCF) \emph{expansion}. These $N$-expansions have been studied in various papers since; see \cite{AW}, \cite{DKW} and \cite{K}. In \cite{KL}, a subclass of these $N$-expansions is introduced, for which the digit set is always finite. These particular $N$-expanions are defined as follows:\smallskip

For $N \in \N_{\geq 2}$ and $\alpha \in \R$ such that $0 < \alpha \leq \sqrt{N}-1$, let $I_{\alpha}:=[\alpha,\alpha+1]$ and $I_{\alpha}^-:=[\alpha,\alpha+1)$. Hereafter we denote by $\N_{\geq k}$ the set of positive integers $n\geq k$. We define the $N$-expansion map $\T:I_{\alpha} \to I_{\alpha}^-$ (or $\I$) as
\begin{equation}\label{alphaNCF}
\T(x):=\frac Nx-d(x),
\end{equation}
where $d: I_{\alpha} \to \N$ is defined by 
$$
d(x):=\left \lfloor \frac Nx -\alpha\right \rfloor,\quad {\text{if either $x\in (\alpha,\alpha +1]$ or both $x=\alpha$ and $N/\alpha- \alpha \not \in \Z$}}
$$
and $$
d(\alpha ) = \left\lfloor \frac{N}{\alpha} -\alpha\right\rfloor -1, \quad \text{if $N/\alpha - \alpha \in\Z$}.
$$
Note that if $N/\alpha -\alpha \in\Z$, we have that $T_{\alpha}(\alpha ) = \alpha +1$. This is the only case in which the range of $\T$ is $\I$ and not $\I^-$.

For a fixed $\alpha \in (0,\sqrt{N}-1]$ and $x \in I_{\alpha}$ we define for $n \in \N$
$$
d_n=d_n(x):=d(T_{\alpha}^{n-1}(x)).
$$
Note that for $\alpha \in (0,\sqrt{N}-1]$ fixed, there are only \emph{finitely} many possibilities for each $d_n$. \smallskip

Applying (\ref{alphaNCF}), we obtain for every $x\in I_{\alpha}$ a continued fraction expansion of the form
\begin{equation}\label{expansion}
x=\T^0(x)=\cfrac{N}{d_1+\T(x)}=
\cfrac{N}{d_1+\cfrac{\displaystyle
N}{\displaystyle d_2+\T^2(x)}}=\dots=\cfrac{N}{d_1+\cfrac{\displaystyle
N}{\displaystyle d_2+\cfrac{N}{d_3 + \ddots }}}\quad,
\end{equation}
which we will throughout this paper write as $x=[d_1,d_2,d_3,\ldots ]_{N, \alpha}$ (note that this expansion is infinite for every $x\in I_{\alpha}$, since $0\not\in I_{\alpha}$); we will call the numbers $d_i,\, i \in \N$, the {\it{partial quotients}} or \emph{digits} of this {\it{$N$-continued fraction expansion}} of $x$; see \cite{DKW, KL}, where these continued fractions (also with a finite set of digits) were introduced and elementary properties were studied (such as the convergence in reference \cite{DKW}).\smallskip

In each \emph{cylinder set} $\Delta_i:=\{x \in I_{\alpha};  d(x)=i \}$ of rank $1$, with $d_{\min}\leq i \leq d_{\max}$, where $d_{max}:=d(\alpha )$ is the largest partial quotient, and $d_{min}:=d(\alpha + 1)$ the smallest one given $N$ and $\alpha$, the map $\T$ obviously has one \emph{fixed point} $f_i$ . As of now we will write simply `cylinder set' for `cylinder set of rank $1$'.\smallskip

It is easy to see that\footnote{For reasons of legibility we will usually omit suffices such as `$(N)$', `$(N,\alpha)$' or `$(N,d)$'.} 
\begin{equation}\label{fixed point}
f_i=f_i(N):=\dfrac{\sqrt{4N+i^2}-i}2, \,\,{\text{for}}\,\, d_{\min} \leq i \leq d_{\max}.
\end{equation}
Note that $N/\alpha - \alpha \in\Z$ if and only if for some $d \in \N_{\leq 2}$ we have that $d+1=\max d_i$ for any $\alpha_0<\alpha$, i.e. $\Delta_{d+1}\neq \emptyset$ and $\alpha=f_{d+1}$.\smallskip

Given $N\in \N_{\geq 2}$, we let $\alpha_{max}=\sqrt{N}-1$ be the largest value of $\alpha$ we consider. The reason for this is that for larger values of $\alpha$ we would have $0$ as a partial quotient as well. Since $\T'(x)=-N/x^2$ and because $0<\alpha\leq \sqrt{N}-1$, we have $|\T'(x)|>1$ on $I_{\alpha}^-$. From this it follows that the fixed points act as {\it{repellers}} and that the maps $\T$ are \emph{expanding} when $0<\alpha\leq \sqrt{N}-1$. This is equivalent to the convergence of the $N$-expansion of all $x\in \I$. \smallskip

Each pair of consecutive cylinders sets $(\Delta_i,\Delta_{i-1})$ is divided by a \emph{discontinuity point} $p_i(N,\alpha)$ of $\T$, satisfying $N/p_i-i=\alpha$, so $p_i=N/(\alpha+i)$. A cylinder set $\Delta_i$ is called \emph{full} if $\T(\Delta_i)=I_{\alpha}^-$ (or $\I$). When a cylinder set is not full, it contains either $\alpha$ (in which case $\T(\alpha)<\alpha+1$) or $\alpha+1$ (in which case $\T(\alpha+1)>\alpha$), and is called \emph{incomplete}. On account of our definition of $T_{\alpha}$, cylinder sets will always be an interval, \emph{never} consist of one single point.

The main object of this paper is the sequence $\T^n(x)$, $n=0,1,2,\ldots$, for $x\in I_{\alpha}$, which is called the \emph{orbit of $x$ under $\T$}. More specifically, we are interested in subsets of $\I$ that we will call \emph{gaps} for such orbits. Before we will give a proper definition of `gap', we will give an example of orbits of points in $\I$ for a pair $\{N,\alpha\}$. Note that, due to the repellence of the fixed points, orbits cannot remain in one cylinder set indefinitely when $x\in \I$ is not a fixed point of $\T$, so any orbit will show an infinite migration between cylinder sets. A naive approach is to compute the orbits of many points of $\I$ and obtain a \emph{plot} of the asymptotic behaviour of these orbits by omitting the first, say hundred, iterations. Figure \ref{fig: N=51, alpha=6} shows such a plot for $N=51$ and $\alpha=6$. It appears that there are parts of $\I$ (illustrated by dashed line segments) that are not visited by any orbit after many iterations of $\T$. 

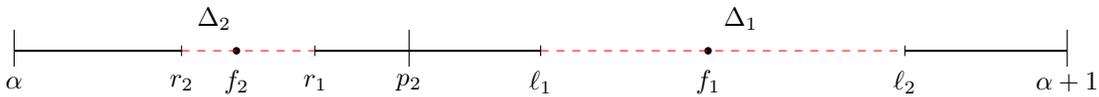
\begin{figure}[!htb]
$$
\begin{tikzpicture}[scale =14] 
\draw[black,fill=black] (.2111,0) circle (.02ex);
\draw[black,fill=black] (.659,0) circle (.02ex);
  \draw [thick] (0,0) -- (0.159,0);
  \draw [dashed,red] (0.159,0) -- (0.2857,0);
  \draw [thick] (0.2857,0) -- (0.5,0);
   \draw [dashed,red] (0.5,0) -- (0.846,0);
\draw [thick] (0.846,0) -- (1,0);
  \draw (.375,-.015) -- (.375,.02);
 \draw (1,-.015) -- (1,.02);
  \draw (0,-.015) -- (0,.02);
   \draw (.159,-.005) -- (.159,.005);
    \draw (.2857,-.005) -- (.2857,.005);
 \draw (.846,-.005) -- (.846,.005);
 \draw (.5,-.005) -- (.5,.005);
\node at (.69,.03) {$\Delta_1$};
\node at (.19,.03) {$\Delta_2$};
\node at (0,-.03) {$\alpha$};
\node at (1,-.03) {$\alpha+1$};
\node at (.375,-.03) {$p_2$};
\node at (.659,-.03) {$f_1$};
\node at (.2111,-.03) {$f_2$};
\node at (.159,-.03) {$r_2$};
\node at (.286,-.03) {$r_1$};
\node at (.5,-.03) {$\ell_1$};
\node at (.846,-.03) {$\ell_2$};
        \end{tikzpicture}
 $$
 \vspace*{-8mm}  \caption{$N=51$, $\alpha=6$ \label{fig: N=51, alpha=6}}
\end{figure}

In fact, setting $\ell_i=\T^i(\alpha)$ and $r_i =\T^i(\alpha +1)$, the orbit of any point -- apart from the fixed points $f_1$ and $f_2$ -- after once having left the interval $(r_2,r_1)\subset \Delta_2$ or $(\ell_1,\ell_2)\subset \Delta_1$ of Figure \ref{fig: N=51, alpha=6}, will \textbf{never} return to it.\smallskip

In order to get a better understanding of the orbits of $N$-expansions, it is useful to consider the graphs of $\T$, which are drawn in the square $\Upsilon_{N,\alpha}:=I_{\alpha} \times I_{\alpha}^- $. This square is divided in rectangular sets of points $\Box_i:=\{(x,y)\in \Upsilon_\alpha :d(x) =i\}$, which are the two-dimensional \emph{fundamental regions} associated with the one-dimensional cylinder sets we already use. We will call these regions shortly \emph{cylinders}. Now consider $(x,\T(x)) \in \Upsilon_{N,\alpha}$. Then $(x,\T(x))$ goes to $(\T(x),\T^2(x))$ under $\T$. Regarding this, $\T$ has one fixed point $F_i:=(f_i,f_i)$ in each $\Box_i$. We will denote the dividing line between $\Box_i$ and $\Box_{i-1}$ by $l_i$, which is the set $\{p_i\} \times [\alpha,\alpha+1)$, with $p_i$ the discontinuity point between $\Delta_i$ and $\Delta_{i-1}$. In case $\T(\Delta_i)=I_{\alpha}^- $, we will call the cylinder $\Box_i$ {\it{full}} and the branch of the graph of $\T$ in $\Box_i$ {\it{complete}}; if a cylinder is not full, we will call it and its associated branch of $\T$ {\it{incomplete}}. We will call the collection of $\Upsilon_{\alpha}$ and its associated branches, fixed points and dividing lines an \emph{arrangement} of $\Upsilon_{\alpha}$. When $\Upsilon_{\alpha}$ is a union of full cylinders, we will call the associated arrangement also full.\smallskip

Figure \ref{fig: N=51, alpha=6.5, x} is an example of such an arrangement, in which a part of the \emph{cobweb} is drawn associated with the orbit we investigated previously. The discontinuity point $p_2=51/8$ is now visible as a dividing line between $\Delta_1$ and $\Delta_2$. 
\begin{figure}[!htb]
$$
\begin{tikzpicture}[scale =6] 
\draw[black,fill=black] (.659,.659) circle (.05ex);
\draw[black,fill=black] (.211,.211) circle (.05ex);
           \draw [dashed] (.375,0) -- (.375,1);
             \draw (0,0) -- (1,1);
     \draw (0,0) -- (0,1);
        \draw (0,0) -- (1,0);
     \draw (0,1) -- (1,1);
     \draw (1,1) -- (1,0);
       \draw [->, blue, thick] (.5,0) -- (.5,.423);
          \draw [->, blue, thick] (.5,.423) -- (.5,.846);
     \draw [->, blue, thick] (.5,.846) -- (.67,.846);
     \draw [->, blue, thick] (.67,.846) -- (.846,.846);
     \draw [->, blue, thick] (.846,.846) -- (.846,.64);
     \draw [->, blue, thick] (.846,.64) -- (.846,.449);
     \draw [->, blue, thick] (.846,.449) -- (.64,.449);
     \draw [->, blue, thick] (.64,.449) -- (.449,.449);
     \draw [->, blue, thick] (.449,.449) -- (.449,.68);
     \draw [->, blue, thick] (.449,.68) -- (.449,.908);
     \draw [->, blue, thick] (.449,.908) -- (.68,.908);
     \draw [->, blue, thick] (.68,.908) -- (.908,.908);
     \draw [->, blue, thick] (.908,.908) -- (.908,.645);
     \draw [->, blue, thick] (.908,.645) -- (.908,.383);
     \draw [->, blue, thick] (.908,.383) -- (.645,.383);
     \draw [->, blue, thick] (.645,.383) -- (.383,.383);
     \draw [->, blue, thick] (.383,.383) -- (.383,.68);
   \draw [->, blue, thick] (.383,.68) -- (.383,.99);
     \draw [->, blue, thick] (.383,.99) -- (.68,.99);
     \draw [->, blue, thick] (.68,.99) -- (.99,.99);
     \draw [->, blue, thick] (.99,.99) -- (.99,.64);
     \draw [->, blue, thick] (.99,.64) -- (.99,.296);
     \draw [->, blue, thick] (.99,.296) -- (.64,.296);
      \draw [->, blue, thick] (.64,.296) -- (.296,.296);
          \draw [->, blue, thick] (.296,.296) -- (.296,.2);
     \draw [->, blue, thick] (.296,.2) -- (.296,.1);
     \draw [->, blue, thick] (.296,.1) -- (.2,.1);
      \draw [->, blue, thick] (.2,.1) -- (.1,.1);
               \node at (.72,.67) {$F_1$};
      \node at (.1414,.22) {$F_2$};
      \node at (.19,.93) {$\Box_2$};
    \node at (.69,.93) {$\Box_1$};
        \node at (-.02,-.04) {$6$};
    \node at (1,-.04) {$7$};
     \node at (-.03,1) {$7$};
    \node at (.375,-.04) {$\tfrac{51}8$};
    \node at (.5,-.04) {$x$};
     \draw[domain=.375:1,smooth,variable=\x] plot ({\x},{51/(\x+6)-7});
   \draw[domain=0:.375,smooth,variable=\x] plot ({\x},{51/(\x+6)-8});
    \end{tikzpicture}
  $$
\vspace*{-8mm} \caption{$N=51$, $\alpha=6$, $x=6.5$ \label{fig: N=51, alpha=6.5, x}}
\end{figure}
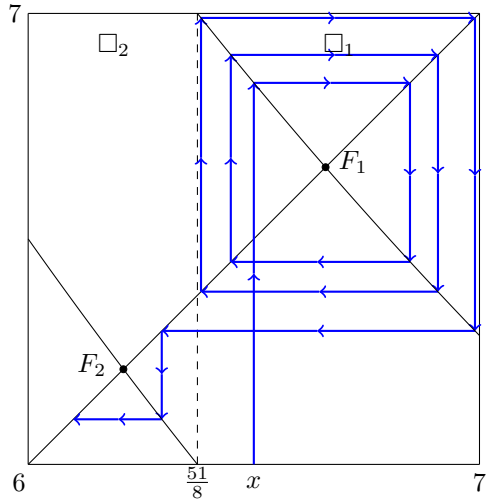

In \cite{DKW} and \cite{KL} the arrangement for $N=4$ and $\alpha=1$ is studied, consisting of two full cylinders $\Box_1$ and $\Box_2$ and not showing any gaps. On the other hand, the demonstration of the interval $(5/2,13/5)$ being a gap of the interval $[2,3]$ in the case $(N,\alpha)=(9,2)$ in \cite{KL} is done without referring to such an arrangement. In this paper, and even more so in the next paper, we will show that arrangements may considerably support the insight in the occurrence of gaps.\smallskip

We will now give a formal definition of gaps, which is slightly delicate, since $T_{\alpha}(I_{\alpha}) = I_{\alpha}^-$ (or $\I$ when $N/\alpha - \alpha \in\Z$). \smallskip

\begin{Definition}
A maximal open interval $(a,b) \subset I_{\alpha}$ is called a \emph{gap} of $I_{\alpha}$ if for almost every $x \in I_{\alpha}$ there is an $n_0 \in \N$ for which $\T^n(x)\notin (a,b)$ for all $n\geq n_0$.
 \end{Definition}
 
\begin{Remark} In the example of Figure \ref{fig: N=51, alpha=6} the intervals $(r_2,r_1)$ and $(\ell_1,\ell_2)$ are gaps and for $x\in (r_2,r_1)\cup(\ell_1,\ell_2)\setminus \{f_1,f_2\}$ there exists an $n_0 = n_0(x)$ such that $\T^n(x)\not\in (r_2,r_1)\cup(\ell_1,\ell_2)$ for $n\in \N_{\geq2}$. The `for almost'\footnote{All `for all' statements in this paper are with respect to Lebesgue measure.} formulation in the definition of `gap' is necessary so as to exclude fixed points and pre-images of fixed points, i.e.~points that are mapped under $\T$ to a fixed point, which may never leave an gap. In Section \ref{four and five} we even find a class of gaps $(a,b)$ such that for uncountably many $x\in \I$ and all $n\in \N\cup\{0\}$ we have $\T^n(x)\in (a,b)$.
\end{Remark}

\begin{Remark}
When we use the word `gap' in relation to arrangements, we mean the gap of the associated interval $\I$.
\end{Remark}
 
\begin{figure}[!htb]
\includegraphics[width=80mm]{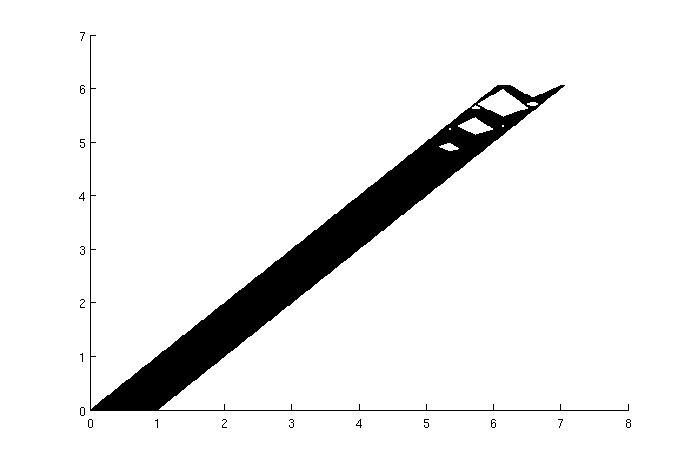}
\vspace*{-6mm} \caption{A simulation of intervals $I_{\alpha}$ with gaps if existent, for $0<\alpha\leq \sqrt{50}-1$ and $N=50$\label{N=50}}
\end{figure}

In \cite{KL} computer simulations were used to get a more general impression of orbits of $N$-expansions. For a lot of values of $\alpha$, plots such as Figure \ref{fig: N=51, alpha=6} were stacked, for $0<\alpha\leq \alpha_{\max}$, so as to obtain graphs such as Figure \ref{N=50}, with the values of $\alpha$ on the vertical axis and at each height the corresponding interval $\I$ drawn. In the same paper, similar graphs are given for $N=9,20,36$ and $100$. In all cases it appears that `gaps' such as in Figure \ref{fig: N=51, alpha=6} appear for values of $\alpha$ equal to or not much smaller than $\alpha_{\max}$. Since the plots in \cite{KL} are based on computer simulations, they do not actually show very small gaps (smaller than pixel size) nor clarify much the connection between the gaps for each $N$. Still, the suggestion is strong that for $\alpha$ sufficiently small there are no gaps. It also seems that for $\alpha$ large enough several disjoint gaps may occur. In Figure \ref{N=50} we see this for $\alpha$ near $\alpha_{\max}=\sqrt{50}-1$. \smallskip

In this paper we will not only investigate conditions for gaplessness, we will also show that simulations such as Figure \ref{N=50} fail to reveal the existence, for certain $N$ and $\alpha$, of one extremely large gap in plots such as Figure \ref{N=50} \emph{below} the last visible gap. In a subsequent paper we will go into another very interesting property of orbits of $N$-expansions that is hardly revealed by simulations such as Figure \ref{N=50}: the existence of large numbers of gaps for large $N$ and $\alpha$ close to $\alpha_{\max}$. But now we will concentrate on gaplessness. 

\begin{Remark} 
When no gaps exist with non-empty intersection with a cylinder set, we call the cylinder set gapless. 
\end{Remark} 

In Section \ref{Full arrangements} we will consider two classes of arrangements that have no gaps: full arrangements and specific arrangements with more than three cylinders. The gaplessness of the latter class, involving the proof of Theorem \ref{no gaps five cylinder case}, for which some preliminary results will be presented shortly, is largely given in Section \ref{Full arrangements}, but involve some intricacies for small values of $N$ so as to finish it at the end of this paper. In Section \ref{two-cylinder cases} we will consider arrangements with two cylinders and in Section \ref{three or more} we will concentrate on arrangements with three cylinders, but will prove a sufficient condition for gaplessness that is valid for arrangements with any number of cylinders larger than $2$. Finally, in Section \ref{four and five} we will prove a result on gaps in certain arrangements with four cylinders and we will finish the proof of Theorem \ref{no gaps five cylinder case}. After that, it is merely a matter of checking that for $N\in \{2,3,4,5,6\}$ all arrangements are gapless. \smallskip

\section{Full arrangements and arrangements with more than four cylinders}
\label{Full arrangements}

When $I_{\alpha}$ consists of full cylinder sets only, we obviously have no gaps. In this situation the mutual relations between $N$, $\alpha$ and $d(\alpha)$ show a great coherence, as expressed in the following theorem:
\begin{Theorem}\label{when I consists of m full cylinder sets}
The interval $I_{\alpha}$ consists of $m$ full cylinder sets, with $m \in \N_{\geq 2}$, if and only if there is a positive integer $k$ such that
\begin{equation}\label{three relations}
\begin{cases}
\alpha=k,\\
N=mk(k+1),\\
d(\alpha)=(m-1)(k+1).
\end{cases}
\end{equation}
\end{Theorem}
 Proof of Theorem \ref{when I consists of m full cylinder sets}: Writing $d:=d(\alpha)$, the interval $I_\alpha$ is the union of $m$ full cylinder sets if and only if
\begin{equation}\label{alfa geheel 1}
\begin{cases} \frac N{\alpha} - d = \alpha +1, \\
\frac N{\alpha+1}-(d-m+1)=\alpha.
\end{cases}
\end{equation}
Note that the first equation in (\ref{alfa geheel 1}) can be written as
$$
N = \alpha^2 + (d+1)\alpha,
$$
while the second equation in (\ref{alfa geheel 1}) equals
$$
N = \alpha^2 + (d+2-m)\alpha + d + 1 - m.
$$
Subtracting the first of these equations from the last we find
\begin{equation}\label{alfa integer}
\alpha = \frac{d+1-m}{m-1}.
\end{equation}
From (5), we have
\begin{equation}\label{eq:alphaASfraction}
\alpha = \frac{-(d+1) + \sqrt{(d+1)^2+4N}}{2},
\end{equation}
which yields that $\alpha$ is either a quadratic irrational or a rational number. Since (6) implies that $\alpha$ is a rational number we find that the integer $(d+1)^2+4N$ must be a square, i.e.\ there exists a positive integer $s$ such that $s^2=(d+1)^2+4N$. Note that $d$ is an even integer if and only if $s^2$ is an odd integer if and only if $s$ is an odd integer. Consequently we find that the numerator of $\alpha$ in~(\ref{eq:alphaASfraction}) is \emph{always} even, and~(\ref{eq:alphaASfraction}) yields that $\alpha$ is a positive integer, say $k$. From the equations in (\ref{alfa geheel 1}) it follows that not only $\alpha = k$ but also $\alpha +1 = k+1$ is a divisor of $N$.\medskip\

From the definition of $\T$ in (\ref{alphaNCF}) (especially the case $N/\alpha-\alpha \in \Z$) we see that
\begin{equation}\label{extraequationTHeorem1}
d = d(\alpha ) = \frac{N}{k} - (k+1).
\end{equation}
On the other hand (\ref{alfa integer}) yields that, since $\alpha =k$,
$$
d = (m-1)(k+1),
$$
and from this and~(\ref{extraequationTHeorem1}) we see that
$$
\frac{N}{k} - (k+1) = (m-1)(k+1),
$$
i.e.~$N = mk(k+1)$.\medskip\

Conversely, let $k$ be a positive integer such that the relations of (\ref{three relations}) hold. Then both $N/ \alpha$ and $N/ (\alpha +1)$ are positive integers, implying that \emph{all} cylinder sets are full. Moreover, since $d=d(\alpha )= d_{\max}$ is given by
$$
d = \frac{N}{\alpha} - \alpha -1 = \frac{mk(k+1)}{k} - k - 1 = (m-1)(k+1),
$$
and $d_{\min} = d(\alpha +1)$ is given by
$$
d(\alpha +1) = \left\lfloor \frac{N}{\alpha +1} - \alpha \right\rfloor = mk - k = (m-1)k,
$$
it follows that there are
$$
d_{\max} - \left( d_{\min} -1\right) = (m-1)(k+1) - (m-1)k +1  = (m-1) + 1 = m
$$
full cylinder sets.\hfill $\Box$

Theorem \ref{when I consists of m full cylinder sets} serves as a starting point for our investigation of orbits of $N$-expansions. The first thing we will do is give some preliminary results (in Subsection \ref{preliminaries}) that we need for proving (in Subsection \ref{large enough branch number}) Theorem \ref{Wilkinson} and Theorem \ref{no gaps five cylinder case} on gaplessness of arrangements with at least five cylinders.  

\subsection{Preliminary results}\label{preliminaries}\smallskip

The first thing to pay attention to is the way $N$ and $\alpha$ and $d(\alpha)$, the value of the largest partial quotient, are interdependent, which is illustrated by the following lemmas:

\begin{Lemma}\label{N and d depending on a}
Given $N$ and $\alpha$, let $d := d(\alpha )$ be the largest possible digit. Then
$$
d\geq N-1\,\, {\text{if and only if}}\,\, \alpha<1.
$$
\end{Lemma}
The proof of this lemma is left to the reader. \smallskip

When $\alpha=\alpha_{\max}=\sqrt{N}-1$, we have 
\begin{align}\label{variance of d}
d(\alpha)=\begin{cases} 
\left \lfloor \frac{2}{\sqrt{2}-1}-(\sqrt{2}-1)\right \rfloor=4& {\text{for}}\,\,N=2;\medskip\\
\left \lfloor \frac{3}{\sqrt{3}-1}-(\sqrt{3}-1)\right \rfloor=3& {\text{for}}\,\,N=3;\medskip\\
\left \lfloor \frac{4}{\sqrt{4}-1}-(\sqrt{4}-1)\right \rfloor-1=2& {\text{for}}\,\,N=4;\medskip\\
\left \lfloor \frac{N}{\sqrt{N}-1}-(\sqrt{N}-1)\right \rfloor=\left \lfloor 2+\frac{\sqrt{N}+1}{N-1}\right \rfloor=2& {\text{for}}\,\, N\in \N_{\geq 5}.\medskip\\
\end{cases}
\end{align}
On the other hand we have, for $N\in\N$, $N\geq 2$ fixed:
$$
\lim_{\alpha \downarrow 0} d(\alpha)=\lim_{\alpha \downarrow 0} \left \lfloor \frac{N}{\alpha}-\alpha\right \rfloor=\infty.
$$
The following lemma provides for a lower bound for the rate of increase of $d(\alpha)$ compared with the rate of decrease of $\alpha$. 
\begin{Lemma}\label{d increases more than twice as fast as alpha decreases}
Let $N\in \N_{\geq 2}$ be fixed and $d:=d(\alpha )$. Then $d$ is constant for $\alpha \in [f_{d+1},f_d)$, and $d$ increases overall more than twice as fast as $\alpha$ decreases.
\end{Lemma}

Proof of Lemma \ref{d increases more than twice as fast as alpha decreases}: Starting from $\alpha_{\max}$, $d$ increases by $1$ each time $\alpha$ decreases beyond a fixed point, i.e. when $N/\alpha - \alpha \in \N$. For the difference between two successive fixed points $f_{d-1}$ and $f_d$ we have
$$
f_{d-1}-f_d=\frac{\sqrt{4N+(d-1)^2}-(d-1)}2-\frac{\sqrt{4N+d^2}-d}2=\frac{\sqrt{4N+(d-1)^2}-\sqrt{4N+d^2}+1}2<\frac12.
$$
This finishes the proof. \hfill $\Box$\smallskip

Closely related to the previous lemma is the following one, the proof of which is left to the reader.
\begin{Lemma}\label{distance between two fixed points for fixed d} 
Let $d \in \N_{\geq 2}$ and $N\in \N_{\geq 2}$ be fixed and let $f_d(N)$ be defined by the equation $N/f_d(N) - d = f_d(N)$ (so $f_d(N)$ is the fixed point of the map $x \mapsto N/x - d$ for $x\in (0,N/d)$). Then
$$
f_{d-1}(N+1)-f_d(N+1)>f_{d-1}(N)-f_d(N).
$$
\end{Lemma}
So, for $d$ fixed, the distance between two consecutive fixed points increases when $N$ increases. We have, in fact, for $d$ fixed:
$$
\lim_{N\to\infty} \left( f_{d-1}(N) - f_d(N)\right) = \tfrac{1}{2};
$$
cf.\ the proof of Lemma \ref{d increases more than twice as fast as alpha decreases}. For $N$ fixed, on the other hand, we have:
$$
\lim_{d\to\infty}  f_d(N) = 0. 
$$

While $d(\alpha)$ is a monotonously non-increasing function of $\alpha$, the number of cylinder sets is not. The reason is obvious: starting from $\alpha=\alpha_{\max}$, the number of cylinder sets changes every time either $\alpha$ or $\alpha+1$ decreases beyond the value of a fixed point; in the first case, the number increases by $1$, and in the second case, it decreases by $1$. Since $\T'(x)=-N/x^2<0$ and $\T''(x)=2N/x^3>0$ on $\I$, $\T(x)$ is decreasing and convex on $\I$, implying that a per saldo increase of the number of cylinder sets. Still, for $N$ and $\alpha$ large enough, it may take a long time of $\alpha$ decreasing from $\alpha_{\max}$ before the amount of cylinder sets stops alternating between two successive numbers $k\in \N_{\geq 2}$ and $k+1$, and starts to alternate between the numbers $k+1$ and $k+2$. As an example, we take $N=100$. When $\alpha$ decreases from $\alpha_{\max}$, the interval $\I$ consists of two cylinder sets until $\alpha$ decreases beyond $f_3$ and cylinder set $\Delta_3$ emerges; then, when $\alpha+1$ decreases beyond $f_1$, cylinder set $\Delta_1$ disappears and so on, until $\alpha$ decreases beyond $f_8$ and $\Delta_9$ emerges while $\Delta_6$ has not yet disappeared. \smallskip

In order to get a grip on counting the number of cylinder sets, the following arithmetic will be useful: a full cylinder set counts for 1, an incomplete left one counts for $N/\alpha-d_{\max}-\alpha$, and an incomplete right one for $\alpha+1-(N/(\alpha+1)-d_{\min})$, giving rise to the following definition: \smallskip

\begin{Definition}\label{branch number}
Let $N \in \N_{\geq 2}$ and $\alpha \in \R$ such that $0 < \alpha \leq \sqrt{N}-1$ and $\T$ the $N$-continued fraction map. The branch number\footnote{The word `branch' refers to the part of the graph of $\T$ on the concerning cylinder set.} $b(N,\alpha)$ is defined as 
\begin{align*}
b(N,\alpha):= &\,d_{\max}-d_{\min}-1\,\,{\text{(the number of full cylinder sets save for the outermost ones)}}  \\
+\,\,&\frac N{\alpha}-d_{\max} -\alpha\,\,{\text{(the length of the image of the leftmost cylinder set)}} \\
+\,\,&\alpha+1-\left (\frac N{\alpha+1}-d_{\min}\right )\,\,{\text{(the length of the image of the rightmost cylinder set),}}
\end{align*}
\end{Definition}
From this the next lemma follows immediately:
\begin{Lemma}\label{branch equation}
For $N\in \N_{\geq 2}$ and $0<\alpha\leq \sqrt{N}-1$ we have
$$
b(N,\alpha)=\frac N{\alpha}-\frac N{\alpha+1}=\frac{N}{\alpha(\alpha+1)}. 
$$
\end{Lemma}
It follows that for fixed $N$ the branch number $b(N,\alpha)$ is a strictly decreasing function of $\alpha$. 

\begin{Remark}
Applying Lemma \ref{branch equation}, we find 
\begin{equation}\label{branch alphamax}
b(N,\alpha_{\max})=\frac{N}{(\sqrt{N}-1)\sqrt{N}}=1+\frac1{\sqrt{N}-1}.
\end{equation}
It follows that $b(N,\alpha)>1$ for all $N\in \N_{\geq 2}$, so the number of cylinder sets is always at least 2. On the other hand, from Lemma \ref{branch equation} it follows that the number of cylinder sets increases to infinity as $\alpha$ decreases from $\alpha_{max}$ to $0$. Actually, we have infinitely many digits \emph{only} when $\alpha = 0$. In this case the corresponding $N$-expansion is the \emph{greedy} $N$-expansion, studied in \cite{AW} and \cite{DKW}.\smallskip
\end{Remark}

The relation $N/(\alpha(\alpha+1))=b$ yields
\begin{equation}\label{formula a in b}
\alpha=\frac{\sqrt{\frac{4N}b +1}-1}2,
\end{equation}
from which we derive that $d(\alpha) =d$ (or $d(\alpha )=d-1$  in case $N/\alpha - \alpha \in\Z$), where $d$ is given by
\begin{equation}\label{formula d in b}
d=\left \lfloor \frac {(b-1)\sqrt{\frac{4N}b +1}+b+1}2\right \rfloor.
\end{equation}

\subsection{Gaplessness when the branch number is large enough}\label{large enough branch number}\smallskip

So far, we merely discussed the way $\I$ is divided in cylinder sets, depending on the values of $N, \alpha, d(\alpha)$ and the branch number $b$. In order to present some first results on sufficient conditions for gaplessness, we will zoom in on some ergodic properties of $\T$.  

\begin{Lemma}\label{before ergodic}
If $\mu$ is an absolutely continuous invariant probability measure for $\T$, then there exists a function $h$ of bounded variation such that 
$$
\mu (A) = \int_{A} h\, d \lambda, \,\, \lambda-\mbox{a.e.},\,\, \mbox{with $\lambda$ the Lebesgue measure},
$$
i.e.~ any absolutely continuous invariant probability measure has a version of its density function of bounded variation.
\end{Lemma}
Proof of Lemma \ref{before ergodic}: Since $\inf |\T'| > 1$, applying Theorem 1 from \cite{LaY} immediately yields the assertion. \hfill$\Box$  

\begin{Theorem}\label{ergodic}
Let $N \in \N_{\geq2}$. Then there is a unique absolutely continuous invariant probability measure $\mu_\alpha$ such that $\T$ is ergodic with respect to $\mu_\alpha$. 
\end{Theorem}
Proof\footnote{see also page 185, Theorem 1 in \cite{LiY}} of Theorem \ref{ergodic}: Let $\mu_\alpha$ be a unique absolutely continuous invariant probability measure for $\T$ and choose its density function $h$ of bounded variation. Then there exists an open interval $J$ such that $h(x) > 0$ for any $x \in J$, since $h$ has at most countably many discontinuity points. Consider $\{T^{n}J \, : \, n \ge 0\}$.  Since $\inf |T'|>1$, there exists an $n_{0}$ such that $T^{n_0}(J)$ includes a discontinuity point.  (If necessary we may choose endpoints of $J$ not in the preimages of discontinuity points of $\T$.) We note that for any measurable subset  $A \subset J$ with $\mu_\alpha(A)>0$ equivalently $\lambda(A)>0$, $\mu(T^{n}A)>0$ for any $n \ge 1$.  Now $T^{n_{0}+1}(J)$ includes two intervals $J_{\ell}$ and $J_{r}$ attached to $\alpha$ and $\alpha + 1$ respectively.  For any measurable subset $B_{0}\subset J_{\ell} \cup J_{r}$ of positive 
$\lambda$-measure, $\mu(B_{0}) > 0$, since otherwise we have a contradiction; $\mu(B_{0}) = 0$ and $\mu(T^{-(n_{0}+1)}(B_{0}) )> 0$ (since there is a $B_{1} \subset J$ such that $T^{n_{0}+1}(B_{1}) = B_{0}, \mu(B_{1})>0$).  This shows that any two absolutely continuous invariant probability measures $\mu_{1}$ and $\mu_{2}$ cannot have disjoint supports (i.e.~they cannot be singular to each other), which is equivalent to the uniqueness of the absolutely continuous invariant probability measure and hence its ergodicity.  \hfill$\Box$

The next result follows directly from Theorem \ref{ergodic}:
\begin{Corollary}\label{after ergodic}
If iteration of $\T$ maps all open subintervals of $I_{\alpha}$ to the interval $I_{\alpha}^{-}$, then $\I$ contains no gaps. 
\end{Corollary}
Proof of Corollary \ref{after ergodic}: The assumption implies that the absolutely continuous invariant probability measure $\mu_\alpha$ is equivalent to the Lebesgue measure, which implies that for any measurable subset 
$A \subset I_{\alpha}$, $\mu(A) =0$ if and only if $\lambda(A)=0$.  Suppose that there is a gap $J$.  Since $J$ is an open interval, we have $\lambda (J) >0$, thus  $\mu(J)>0$. Since $\mu(I_{\alpha})<\infty$ implies a.e.~$x \in J$, there exists infinitely many positive integers $n$ such that $T^{n}(x) \in J$ (by the Poincar\'{e} recurrence theorem), which contradicts the assumption that there is a gap. \hfill$\Box$

Before we present the first of two theorems on gaplessness, we note that in the case $N=2$, the condition $|\T'(x)|>2$ for all $x \in I_{\alpha}$ is not satisfied for any $\alpha \in (0,\sqrt{2}-1]$.

\begin{Theorem}\label{Wilkinson}
Let $N\in \N_{\geq 3}$, and let $0< \alpha \leq \sqrt{N}-1$. Let $|\T'(x)|> 2$ for all $x \in I_{\alpha}$. Then $I_{\alpha}$ contains no gaps.
\end{Theorem}

Proof of Theorem \ref{Wilkinson}: The condition implies $N/(\alpha+1)^2>2$, yielding $\alpha < \sqrt{N/2}-1$. From Lemma \ref{branch equation} it follows that 
$$
b(N,\alpha)> \frac {2\sqrt{2N}}{\sqrt{2N}-2},
$$
which is larger than $2$ for all $N \in \N_{\geq 3}$. So $I_{\alpha}$ consists of at least three cylinder sets. Since $|T_{\alpha}^{\prime}(x)| > 2$ for all $x\in \I$, there exists an $\varepsilon >0$ such that for any open interval $J_0$ that is contained in a cylinder set of $\T$ we have
$$
\left| \T(J_0)\right| \geq (2+\varepsilon )|J_0|,
$$
where $|J|$ denotes the length (i.e. Lebesgue measure) of an interval $J$.\smallskip

If $\T(J_0)$ contains two consecutive discontinuity points $p_{i+1},p_i$ of $\T$, then
$$
(p_{i+1},p_i)\subset \T(J_0),
$$
and we immediately have that
$$
I^O_{\alpha}:=(\alpha ,\alpha+1)=\T(p_{i+1},p_i)\subset \T^2(J_0).
$$
If $\T(J_0)$ contains only one discontinuity point $p$ of $\T$, then $\T(J_0)$ is the disjoint union of two subintervals located in two adjacent cylinder sets:
$$
\T(J_0)=J_1^{\prime}\cup J_2^{\prime}.
$$
Obviously, 
$$
|\T(J_0)| = |J_1^{\prime}| + |J_2^{\prime}|.
$$
Now select the larger of these two intervals $J_1^{\prime}$, $J_2^{\prime}$, and call this interval $J_1$. Then
$$
|J_1| \geq (1+\tfrac{\varepsilon}{2})|J_0|.
$$
In case $\T(J_0)$ does not contain any discontinuity point of $\T$, we set $J_1=\T(J_0)$. Induction yields that there exists an $\ell \in \N$ such that
$$
|J_{\ell}| \geq \left( 1 + \frac{\varepsilon}{2}\right)^{\ell} \left| J_0\right|,
$$
whenever $\T(J_{\ell - 1})$ includes no more than one discontinuity point of $\T$. But then there must exist a $k \in \N$ such  that $\T(J_k)$ contains two (or more) consecutive discontinuity points of $\T$, and we find that $\T^2(J_k)=I_{\alpha}^O$. Applying Corollary \ref{after ergodic}, we conclude that there is no gap in $\I$. \hfill$\Box$\smallskip

The next theorem, which is partly a corollary of the previous one, gives an even more explicit condition for gaplessness. 
\begin{Theorem}\label{no gaps five cylinder case}
Let $\I$ consist of five cylinder sets or more. Then $\I$ has no gaps.
\end{Theorem}

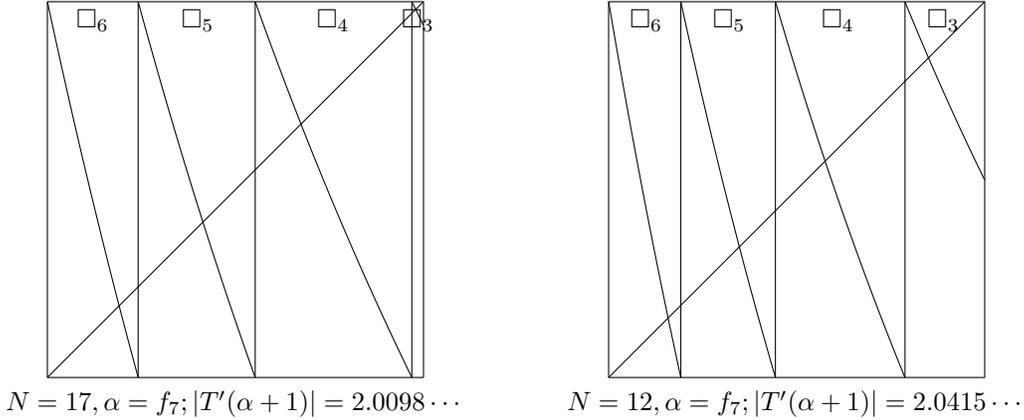
\begin{figure}[!htb]
\minipage{0.5\textwidth}
$$
\begin{tikzpicture}[scale =5] 
      \draw (.24134,0) -- (.24134,1);
       \draw (.5525,0) -- (.5525,1);
        \draw (.969,0) -- (.969,1);
                     \draw (0,0) -- (1,0);
             \draw (0,0) -- (1,1);
     \draw (0,0) -- (0,1);
     \draw (0,1) -- (1,1);
     \draw (1,1) -- (1,0);
        \node at (.12,.95) {$\Box_6$};
              \node at (.4,.95) {$\Box_5$};
    \node at (.76,.95) {$\Box_4$};
              \node at (.985,.95) {$\Box_3$};
               \draw[domain=.969:1,smooth,variable=\x] plot ({\x},{17/(\x+1.9083)-4.9083});
             \draw[domain=.5525:.969,smooth,variable=\x] plot ({\x},{17/(\x+1.9083)-5.9083});
                  \draw[domain=.24134:.5525,smooth,variable=\x] plot ({\x},{17/(\x+1.9083)-6.9083});
             \draw[domain=0:.24134,smooth,variable=\x] plot ({\x},{17/(\x+1.9083)-7.9083});
                    \node at (.5,-.07) {$N=17, \alpha=f_7; |T'(\alpha+1)|=2.0098\cdots$};
                       \node at (.5,-.25) {${}$};    
    \end{tikzpicture}
 $$
\endminipage
\minipage{0.5\textwidth}
$$
\begin{tikzpicture}[scale =5] 
      \draw (.192,0) -- (.192,1);
       \draw (.4435,0) -- (.4435,1);
        \draw (.788,0) -- (.788,1);
                     \draw (0,0) -- (1,0);
             \draw (0,0) -- (1,1);
     \draw (0,0) -- (0,1);
     \draw (0,1) -- (1,1);
     \draw (1,1) -- (1,0);
        \node at (.1,.95) {$\Box_6$};
              \node at (.32,.95) {$\Box_5$};
    \node at (.61,.95) {$\Box_4$};
              \node at (.89,.95) {$\Box_3$};
               \draw[domain=.788:1,smooth,variable=\x] plot ({\x},{12/(\x+1.4244)-4.4244});
             \draw[domain=.4435:.788,smooth,variable=\x] plot ({\x},{12/(\x+1.4244)-5.4244});
                  \draw[domain=.192:.4435,smooth,variable=\x] plot ({\x},{12/(\x+1.4244)-6.4244});
             \draw[domain=0:.192,smooth,variable=\x] plot ({\x},{12/(\x+1.4244)-7.4244});
                    \node at (.5,-.07) {$N=12, \alpha=f_7; |T'(\alpha+1)|=2.0415\cdots$};
                      \node at (.5,-.25) {${}$};   
    \end{tikzpicture}
 $$
\endminipage
\vspace*{-8mm} \caption{Two arrangements illustrating the gaplessness of arrangements with five cylinders or more on account of Theorem \ref{Wilkinson} \label{fig: 17 and 12}}
\end{figure}

Proof of Theorem \ref{no gaps five cylinder case}, part I: Let $\I$ consist of five cylinder sets or more. Then $b(N,\alpha)>3$, implying 
$$
\alpha<\frac12\sqrt{\tfrac{4N}3+1}-\frac12\,\,{\text{(cf.~(\ref{formula a in b})), in which case}}\,\,|\T'(\alpha+1)|>3-\frac{3\sqrt{12N+9}-9}{2N}.
$$
The second inequality yields that for $N \in \N_{\geq 18}$ we have $|\T'(\alpha+1)|>2$ and, applying Theorem \ref{Wilkinson}, $\I$ is gapless. Now suppose $N\in \{12,\ldots,17\}$. Then $b(N\alpha)=3$ involves arrangements with four cylinders. In each of these cases, the smallest $\alpha$ such that $\I$ has not yet (i.e.~decreasing from $\alpha_{\max}$) consisted of five cylinder sets is $f_7$. In all six cases (two of which are illustrated in Figure \ref{fig: 17 and 12}) we have  $|\T'(f_7+1)|>2$, yielding the gaplessness of $\I$ for arrangements with five or more cylinders in case $N\in \{12,\ldots,17\}$. This finishes the proof of Theorem \ref{no gaps five cylinder case} for $N\in \N_{\geq 12}$. For $N\in \{2,\ldots,11\}$ a similar approach does not work. We will use some ideas that we will introduce and develop in the next sections and will finish the proof of Theorem \ref{no gaps five cylinder case} at the end of Section \ref{four and five}.

In the following we will go into conditions for gaplessness of arrangements consisting of less than five cylinders. We will start with two cylinders and will use the results for arrangements with three and four cylinders. \smallskip

\begin{Remark}
Since $b(N,\alpha)$ is a strictly decreasing function of $\alpha$ (cf. Lemma \ref{branch equation}) and $b(N,\alpha_{\max})=1+1/(\sqrt{N}-1)$,  the condition $I_{\alpha}=\Delta_{d}\cup\Delta_{d-1}$ is never satisfied in case $N \in \{2,3\}$.
\end{Remark} 

\section{Gaplessness when $\I$ consists of two cylinder sets}\label{two-cylinder cases}

In general, when the branch number is not much larger than $1$ (which is when $\alpha$ is not much smaller than $\alpha_{\max}$), the overall expanding power of $\T$, determined by $\T'$ (or $|\T'|$, which we will often use), is not enough to exclude the existence of gaps; we shall elaborate on this in a subsequent article. However, in the case of two cylinder sets $\I=\Delta_d\cup\Delta_{d-1}$, there is a very clear condition under which this power suffices: 

\begin{Theorem}\label{gaplessness two cylinders} 
Let $I_{\alpha}=\Delta_{d}\cup\Delta_{d-1}$. If $\T(\alpha)\geq f_{d-1}$ and $\T(\alpha+1)\leq f_{d}$, then $I_{\alpha}$ is gapless.
\end{Theorem}

Although the statement of \ref{gaplessness two cylinders} is intuitively clear, for the proof of Theorem \ref{gaplessness two cylinders} we need several results and lemmas that we will prove first. Then, immediately following Remark \ref{value root} on page \pageref{value root}, we will prove Theorem \ref{gaplessness two cylinders} itself.\smallskip

\begin{Remark}
If either $\T(\alpha)<f_{d-1}$ or $\T(\alpha+1) >f_d$, it is easy to see that $(\T(\alpha),\T^2(\alpha))$ or $(\T^2(\alpha+1),\T(\alpha+1))$ is a gap, respectively.
\end{Remark}

Since arrangements under the condition of Theorem \ref{gaplessness two cylinders} play an important role in this section, we introduce the following notations:

\begin{Definition}\label{family arrangements 1}
Let $N \in \N_{\geq 4}$ be fixed. For $d \in \N_{\geq 2}$, we define $\mathcal F(d)$ as the family of all arrangements $\Upsilon_{N,\alpha}$ such that $\I=\Delta_d\cup\Delta_{d-1}$, $\T(\alpha)\geq f_{d-1}$ and $\T(\alpha+1) \leq f_{d}$. We will write $\F^*(d)$ in case $\alpha$ satisfies the equation $\T(\alpha)= f_{d-1}$, the root of which we will henceforth denote by $\alpha(N,d)$.
\end{Definition}

\begin{Remark}
Note that for each $N\in \N_{\geq 4}$ and $d\in \N_{\geq 2}$ we have that $\alpha(N,d)$, if it exists, is the only value of $\alpha$ such that $\F^*(d)$ is not void. 
\end{Remark}

If the expanding power of $\T$ is large enough to exclude the existence of gaps for \emph{the largest} $\alpha$ for which an arrangement in $\F(d)$ exists, there will not be gaps in any arrangement in $\F(d)$. We will now first show how to find these largest $\alpha$, which takes some effort. When we have finished that, we will go into the expanding power of $|\T'|$ in these arrangements with largest $\alpha$.\smallskip

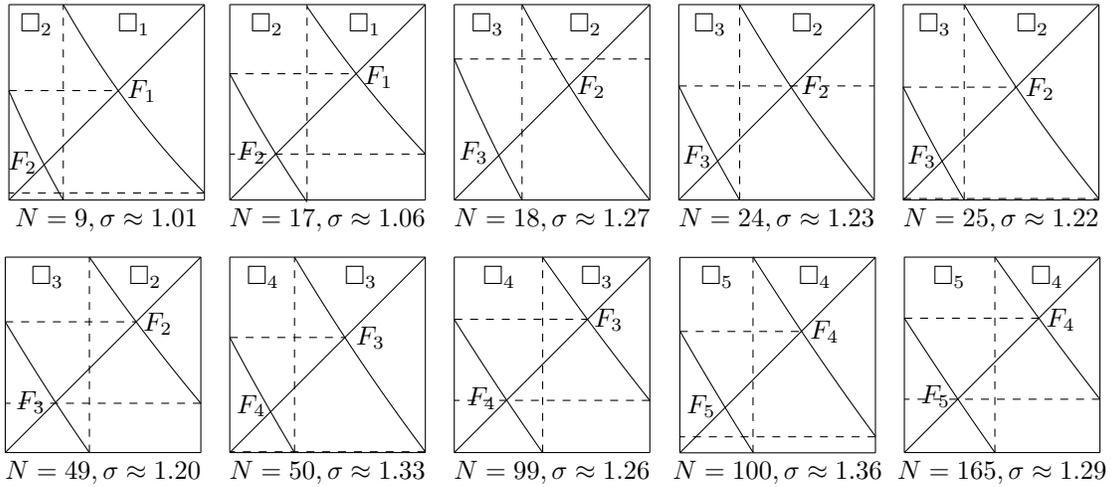
\begin{figure}[!htb]
\minipage{0.2\textwidth}
$$
\begin{tikzpicture}[scale =2.6] 
    \draw (0,0) -- (1,0);
    \draw [dashed] (.278,0) -- (.278,1);
      \draw [dashed] (0,.56) -- (.56,.56);
        \draw [dashed] (1,.036) -- (0,.036);
     \draw (0,0) -- (1,1);
     \draw (0,0) -- (0,1);
     \draw (0,1) -- (1,1);
     \draw (1,1) -- (1,0);
        \draw[domain=0.278:1,smooth,variable=\x] plot ({\x},{9/(\x+1.982)-2.982});
            \draw[domain=0:.278,smooth,variable=\x] plot ({\x},{9/(\x+1.982)-3.982});
         \node at (.64,.9) {$\Box_1$};
      \node at (.14,.9) {$\Box_2$};
         \node at (.68,.56) {$F_1$};
      \node at (.07,.19) {$F_2$};
      \node at (.5,-.1) {$N=9, \sigma\approx1.01$};   
    \end{tikzpicture}
 $$
\endminipage
\minipage{0.2\textwidth}
$$
\begin{tikzpicture}[scale =2.6] 
    \draw (0,0) -- (1,0);
    \draw [dashed] (.395,0) -- (.395,1);
        \draw [dashed] (0,.6462) -- (.6462,.6462);
   \draw [dashed] (1,.235) -- (0,.235);
     \draw (0,0) -- (1,1);
     \draw (0,0) -- (0,1);
     \draw (0,1) -- (1,1);
     \draw (1,1) -- (1,0);
        \draw[domain=0.395:1,smooth,variable=\x] plot ({\x},{17/(\x+3.0071)-4.0071});
            \draw[domain=0:.395,smooth,variable=\x] plot ({\x},{17/(\x+3.0071)-5.0071});
         \node at (.69,.9) {$\Box_1$};
      \node at (.19,.9) {$\Box_2$};
        \node at (.76,.63) {$F_1$};
      \node at (.11,.24) {$F_2$};
      \node at (.5,-.1) {$N=17, \sigma\approx1.06$};   
    \end{tikzpicture}
 $$
\endminipage
\minipage{0.2\textwidth}
$$
\begin{tikzpicture}[scale =2.6] 
    \draw (0,0) -- (1,0);
      \draw [dashed] (0,.7215) -- (1,.7215);
    \draw [dashed] (.3465,0) -- (.3465,1);
     \draw (0,0) -- (1,1);
     \draw (0,0) -- (0,1);
     \draw (0,1) -- (1,1);
     \draw (1,1) -- (1,0);
        \draw[domain=0.3465:1,smooth,variable=\x] plot ({\x},{18/(\x+2.772)-4.772});
            \draw[domain=0:.3465,smooth,variable=\x] plot ({\x},{18/(\x+2.772)-5.772});
         \node at (.67,.9) {$\Box_2$};
      \node at (.17,.9) {$\Box_3$};
         \node at (.7,.58) {$F_2$};
      \node at (.1,.24) {$F_3$};
      \node at (.5,-.1) {$N=18, \sigma\approx1.27$};   
    \end{tikzpicture}
 $$
\endminipage
\minipage{0.2\textwidth}
$$
\begin{tikzpicture}[scale =2.6] 
    \draw (0,0) -- (1,0);
    \draw [dashed] (.3114,0) -- (.3114,1);
     \draw [dashed] (0,.5836) -- (1,.5836);
     \draw (0,0) -- (1,1);
     \draw (0,0) -- (0,1);
     \draw (0,1) -- (1,1);
     \draw (1,1) -- (1,0);
        \draw[domain=.3114:1,smooth,variable=\x] plot ({\x},{24/(\x+3.4244)-5.4244});
            \draw[domain=0:.3114,smooth,variable=\x] plot ({\x},{24/(\x+3.4244)-6.4244});
         \node at (.66,.9) {$\Box_2$};
      \node at (.16,.9) {$\Box_3$};
        \node at (.7,.58) {$F_2$};
      \node at (.095,.22) {$F_3$};
      \node at (.5,-.1) {$N=24, \sigma\approx1.23$};   
    \end{tikzpicture}
 $$
\endminipage
\minipage{0.2\textwidth}
$$
\begin{tikzpicture}[scale =2.6] 
    \draw (0,0) -- (1,0);
    \draw [dashed] (.3118,0) -- (.3118,1);
         \draw [dashed] (0,.5774) -- (.5774,.5774);
     \draw [dashed] (0,.0074) -- (1,.0074);
     \draw (0,0) -- (1,1);
     \draw (0,0) -- (0,1);
     \draw (0,1) -- (1,1);
     \draw (1,1) -- (1,0);
        \draw[domain=.3118:1,smooth,variable=\x] plot ({\x},{25/(\x+3.5216)-5.5216});
            \draw[domain=0:.3118,smooth,variable=\x] plot ({\x},{25/(\x+3.5216)-6.5216});
         \node at (.66,.9) {$\Box_2$};
      \node at (.16,.9) {$\Box_3$};
       \node at (.7,.57) {$F_2$};
      \node at (.095,.22) {$F_3$};
      \node at (.5,-.1) {$N=25, \sigma\approx1.22$};   
    \end{tikzpicture}
  $$
\endminipage
\\
\minipage{0.2\textwidth}
$$
\begin{tikzpicture}[scale =2.6] 
    \draw (0,0) -- (1,0);
       \draw [dashed] (.43,0) -- (.43,1);
        \draw [dashed] (1,.2523) -- (0,.2523);
          \draw [dashed] (0,.669) -- (.669,.669);
     \draw (0,0) -- (1,1);
     \draw (0,0) -- (0,1);
     \draw (0,1) -- (1,1);
     \draw (1,1) -- (1,0);
        \draw[domain=0.43:1,smooth,variable=\x] plot ({\x},{49/(\x+5.402)-7.402});
            \draw[domain=0:.43,smooth,variable=\x] plot ({\x},{49/(\x+5.402)-8.402});
         \node at (.715,.9) {$\Box_2$};
      \node at (.215,.9) {$\Box_3$};
       \node at (.78,.66) {$F_2$};
      \node at (.13,.26) {$F_3$};
      \node at (.5,-.1) {$N=49, \sigma\approx1.20$};   
    \end{tikzpicture}
 $$
\endminipage
\minipage{0.2\textwidth}
$$
\begin{tikzpicture}[scale =2.6] 
    \draw (0,0) -- (1,0);
    \draw [dashed] (.331,0) -- (.331,1);
        \draw [dashed] (1,.0043) -- (0,.0043);
          \draw [dashed] (0,.589) -- (.589,.589);
     \draw (0,0) -- (1,1);
     \draw (0,0) -- (0,1);
     \draw (0,1) -- (1,1);
     \draw (1,1) -- (1,0);
        \draw[domain=0.331:1,smooth,variable=\x] plot ({\x},{50/(\x+5.1396)-8.1396});
            \draw[domain=0:.331,smooth,variable=\x] plot ({\x},{50/(\x+5.1396)-9.1396});
         \node at (.67,.9) {$\Box_3$};
      \node at (.17,.9) {$\Box_4$};
       \node at (.72,.59) {$F_3$};
      \node at (.11,.24) {$F_4$};
      \node at (.5,-.1) {$N=50, \sigma\approx1.33$};   
    \end{tikzpicture}
 $$
\endminipage
\minipage{0.2\textwidth}
$$
\begin{tikzpicture}[scale =2.6] 
    \draw (0,0) -- (1,0);
    \draw [dashed] (.452,0) -- (.452,1);
        \draw [dashed] (1,.267) -- (0,.267);
          \draw [dashed] (0,.682) -- (.682,.682);
     \draw (0,0) -- (1,1);
     \draw (0,0) -- (0,1);
     \draw (0,1) -- (1,1);
     \draw (1,1) -- (1,0);
        \draw[domain=0.452:1,smooth,variable=\x] plot ({\x},{99/(\x+7.881)-10.881});
            \draw[domain=0:.452,smooth,variable=\x] plot ({\x},{99/(\x+7.881)-11.881});
         \node at (.73,.9) {$\Box_3$};
      \node at (.23,.9) {$\Box_4$};
        \node at (.79,.68) {$F_3$};
      \node at (.14,.28) {$F_4$};
      \node at (.5,-.1) {$N=99, \sigma\approx1.26$};   
    \end{tikzpicture}
 $$
\endminipage
\minipage{0.2\textwidth}
$$
\begin{tikzpicture}[scale =2.6] 
    \draw (0,0) -- (1,0);
    \draw [dashed] (.374,0) -- (.374,1);
        \draw [dashed] (1,.082) -- (0,.082);
          \draw [dashed] (0,.621) -- (.621,.621);
     \draw (0,0) -- (1,1);
     \draw (0,0) -- (0,1);
     \draw (0,1) -- (1,1);
     \draw (1,1) -- (1,0);
        \draw[domain=0.374:1,smooth,variable=\x] plot ({\x},{100/(\x+7.576)-11.576});
            \draw[domain=0:.374,smooth,variable=\x] plot ({\x},{100/(\x+7.576)-12.576});
         \node at (.69,.9) {$\Box_4$};
      \node at (.19,.9) {$\Box_5$};
         \node at (.74,.62) {$F_4$};
      \node at (.11,.24) {$F_5$};
      \node at (.5,-.1) {$N=100, \sigma\approx1.36$};   
    \end{tikzpicture}
 $$
\endminipage
\minipage{0.2\textwidth}
$$
\begin{tikzpicture}[scale =2.6] 
    \draw (0,0) -- (1,0);
    \draw [dashed] (.463,0) -- (.463,1);
        \draw [dashed] (1,.273) -- (0,.273);
          \draw [dashed] (0,.6875) -- (.6875,.6875);
     \draw (0,0) -- (1,1);
     \draw (0,0) -- (0,1);
     \draw (0,1) -- (1,1);
     \draw (1,1) -- (1,0);
        \draw[domain=0.463:1,smooth,variable=\x] plot ({\x},{165/(\x+10.3125)-14.3125});
 \draw[domain=0:.463,smooth,variable=\x] plot ({\x},{165/(\x+10.3125)-15.3125});
         \node at (.73,.9) {$\Box_4$};
      \node at (.23,.9) {$\Box_5$};
            \node at (.8,.68) {$F_4$};
      \node at (.16,.29) {$F_5$};
      \node at (.5,-.1) {$N=165, \sigma\approx1.29$};   
    \end{tikzpicture}
 $$
\endminipage
\caption{Arrangements in $\mathcal F(d)$, $d \in \{2,3,4,5\}$, where $\alpha$ is maximal} \label{fig: optimal pd 2}
\end{figure}

For $4 \leq N \leq 8$, with $d=2$, we have $T_{\alpha_{\max}}(\alpha_{\max})>f_1$, while $T_{\alpha_{\max}}(\alpha_{\max}+1)=\alpha_{\max}<f_2$. Hence we see $\Upsilon_{N,\alpha_{\max}}\in \mathcal F(2)$ and $\mathcal F(2)\neq \emptyset$. For $N\in \N_{\geq 9}$ we have $T_{\alpha_{\max}}(\alpha_{\max})<f_1$. When $d=2$ we can find $\alpha$ such that $\Upsilon_{N,\alpha}$ in $\F^*(2)$ for each $9 \leq N \leq 17$; see Figure \ref{fig: optimal pd 2}, where ten arrangements in various $\F(d)$ are drawn. Underneath each arrangement we have mentioned an approximation of $\sigma(\alpha):=|\T'(\alpha+1)|$, which we will later return to. This $\sigma$ is important, because it is the expanding power on the rightmost cylinder set that may be too weak to exclude gaps. \smallskip

When $d=2$ and $N\in \N_{\geq 18}$, the condition $\T(\alpha)=f_{d-1}$ yields $\T(\alpha+1)> f_d$, and $d$ has to \emph{increase by $1$} so as to find an arrangement in $\F(3)$. When $d=3$, for $18\leq N\leq 24$ we find that the largest $\alpha$ is $f_{d-2}-1$, in which case $\T(\alpha+1)=\alpha<f_d$ and $\T(\alpha)>f_{d-1}$ (so in this case the arrangement with the largest $\alpha$ is in $\F(3)$ but not in $\F^*(3)$); for $25\leq N\leq 49$, the largest $\alpha$ is such that $\T(\alpha)=f_2$. When $N\in \N_{\geq 50}$, the family $\F(3)$ is empty and $d$ has to increase further; see Figure \ref{fig: optimal pd 2} once more. In the proof of Lemma \ref{each arrangement in F*^} this approach (of exhausting $\F(d)$ for successive values of $N$ and going to $\F(d+1)$ for larger values of $N$) will be formalised into a proof by induction. Due to (\ref{formula d in b}) such an increase is always possible, no matter how large $d$ and $N$ become.

Note that this inductive approach works since for each $d$ only \emph{finitely} many $N$ exist such that there are $\alpha$ with $\Upsilon_{N,\alpha}\in {\mathcal{F}}(d)$. To see why this claim holds, note that for fixed $N$ and $d$, the smallest $\alpha$ for which $d=d(\alpha )=d_{\rm max}$ is $\alpha_d$, given by
$$
\alpha_d=f_{d+1}=\frac{\sqrt{4N+(d+1)^2}-(d+1)}{2};
$$
cf.~(\ref{fixed point}). For this $\alpha$ it is not necessarily so that $I_{\alpha_d}=\Delta_d\cup \Delta_{d-1}$, i.e.\ that $I_{\alpha_d}$ consists of two cylinder sets (e.g.~if $N=2$ and $d=5$, there are five cylinder sets). However, if $b(N,\alpha_d)\leq 2$ we know that $I_{\alpha_d}$ exists of two cylinder sets, the left one of which is full. According to Lemma~\ref{branch equation}, the branch number $b(N,\alpha_d)$ satisfies
$$
b(N,\alpha_d) = \frac{N}{f_{d+1}(f_{d+1}+1)} = \frac{4N}{4N+(d+1)^2-2d\sqrt{4N+(d+1)^2}+d^2-1}.
$$
Keeping $d$ fixed and letting $N\to\infty$, we find
$$
\lim_{N\to\infty} b(N,\alpha_d) = 1.
$$
In view of this and Lemma \ref{distance between two fixed points for fixed d}  (and the results mentioned directly thereafter), we choose $N$ sufficiently large, such that for $\alpha \geq \alpha_d$ we have $b(N,\alpha )< 5/4$ and $f_{d-1}-f_d > 1/4$. \smallskip

Now suppose that for such a sufficiently large value of $N$ there exists an $\alpha \geq \alpha_d$, such that $\alpha\in {\mathcal{F}}(d)$. Then by Definition \ref{branch number} of branch number and the assumption that $\alpha\in {\mathcal{F}}(d)$, we have that
\begin{equation*}
b(N,\alpha ) \geq  1+ f_{d-1} - f_d  >  1\tfrac{1}{4},
\end{equation*}
which is \emph{impossible} since for $N$ sufficiently large, $d$ fixed and $\alpha\geq \alpha_d$ we have
$$
b(N,\alpha ) \leq b(N,\alpha_d) < 1\tfrac{1}{4}.
$$
It follows that for $d$ fixed and $N$ sufficiently large, ${\mathcal{F}}(d)$ is void.\smallskip

We will prove (in Lemma \ref{each arrangement in F*^}) that when $N\in \N_{\geq 25}$ there exists a minimal $d\in \N_{\geq 3}$ such that the arrangement in $\F(d)$ with $\alpha$ maximal lies in $\F^*(d)$. Before we will prove this, we will explain the relation between $d$ and $N$ for arrangements in $\F^*(d)$.\smallskip

In Figure \ref{fig: optimal pd 2} we see that for $N\in \{49,99,165\}$ the arrangements in $\F^*$ are very similar, and that the arrangement for $N=100$ is more similar to these than the arrangement for $N=50$. Moreover, the last three arrangements look hardly curved. This is easy to understand, considering the following equations, where $b(N,\alpha)=b$ is fixed:
$$
|\T'(\alpha)|=\frac{b(\sqrt{4bN+b^2}+2N+b)}{2N}\quad{\text{and}}\quad |\T'(\alpha)|-|\T'(\alpha+1)| =\frac{b\sqrt{4bN+b^2}}{N}.
$$

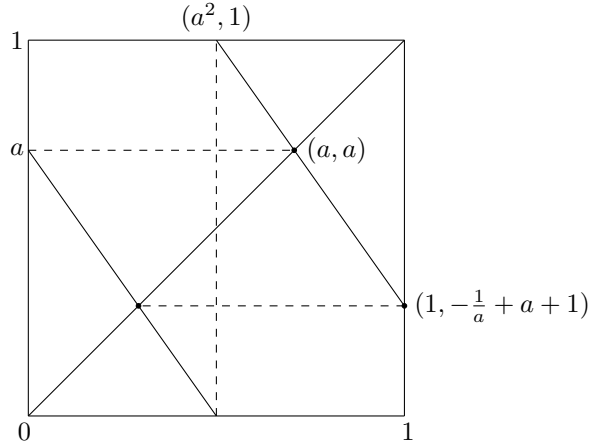
\begin{figure}[!htb]
$$
\begin{tikzpicture}[scale =5] 
\draw[black,fill=black] (.293,.293) circle (.04ex);
\draw[black,fill=black] (.707,.707) circle (.04ex);
\draw[black,fill=black] (1,.293) circle (.04ex);
 \draw [dashed] (.5,0) -- (.5,1);
 \draw [dashed] (0,.707) -- (.707,.707);
  \draw [dashed] (1,.293) -- (.293,.293);
         \draw (0,0) -- (1,1);
     \draw (0,0) -- (0,1);
     \draw (0,1) -- (1,1);
     \draw (1,1) -- (1,0);
        \draw (0,0) -- (1,0);
    \node at (.82,.707) {$(a,a)$};
      \node at (1.26,.293) {$(1,-\tfrac1a+a+1)$};
  \draw[domain=.5:1,smooth,variable=\x] plot ({\x},{1.707-1.4142*\x});
    \draw[domain=0:.5,smooth,variable=\x] plot ({\x},{.707-1.4142*\x});
   \node at (-.03,.71) {$a$};    
      \node at (.5,1.06) {$(a^2,1)$};   
     \node at (-.01,-.04) {$0$};   
         \node at (1.01,-.04) {$1$};  
           \node at (-.03,1) {$1$};  
            \end{tikzpicture}
  $$
\caption{The `limit graph' of $\T$, translated over $(-\alpha,-\alpha)$, under the conditions $\I=\Delta_d\cup\Delta_{d-1}$ and $N/\alpha-d=f_{d-1}$ for $N\to \infty$ (and $\alpha, d \to \infty$ accordingly)} {\label{fig: limit arrangement}}
\end{figure}

Since $(b\sqrt{4bN+b^2})/N$ is a decreasing function of $N$, approaching $0$ from above as $N \to \infty$, the second equation implies that for a fixed branch number $b$ the branches become less curved as $N$ increases; i.e., the curves approach linearity as $N \to \infty$ and $b$ is fixed. Although in $\mathcal F^*(d)$ the branch number is not so much fixed as bounded between $1$ and $2$, we have a similar decrease of curviness as $N$ increases. The arrangements for $N \in \{49,99,165\}$ in Figure \ref{fig: optimal pd 2} suggest that (assuming $\T(\alpha)=f_{d-1}$, i.e.~$\alpha = \alpha(N,d)$)) when $N \to \infty$ (and $d \to \infty$ and $\alpha \to \infty$ accordingly), the difference $f_d-\T(\alpha+1)$ tends to $0$, yielding a `limit graph' of $\T$ that consists of two parallel line segments (the straightened branch curves of $\T$); see Figure \ref{fig: limit arrangement}, obtained by translating the graph over $(-\alpha,-\alpha)$. In this situation we have both $a:=\T(\alpha)\Mod \alpha =f_{d-1}\Mod \alpha$ and $\T(\alpha+1)=f_d$ (also $\Mod \alpha$ in Figure \ref{fig: limit arrangement}). Because in the limit both parts of the graph are linear with the same slope, we also have that $(0,a+1)$ lies on the prolonged right line segment, from which we derive that the line segments have slope $-1/a$. The line with equation $y=-x/a+a+1$ intersects the line $y=1$ at $(a^2,1)$ (so the dividing line is $x=a^2$) and intersects the line $x=1$ in $(1,-1/a+a+1)$, yielding the point $(-1/a+a+1,-1/a+a+1)$ on the line through $(0,a)$ with equation $y=-x/a+a$ (since $\T(\alpha+1)=f_d)$. From this we derive $2a^2=1$, so $a=\sqrt{1/2}$. 

From Figure \ref{fig: limit arrangement} we almost immediately find that the branch number for the limit case is $\sqrt{2}$ and that the dividing line is at $1/2$. We use this \emph{heuristic} to find a formula describing the relation between $N$ and $d$ for arrangements in $\F^*(d)$ very precisely. Note that for arrangements similar to the limit graph we have
$$
1 + f_{d-1} - f_d =\frac{\sqrt{4N+(d-1)^2}-\sqrt{4N+d^2}+1}2+1\approx b(N,\alpha)\approx \sqrt{2},
$$
from which we derive
\begin{equation}\label{d and N limit}
N\approx(4+3\sqrt{2})(d^2-d)+2\quad{\text{or}}\quad d\approx\frac12\left(1+\sqrt{(6\sqrt{2}-8)(N-2)+1}\right).
\end{equation}

If, for $d$ fixed, we determine arrangements in $\F^*(d)$ such that the difference $f_d-\T(\alpha+1)$ is positive and as small as possible according to our heuristic, the best function seems to be $N=(4+3\sqrt{2})(d^2-d)$, yielding the right $N$ (after rounding off to the nearest integer) for $d \in \{3,\ldots,500\} \setminus \{9,50,52,68,69,80,97,129,$ $167,185,210,231,289,330,416,440,444,479,485\}$, in all of which cases the rounding off should have been up instead of down. For $d=2$ we find $N=\lceil 2(4+3\sqrt{2})\rceil=17$, for $d=3$ we find $N=\lfloor 6(4+3\sqrt{2})\rfloor=49$, for $d=4$ we find $N=\lceil 12(4+3\sqrt{2})\rceil=99$ and or $d=5$ we find $N=\lceil 20(4+3\sqrt{2})\rceil=165$; see Figure \ref{fig: optimal pd 2} once more.\smallskip

Although we do not know generally when rounding off to the nearest integer yields the right $N$, with (\ref{d and N limit}) we can find a very good overall indication of the relation between $d$ and $N$ for arrangements in $\F^*(d)$ by looking at the difference between the image of $\alpha(N,d)+1$ and $f_d$; see Definition \ref{family arrangements 1}. With some straightforward calculations we find that 
\begin{equation}\label{a(N,d) condition}
\alpha(N,d)=\frac{N\left(\sqrt{4N+(d-1)^2}-(d+1)\right)}{2(N-d)}.
\end{equation}
Applying (\ref{a(N,d) condition}), we write $f_d(N)-(N/(\alpha(N,d)+1)-(d-1))$ as  
\begin{equation*}
j_d(N):= \frac{(N^2+dN+d)\sqrt{4N+d^2}-N^2\sqrt{4N+(d-1)^2}-(N^2-d(d-4)N-d(d-2))}{2(N^2+dN+d)}
\end{equation*}
and, more generally, define 
\begin{equation}\label{j-d}
j_d(x):= \frac{(x^2+dx+d)\sqrt{4x+d^2}-x^2\sqrt{4x+(d-1)^2}-(x^2-d(d-4)x-d(d-2))}{2(x^2+dx+d)}
\end{equation}
for $x \in [25,\infty)$.

We note that $\I=\Delta_d\cup\Delta_{d-1}$ is equivalent to $N/(\alpha+1)-(d-1)\geq \alpha$. In case $\alpha=\alpha(N,d)$, we have 
\begin{equation}\label{distance t(a+1) to a}
\frac{N}{\alpha(N,d)+1}-(d-1)=\frac {N^2\sqrt{4N+(d-1)^2}-\left((d-1)N^2+2d(d-2)N+2d(d-1)\right )}{2(N^2+dN+d)}.
\end{equation}
Applying (\ref{distance t(a+1) to a}), for the difference $h_d(N) :=  N/(\alpha(N,d)+1) - (d-1) - \alpha(N,d)$ we write
\begin{equation}\label{large expression 1a}
h_d(N):=\frac{2N^3+4dN^2+(2d^3-5d^2+3d)N+2d^2(d-1)-dN(2N+1)\sqrt{4N+(d-1)^2}}{2(N-d)(N^2+dN+d)},
\end{equation}
and, more generally, define 
\begin{equation}\label{large expression 1b}
h_d(x):=\frac{2x^3+4dx^2+(2d^3-5d^2+3d)x+2d^2(d-1)-dx(2x+1)\sqrt{4x+(d-1)^2}}{2(x-d)(x^2+dx+d)},
\end{equation}
for $x \in [25,\infty)$.

Now we can prove the lemma that is illustrated by the arrangements for $N\in \{17,49,99,165\}$ in Figure \ref{fig: optimal pd 2}. In order to so, we define for fixed $d \in \N_{\geq 2}$
$$
S(d):=\{N\in \N_{\geq 4}:\I=\Delta_d\cup\Delta_{d-1}, \T(\alpha) = f_{d-1}(N)\,{\text{and}}\,\, \T(\alpha+1) \leq f_{d}(N)\}
$$
and
$$
M_d:=\max S(d).
$$

\begin{Lemma}\label{boundary values for d and N}
Let $d \in \N_{\geq 2}$. Then $M_d\in \{\lfloor (4+3\sqrt{2})(d^2-d)\rfloor, \lceil (4+3\sqrt{2})(d^2-d)\rceil\}$. 
\end{Lemma}
Proof of Lemma \ref{boundary values for d and N}: First we note that for $d=2$, we have that $\lceil (4+3\sqrt{2})(d^2-d)\rceil$ equals $17$, which corresponds with what we had already calculated and drawn in Figure  \ref{fig: optimal pd 2}. Now let $d\in \N_{\geq 3}$.  First we have to show that $h_d(M_d)>0$ for $M_d\in\{\lfloor (4+3\sqrt{2})(d^2-d)\rfloor,\lceil (4+3\sqrt{2})(d^2-d)\rceil\}$, for this assures us that $\I=\Delta_d\cup\Delta_{d-1}$. We will leave this to the reader; it is merely very cumbersome to show, while technically straightforward\footnote{We have throughout this paper frequently used {\it{(Wolfram) Mathematica}} for making intricate calculations, all of which are nonetheless algebraically basic. In relevant cases we think it will be evident if we did.}. \smallskip

The only thing left to do is showing that 
\begin{equation}\label{two conditions a}
\begin{cases}
j_d((4+3\sqrt{2})(d^2-d))>0;\\
j_d((4+3\sqrt{2})(d^2-d)+1)<0,
\end{cases}
\end{equation}
since the first equation implies that $j_d(\lfloor(4+3\sqrt{2})(d^2-d)\rfloor)>0$, while the second implies that $j_d(\lceil(4+3\sqrt{2})(d^2-d)\rceil+1)<0$. The work to be done is as cumbersome and straightforward as the previous work to be done for this proof and is left to the reader as well. \hfill$\Box$\smallskip

Before we will show that for $N\in \N_{\geq 25}$ there are a $d\in \N_{\geq 3}$ and an $\alpha$ such that $\Upsilon_{N,\alpha} \in \F^*(d)$, we will prove the following lemma:

\begin{Lemma}\label{increase of distance to alpha}
Let $d \in \N_{\geq 3}$. Let $N \in \N_{\geq 25}$ be such that $I_{\alpha(M_d,d)}=\Delta_d\cup\Delta_{d-1}$ and $T_{\alpha (M_d,d)}(\alpha (M_d,d))=f_{d-1}$ for $M_d \in \{N,N+1\}$. Then 
$$
T_{\alpha(N+1,d)}(\alpha(N+1,d)+1)-\alpha(N+1,d)>T_{\alpha(N,d)}(\alpha(N,d)+1)-\alpha(N,d),
$$
i.e. $h_d(N+1) > h_d(N)$.
\end{Lemma}

Proof of Lemma \ref{increase of distance to alpha}: We want to show that $h_d(N)$ from (\ref{large expression 1a}) is an increasing sequence, and do so by calculating the derivative of
with $x \in [25,\infty)$, and then showing that $h_d'(x)>0$ on $[25,\infty)$. Although a little bit intricate, the work is straightforward and is left to the reader. \hfill$\Box$\smallskip

Now we can prove the following lemma:

\begin{Lemma}\label{each arrangement in F*^}
Let $N \in \{9,\ldots,17,25,26,\ldots\}$. Then there are $d \in \N_{\geq 2}$ and $\alpha \in (0,\sqrt{N}-1)$ such that $\I=\Delta_d\cup\Delta_{d-1}$, $\T(\alpha) = f_{d-1}$ and $\T(\alpha+1)\leq f_d$ (i.e.\ $\alpha = \alpha (N,d)$).
\end{Lemma}
Proof of Lemma \ref{each arrangement in F*^}: We will use induction on $d$. For $N \in \{9,\ldots,17,25,26,\ldots,99\}$ and $d \in \{2,3,4\}$ we refer to Figure \ref{fig: optimal pd 2} and leave the calculations to the reader. Specifically, we have for $50\leq N \leq 99$ that $\Upsilon_{N,\alpha(N,4)} \in \F^*(4)$. It is easily seen that $\Upsilon_{99,\alpha(99,5)} \in \F^*(5)$ as well. Due to Lemma \ref{boundary values for d and N}, there is an $N_5>99$ such that $\Upsilon_{N_5,\alpha(N_5,5)} \in \F^*(5)$. Applying Lemma (\ref{increase of distance to alpha}), we see that for all $N \in \{99,\ldots,N_5\}$ we have $\Upsilon_{N,\alpha(N,5)} \in \F^*(5)$. For the induction step, let $d\in \N_{\geq 5}$ be such that there is an $\alpha$ for which $\Upsilon_{N_d,\alpha} \in \F^*(d)$, where $N_d$ is the largest such $N$ possible, cf.~Lemma \ref{boundary values for d and N}. If we can show that for this $N_d$ there is an $\alpha'$ such that $\Upsilon_{N_d,\alpha'} \in \F^*(d+1)$, we are finished. This can be done by showing that 
\begin{equation}\label{two conditions b}
h_{d+1}((4+3\sqrt{2})(d^2-d)-1)>0,
\end{equation}
for this implies $h_{d+1}(N_d)>0$, in which case $\alpha'$ is such that $N_d/\alpha'-(d+1)=f_d$, i.e.~$\alpha^{\prime} = \alpha(N,d+1)$. Although intricate, the calculations are straightforward and are left to the reader. \hfill$\Box$

\begin{Remark}
Although Lemma \ref{each arrangement in F*^} is about $N$ in the first place, our approach is actually based on increasing $d$ and then determining all $N$ such that arrangements $\Upsilon_{N,\alpha} \in \F^*(d)$ exist. The proof of Lemma \ref{each arrangement in F*^} yields the arrangements with the \emph{smallest} $d$ (and therefore the largest $\alpha$) for which $\Upsilon_{N,\alpha} \in \F^*(d)$, as illustrated by the last five arrangements of Figure \ref{fig: optimal pd 2}. 
\end{Remark}

\begin{Example}
For $d=4$ we have $M_d-1=\lfloor(4+3\sqrt{2})(d^2-d)\rfloor=98$. Then $\Upsilon_{M_d-1,\alpha(M_d-1,4)} \in \F^*(4)$ and $\Upsilon_{M_d,\alpha(M_d,4)} \in \F^*(4)$, while $\Upsilon_{M_d+1,\alpha(M_d+1,5)} \in \F^*(5)$; see Figure \ref{fig: optimal pd 2}. It follows immediately from our construction of $\alpha(M_d+1,5)$ that this is the largest $\alpha$ such that $\Upsilon_{M_d+1,\alpha} \in \F(5)$.\smallskip
\end{Example}

With manual calculations we can quickly calculate the expanding power of $\T$ in $\alpha+1$ for arrangements in $\F$ and $\F^*$ and $N$ not too large, say $N \in \N_{\leq 49}$, where the smallest values are found where $\alpha$ is as large as possible. The next proposition gives a lower bound for $|\T'(\alpha+1)|$ for such arrangements for most $N$.

\begin{Proposition}\label{min slope}
Let $N \in \{18\}\cup\{50,51,\ldots\}\setminus \{95,\ldots,99\}$ and $\alpha \in (0, \sqrt{N}-1]$ such that $I_{\alpha}=\Delta_{d}\cup\Delta_{d-1}$ for some $d\in\N$, $d\in \N_{\geq 2}$. Furthermore, suppose that $\T(\alpha)\geq f_{d-1}$ and $\T(\alpha+1) \leq f_{d}$. Then $|\T'(\alpha+1)|>\sqrt[3]{2}=1.259921\cdots$.
\end{Proposition}
Proof of Proposition \ref{min slope}: Considering Lemma \ref{each arrangement in F*^}, we can confine ourselves to arrangements in $\F^*$ with $\alpha$ as large as possible. For $\alpha=\alpha(N,d)$ (cf.\ (\ref{a(N,d) condition})) we can write $|\T'(\alpha+1)|=N/(\alpha+1)^2$ as
\begin{equation}\label{expression slope 1}
k_d(N)=\frac{2N^4+(d-1)^2N^3+2d(d-1)N^2+2d^2N+((d-1)N^3+2dN^2)\sqrt{4N+(d-1)^2}}{2(N^4+2dN^3+d(d+2)N^2+2d^2N+d^2)}.
\end{equation}
It is not hard to find that, for $d$ \emph{fixed}, $k_d$ is a decreasing sequence, with $\lim_{N \to \infty} k_d(N)=1$. However, from (\ref{d and N limit}) it follows that if $N\to\infty$ we have that also $d\to\infty$ in a precise manner. Due to the previous lemmas, for each $d$ we can confine ourselves to considering only $N/(\alpha+1)^2$ for the largest $N$ and $\alpha$ such that $\Upsilon_{N,\alpha} \in \F^*(d)$. Applying Lemma \ref{boundary values for d and N}, an easy way to check if indeed $|T_{\alpha(N,d)}'(\alpha(N,d)+1)|>\sqrt[3]{2}$ is considering $k_d(x)$, with $x \in [100,\infty)$, and then calculating $k_d((4+3\sqrt{2})(d^2-d)+1)$ for $d\in \N_{\geq 5}$, which is amply larger than $\sqrt[3]{2}=1.2599\cdots$. For the remaining cases $d=3$ and $N=18$ and for $d=4$ and $N \in \{50,51,\ldots,94\}$ it is easily checked manually that indeed $|T_{\alpha(N,d)}'(\alpha(N,d)+1)|>\sqrt[3]{2}$. \hfill $\Box$

\begin{Remark}
Considering our previous remarks concerning arrangements in $\F^*$, it may be clear that $\lim_{N\to \infty} N/(\alpha(N,d)+1)^2=\sqrt{2}$.
\end{Remark}
\begin{Remark}\label{value root}
The value $\sqrt[3]{2}$ in the proof of Proposition \ref{min slope} relates to the proof of Theorem \ref{gaplessness two cylinders} and also to the proofs of Proposition \ref{gaplessness three cylinders easy} and Theorem \ref{Wilkinson}, where the numbers $\sqrt{2}$ and $2$ have a similar importance. Considering the proof of Proposition \ref{min slope}, we could actually replace $\sqrt[3]{2}$ by the smallest possible value, given by
$$
\frac{94}{(\alpha(94)+1)^2}=\frac{20480015+320305\sqrt{385}}{21233664}=1.2604\cdots.
$$ 
\end{Remark}

Finally we are ready to prove Theorem \ref{gaplessness two cylinders}, stating that $\I=\Delta_{d}\cup\Delta_{d-1}$, with $d:=d(\alpha)$, is gapless if $\T(\alpha)\geq f_{d-1}$ and $\T(\alpha+1)\leq f_{d}$. Considering Remark \ref{value root} the value $\sqrt[3]{2}$ in Proposition \ref{min slope} can be replaced by $1.26$, the third power of which is $2.000376$. We will use this to stress that the gaplessness of Theorem \ref{gaplessness two cylinders} is actually relatively ample and does not require infinitesimal estimations.\smallskip

Proof of Theorem \ref{gaplessness two cylinders}: First we note that the conditions imply $N \in \N_{\geq 4}$. Now let $\Upsilon_{N,\alpha} \in \F(d)$ and let $K\subset \I$ be any open interval. Since $K$ expands under $\T$, there is an $n \in \N\cup\{0\}$ such that $\T^n(K)$ contains for the first time a fixed point or the discontinuity point $p_d$, in the former case of which we are finished. So we assume that $T_{\alpha}^n(K)\cap \Delta_d=(b,p_d]=:L$, with $f_d <b < p_d$. Note that $\T(L)=[\alpha,\T(b))\subset [\alpha ,f_d)$. For $\T^2(L)=(\T^2(b),\T(\alpha)]$,  we similarly may assume that $f_{d-1}<\T^2(b)<\alpha+1$ (since otherwise $f_{d-1}\in \T^2(L)$, and again we are done).\smallskip

Now suppose that $\T^3(L)$ contains $p_d$, excluding $f_d \in \T^3(L)$. Then $\T^3(L) = L_1\cup M_1$, with $L_1 = [\T^2(\alpha ),p_d]$ and $M_1=(p_d,\T^3(b))$. First we confine ourselves to $N \in \{18\}\cup\{50,\ldots\} \setminus \{95,\ldots,99\}$. Since then $|\T^3(L)|>2.000376|L|$ (cf.\ Remark \ref{value root}), we have certainly $|L_1|>1.001|L|$ or $|M_1|>1.001|L|$. If we consider the images of $L_1$ and $M_1$ under $\T$, $\T^2$ and $\T^3$ similarly as we did with the images of $L$, we find that due to expansiveness (see the proof of Theorem \ref{Wilkinson}) there \textit{must} be an $m$ such that $f_d \in \T^{3m}(L_1)$ or $f_{d-1} \in \T^{3m}(M_1)$ and we are finished. If $\T^3(L)$ does \emph{not} contain $p_d$, the expansion of $L$ will only go on longer, yielding even larger $L_1'$ and $M_1'$ and the reasoning would only be stronger that no gaps can exist.

For $N \in \{4,\ldots,17,19,20,\ldots,49,95,96,\ldots,99\}$ a similar approach can be taken, but there is no useful general lower bound for $|\T'(x)|$ on $\I$. For these cases, however, the moderate expanding power in $\Delta_{d-1}$ is easily made up for by a relatively strong expanding power in $\Delta_d$, and the gaplessness is easily, although tediously, checked by hand (cf.\ Examples \ref{example N=7} and \ref{example N=99} below). This finishes the proof of Theorem \ref{gaplessness two cylinders}.\\
${}$ \hfill $\Box$\smallskip

\begin{Example}\label{example N=7}
In case $N=7$, there exist $\alpha \in (0,\sqrt{7}-1]$ for which $\I=\Delta_2\cup\Delta_1$. The largest $\alpha$ for which $\Upsilon_{7,\alpha} \in \F(2)$ is $\alpha_{\max}=\sqrt{7}-1$, in which case $|\T'(\alpha+1)|=1$. However, $|\T'(f_2)|=2.0938\cdots>2$, and the approach taken above works if only for the expanding power of $\T$ on $[\alpha,f_2)$.
\end{Example}
\begin{Example}\label{example N=99}
In case $N=99$, we have $\I=\Delta_4\cup\Delta_3$, and $\Upsilon_{99,\alpha} \in \F^*$ for $\alpha=99(\sqrt{405}-5)/190=7.8807\cdots$. Then $|\T'(\alpha+1)|=1.2552\cdots$, $|\T'(f_3)|=1.3503\cdots$ and $|\T'(f_4)|=1.4908\cdots$. So for an interval $(p_4,x)$, with $x \in (p_d,f_3)$, assuming that $f_4\not\in T_{\alpha}^3(p_4,x)$, we have $|\T^3(p_4,x)|> 1.3503\cdots\times1.2552\cdots\times1.4908\cdots\times|(p_4,x)|\gg2|(p_4,x)|$, implying enough expanding power for $\T^3$ to exclude the existence of gaps.
\end{Example}

\begin{Remark}
We can also prove that $|\T'(x)|>\sqrt{2}$ on $\Delta_d$ for all arrangements under the assumptions of Theorem \ref{gaplessness two cylinders}, but we cannot do without knowledge about the slope on $\Delta_{d-1}$. 
\end{Remark}

Next we will make preparations for formulating a sufficient condition for gaplessness in case $\I$ consists of three cylinder sets. Proving it involves more subtleties on the one hand, but will have a lot of similarities with the two-cylinder set case on the other hand. Once we have finished that, not much work remains to be done for gaplessness in case $\I$ consists of four or five cylinder sets. \smallskip

 \section{A sufficient condition for gaplessness when $\I$ consists of three or four cylinder sets}\label{three or more}\smallskip
 
When $I_{\alpha}=\Delta_d\cup\ldots \cup\Delta_{d-m}$, with $m \in \{2,3\}$, there is a sufficient condition for gaplessness that resembles the condition for gaplessness in case $\I$ consists of two cylinder sets a lot:

\begin{Theorem}\label{gaplessness three cylinders or more}
Let $I_{\alpha}=\Delta_d\cup\ldots \cup\Delta_{d-m}$, with $m\in \{2,3\}$. Then $\I$ is gapless if 
$$
\T(\alpha)\geq f_{d-1}\,\,{\text{\bf{or}}}\,\,\T(\alpha+1)\leq f_{d-m+1}.
$$
\end{Theorem}

We will prove this theorem in parts. In Subsection \ref{three-cylinder cases} we will prove Theorem \ref{gaplessness three cylinders or more} for $m=2$; in Subsection \ref{four} we will extend the result of Subsection \ref{three-cylinder cases} to $m=3$; considering Theorem \ref{no gaps five cylinder case}, extension to larger $m$ is not useful. 

\begin{Remark}
The difference between the `and' of Theorem \ref{gaplessness two cylinders} and the `or' of Theorem \ref{gaplessness three cylinders or more} has to do with the existence, in the latter case, of at least one full cylinder set. 
\end{Remark}

\subsection{Gaplessness when $\I$ consists of three cylinder sets}\label{three-cylinder cases}\smallskip

Since we have $m=2$, the condition $\T(\alpha)\geq f_{d-1}$ can be split in 
\begin{equation}\label{conditions}
\begin{cases}
1.\,\,\T(\alpha+1)\leq f_{d-1} \leq \T(\alpha);\\
2.\,\,f_{d-1}\leq  \T(\alpha) \leq \T(\alpha+1);\\
3.\,\,f_{d-1}\leq \T(\alpha+1) \leq \T(\alpha);
\end{cases}
\end{equation}

of course the condition $\T(\alpha+1)\leq f_{d-1}$ can be split in a similar way. We will prove Theorem \ref{gaplessness three cylinders or more} by proving gaplessness according to this distinction in three cases, associated with Lemma \ref{gaplessness three cylinders easy}, \ref{wild case} and \ref{gaplessness three cylinders hard} respectively. The first of these is not very hard to prove:

\begin{Lemma}\label{gaplessness three cylinders easy} 
Let $I_{\alpha}=\Delta_d\cup\Delta_{d-1}\cup \Delta_{d-2}$. If $\T(\alpha)\geq f_{d-1}$ {\bf{and}} $\T(\alpha+1)\leq f_{d-1}$, then $I_{\alpha}$ is gapless.
\end{Lemma}

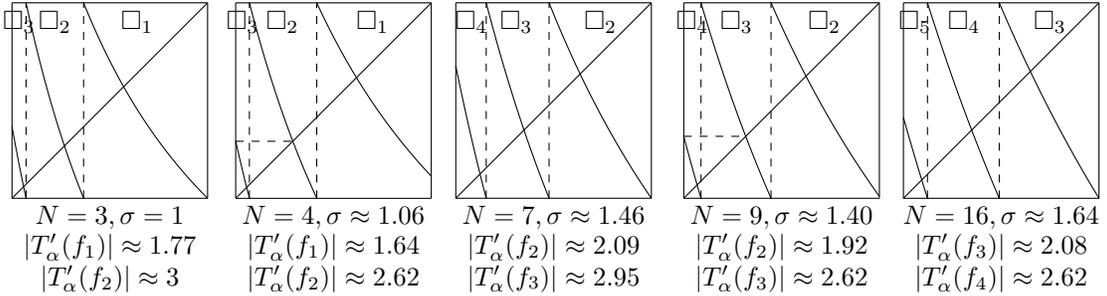
\begin{figure}[!htb]
\minipage{0.2\textwidth}
$$
\begin{tikzpicture}[scale =2.6] 
    \draw (0,0) -- (1,0);
    \draw [dashed] (.072,0) -- (.072,1);
       \draw [dashed] (.366,0) -- (.366,1);
        \draw (0,0) -- (1,1);
     \draw (0,0) -- (0,1);
     \draw (0,1) -- (1,1);
     \draw (1,1) -- (1,0);
        \draw[domain=0.366:1,smooth,variable=\x] plot ({\x},{3/(\x+.732)-1.732});
        \draw[domain=0.072:.366,smooth,variable=\x] plot ({\x},{3/(\x+.732)-2.732});
            \draw[domain=0:.072,smooth,variable=\x] plot ({\x},{3/(\x+.732)-3.732});
         \node at (.64,.9) {$\Box_1$};
      \node at (.22,.9) {$\Box_2$};
          \node at (.036,.9) {$\Box_3$};
              \node at (.5,-.1) {$N=3, \sigma=1$};   
                 \node at (.5,-.25) {$|\T'(f_1)|\approx1.77$};  
                     \node at (.5,-.41) {$|\T'(f_2)|\approx 3$};       
    \end{tikzpicture}
 $$
\endminipage
\minipage{0.2\textwidth}
$$
\begin{tikzpicture}[scale =2.6] 
    \draw (0,0) -- (1,0);
    \draw [dashed] (.07,0) -- (.07,1);
       \draw [dashed] (.4144,0) -- (.4144,1);
           \draw [dashed] (0,.292) -- (.292,.292);
        \draw (0,0) -- (1,1);
     \draw (0,0) -- (0,1);
     \draw (0,1) -- (1,1);
     \draw (1,1) -- (1,0);
        \draw[domain=0.4144:1,smooth,variable=\x] plot ({\x},{4/(\x+.944)-1.944});
        \draw[domain=0.07:.4144,smooth,variable=\x] plot ({\x},{4/(\x+.944)-2.944});
            \draw[domain=0:.07,smooth,variable=\x] plot ({\x},{4/(\x+.944)-3.944});
         \node at (.7,.9) {$\Box_1$};
      \node at (.24,.9) {$\Box_2$};
          \node at (.035,.9) {$\Box_3$};
              \node at (.5,-.1) {$N=4, \sigma\approx1.06$};   
                 \node at (.5,-.25) {$|\T'(f_1)|\approx1.64$};  
                     \node at (.5,-.41) {$|\T'(f_2)|\approx 2.62$};       
    \end{tikzpicture}
 $$
\endminipage
\minipage{0.2\textwidth}
$$
\begin{tikzpicture}[scale =2.6] 
    \draw (0,0) -- (1,0);
          \draw [dashed] (.477,0) -- (.477,1);
             \draw [dashed] (.1555,0) -- (.1555,1);
        \draw (0,0) -- (1,1);
     \draw (0,0) -- (0,1);
     \draw (0,1) -- (1,1);
     \draw (1,1) -- (1,0);
          \draw[domain=0.477:1,smooth,variable=\x] plot ({\x},{7/(\x+1.1926)-3.1926});
        \draw[domain=0.1555:.477,smooth,variable=\x] plot ({\x},{7/(\x+1.1926)-4.1926});
            \draw[domain=0:.1555,smooth,variable=\x] plot ({\x},{7/(\x+1.1926)-5.1926});
           \node at (.74,.9) {$\Box_2$};
      \node at (.31,.9) {$\Box_3$};
        \node at (.08,.9) {$\Box_4$};
           \node at (.5,-.1) {$N=7, \sigma\approx1.46$};   
         \node at (.5,-.25) {$|\T'(f_2)|\approx 2.09$};      
             \node at (.5,-.41) {$|\T'(f_3)|\approx 2.95$};      
    \end{tikzpicture}
  $$
\endminipage
\minipage{0.2\textwidth}
$$
\begin{tikzpicture}[scale =2.6] 
    \draw (0,0) -- (1,0);
          \draw [dashed] (.446,0) -- (.446,1);
             \draw [dashed] (.088,0) -- (.088,1);
                  \draw [dashed] (0,.317) -- (.317,.317);
        \draw (0,0) -- (1,1);
     \draw (0,0) -- (0,1);
     \draw (0,1) -- (1,1);
     \draw (1,1) -- (1,0);
          \draw[domain=0.446:1,smooth,variable=\x] plot ({\x},{9/(\x+1.5374)-3.5374});
        \draw[domain=0.088:.446,smooth,variable=\x] plot ({\x},{9/(\x+1.5374)-4.5374});
            \draw[domain=0:.088,smooth,variable=\x] plot ({\x},{9/(\x+1.5374)-5.5374});
           \node at (.72,.9) {$\Box_2$};
      \node at (.27,.9) {$\Box_3$};
        \node at (.044,.9) {$\Box_4$};
           \node at (.5,-.1) {$N=9, \sigma\approx1.40$};   
         \node at (.5,-.25) {$|\T'(f_2)|\approx 1.92$};      
             \node at (.5,-.41) {$|\T'(f_3)|\approx 2.62$};      
    \end{tikzpicture}
 $$
\endminipage
\minipage{0.2\textwidth}
$$
\begin{tikzpicture}[scale =2.6] 
    \draw (0,0) -- (1,0);
    \draw [dashed] (.123,0) -- (.123,1);
        \draw [dashed] (.49,0) -- (.49,1);
          \draw (0,0) -- (1,1);
     \draw (0,0) -- (0,1);
     \draw (0,1) -- (1,1);
     \draw (1,1) -- (1,0);
        \draw[domain=0.49:1,smooth,variable=\x] plot ({\x},{16/(\x+2.1231)-5.1231});
             \draw[domain=0.123:.49,smooth,variable=\x] plot ({\x},{16/(\x+2.1231)-6.1231});
    \draw[domain=0:.123,smooth,variable=\x] plot ({\x},{16/(\x+2.1231)-7.1231});
         \node at (.06,.9) {$\Box_5$};
      \node at (.31,.9) {$\Box_4$};
       \node at (.75,.9) {$\Box_3$};
      \node at (.5,-.1) {$N=16, \sigma\approx1.64$};   
    \node at (.5,-.25) {$|\T'(f_3)|\approx 2.08$}; 
      \node at (.5,-.41) {$|\T'(f_4)|\approx 2.62$};           
    \end{tikzpicture}
 $$
\endminipage
\caption{Arrangements with largest $\alpha$ such that there is a $d$ with $\I=\Delta_d\cup\Delta_{d-1}\cup\Delta_{d-2}$ under the condition $\T(\alpha)\geq f_{d-1}$ \emph{and} $\T(\alpha+1)\leq f_{d-1}$} \label{fig: optimal pd 3}
\end{figure}

Proof of Lemma \ref{gaplessness three cylinders easy}: The assumptions imply that $b(N,\alpha)>2$, yielding $\sigma(\alpha)=|\T'(\alpha+1)|> \sqrt{2}$ for $N\in \N_{\geq 17}$. If $N\in \N_{\geq 17}$, we let $K\subset \I$ be any open interval. Since $K$ expands under $\T$, there is an $n\in \N\cup\{0\}$) such that $\T^n(K)$ contains a fixed point or a discontinuity point $p_{d-i}$ (with $i\in \{0,1\}$) , in the former case of which we are finished. So we assume that $\T^n(K)\supset L$, where $L=(b,p_{d-i}]$, with $f_{d-i} <b < p_{d-i}$, with $i \in \{0,1\}$. If $\T(L)$ contains a fixed point, we are finished. If $\T(L)$ does not contain a fixed point, then it cannot contain a discontinuity point, and we have that $|\T^2(L)|>2|L|$, implying enough expanding power of $\T$ to ensure gaplessness of at least one cylinder set (which might be non-full). Since both $\T(\alpha)\geq f_{d-1}$ and $\T(\alpha+1)\leq f_{d-1}$, it follows that $\I$ is gapless. For $2\leq N \leq 16$ the slopes on $\I$ may differ considerably: for some $N$, such as $N=7$ and $N=16$ we also have $\sigma>\sqrt{2}$, but when this is not he case, the steepness left of $f_{d-2}$  is amply larger then $\sqrt{2}$; see Figure \ref{fig: optimal pd 3} for some examples where $\alpha$ is as large as possible. This finishes the proof of Lemma \ref{gaplessness three cylinders easy} (cf.~case 1 in (\ref{conditions})). \hfill $\Box$

If $\I=\Delta_d\cup\Delta_{d-1}\cup\Delta_{d-2}$ under the condition $\T(\alpha)\geq f_{d-1}$ and $\T(\alpha+1) > f_{d-1}$ or under the condition $\T(\alpha) < f_{d-1}$ and $\T(\alpha+1) \leq f_{d-1}$, $\I$ is gapless as well, but this is much harder to prove. The following definition will be convenient:

\begin{Definition}\label{small set}
Let $\I=\Delta_d\cup\ldots\cup\Delta_{d-m}$, and $1\leq m\leq d-1$. If $\T(\alpha)\leq f_{d-1}$ or $\T(\alpha+1)\geq f_{d-m+1}$, the cylinder set $\Delta_d$ respectively $\Delta_{d-m}$ is called \emph{small}.
\end{Definition}

Taking a similar approach as in the proof of Theorem \ref{gaplessness two cylinders}, one can show that the map $\T$ has enough expansive power to ensure that for any open interval $K\subset \I$ there exists a non-negative integer $n$ such that $\T^n(K)$ contains a fixed point. If this fixed point is in a non-small or even full cylinder set, we are done (as in the proofs of Theorem \ref{gaplessness two cylinders} and Lemma \ref{gaplessness three cylinders easy}). However, if this fixed point is from the \emph{small} cylinder set, then it only follows that every point of the small cylinder set is in the orbit under $\T$ of some point in $K$. Note this implies that the \emph{small} cylinder set is \emph{gapless}. So we may assume that the small cylinder set is gapless. Let us assume that the left cylinder set is small. We define $L:=\T(\Delta_d)\setminus\Delta_d$. Since $\Delta_d$ is gapless, we have $L=(p_d,T(\alpha)]\subset (p_d,f_{d-1})$, so $\T(L)=[\T^2(\alpha),\alpha+1)$. If $\T^2(\alpha)\leq f_{d-2}$, we are finished, so we assume that $\T^2(\alpha) > f_{d-2}$. We then have $\T^2(L)=(\T(\alpha+1),\T^3(\alpha)]$. If $\T^3(\alpha)\geq f_{d-1}$ we are finished, since then $f_{d-1} \in \T^2(L)$.\smallskip

 The question arises whether it is possible to keep avoiding fixed points if we go on with letting $\T$ work on $L$ and its images (or similarly, when the right cylinder set is small, some interval $R:=\T(\Delta_1)\setminus\Delta_1$). We will argue that this is not possible in the two most plausible cases for gaps to exist, involving the least expansion. 

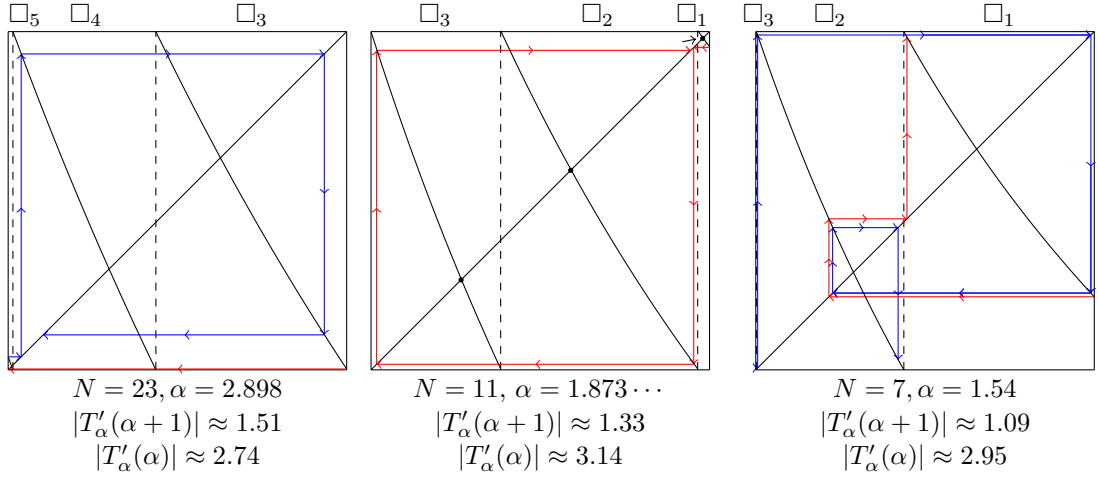
\begin{figure}[!htb]
  \minipage{0.33\textwidth}
$$
\begin{tikzpicture}[scale =4.5]
  \draw (0,0) --   (1,0);
             \draw [dashed] (.4363,0) -- (.4363,1);
                  \draw [dashed] (.014,0) -- (.014,1);
       \draw (0,0) -- (1,1);
     \draw (0,0) -- (0,1);
     \draw (0,1) -- (1,1);
     \draw (1,1) -- (1,0);
         \draw [->,blue](0,.0385) -- (.0385,.0385);
                \draw [->,blue](.0385,.0385) -- (.0385,.48);
                 \draw [->,blue](.0385,.48) -- (.0385,.934);
                            \draw [->,blue](.0385,.934) -- (.48,.934);
        \draw [->,blue](.48,.934) -- (.934,.934);
     \draw [->,blue](.934,.934) -- (.934,.52);
        \draw [->,blue](.934,.52) -- (.934,.1034);
                 \draw [->,blue](.934,.1034) -- (.52,.1034);
                 \draw [->,blue](.52,.1034) -- (.1034,.1034);
                                 \draw [->,red] (1,.00246) -- (.5,.00246);
      \draw [->,red] (.5,.00246) -- (.00246,.00246);
                    \node at (.72,1.05) {$\Box_3$};
                          \node at (.23,1.05) {$\Box_4$};
                            \node at (.05,1.05) {$\Box_5$};
                       \draw[domain=.4363:1,smooth,variable=\x] plot ({\x},{23/(\x+2.898)-5.898});
          \draw[domain=.014:.4363,smooth,variable=\x] plot ({\x},{23/(\x+2.898)-6.898});
                 \draw[domain=0:.014,smooth,variable=\x] plot ({\x},{23/(\x+2.898)-7.898});
                                  \node at (.5,-.06) {$N=23, \alpha=2.898$}; 
                      \node at (.5,-.16) {$|\T'(\alpha+1)|\approx1.51$}; 
                     \node at (.5,-.26) {$|\T'(\alpha)|\approx2.74$};  
    \end{tikzpicture}
  $$
\endminipage
\minipage{0.33\textwidth}
$$
\begin{tikzpicture}[scale =4.5] 
  \draw (0,0) -- (1,0);
         \draw [dashed] (.383,0) -- (.383,1);
        \draw [dashed] (.965,0) -- (.965,1);
                \draw[black,fill=black] (.266,.266) circle (.04ex);
\draw[black,fill=black] (.59,.59) circle (.04ex);
\draw[black,fill=black] (.98,.98) circle (.04ex);
\draw [->] (.92,.97) -- (.96,.978);
             \draw (0,0) -- (1,1);
     \draw (0,0) -- (0,1);
     \draw (0,1) -- (1,1);
     \draw (1,1) -- (1,0);
         \draw [->,red] (1,.9534) -- (.977,.9534);
             \draw [->,red] (.977,.9534) -- (.9534,.9534);
    \draw [->,red] (.9534,.9534) -- (.9534,.485);
    \draw [->,red] (.9534,.485) -- (.9534,.0165);
    \draw [->,red] (.9534,.0165) -- (.485,.0165);
    \draw [->,red] (.485,.0165) -- (.0165,.0165);
     \draw [->,red] (.0165,.0165) -- (.0165,.48);
    \draw [->,red] (.0165,.48) -- (.0165,.9446);
    \draw [->,red] (.0165,.9446) -- (.48,.9446);
    \draw [->,red] (.48,.9446) -- (.9446,.9446);
     \draw[domain=0.965:1,smooth,variable=\x] plot ({\x},{11/(\x+1.873)-2.873});
           \draw[domain=.383:.965,smooth,variable=\x] plot ({\x},{11/(\x+1.873)-3.873});
         \draw[domain=0:.383,smooth,variable=\x] plot ({\x},{11/(\x+1.873)-4.873});
                  \node at (.95,1.05) {$\Box_1$};
      \node at (.67,1.05) {$\Box_2$};
  \node at (.19,1.05) {$\Box_3$};
            \node at (.5,-.06) {$N=11$, $\alpha=1.873\cdots$};  
             \node at (.5,-.16) {$|\T'(\alpha+1)|\approx1.33$}; 
                     \node at (.5,-.26) {$|\T'(\alpha)|\approx3.14$};  
    \end{tikzpicture}
  $$
\endminipage
\minipage{0.33\textwidth}
$$
\begin{tikzpicture}[scale =4.5]
  \draw (0,0) --   (1,0);
             \draw [dashed] (.4374,0) -- (.4374,1);
                  \draw [dashed] (.002,0) -- (.002,1);
       \draw (0,0) -- (1,1);
     \draw (0,0) -- (0,1);
     \draw (0,1) -- (1,1);
     \draw (1,1) -- (1,0);
          \draw [->,blue](0,.0055) -- (.0055,.0055);
        \draw [->,blue](.0055,.0055) -- (.0055,.5);
                \draw [->,blue](.0055,.5) -- (.0055,.99);
        \draw [->,blue](.0055,.99) -- (.5,.99);
              \draw [->,blue](.5,.99) -- (.99,.99);
       \draw [->,blue](.99,.99) -- (.99,.6);
        \draw [->,blue](.99,.6) -- (.99,.227);
                \draw [->,blue](.99,.227) -- (.6,.227);
        \draw [->,blue](.6,.227) -- (.227,.227);
              \draw [->,blue](.5,.99) -- (.99,.99);
                \draw [->,blue](.99,.6) -- (.99,.227);
                \draw [->,blue](.99,.227) -- (.6,.227);
        \draw [->,blue](.6,.227) -- (.227,.227);
              \draw [->,blue](.227,.227) -- (.227,.32);
                \draw [->,blue](.227,.32) -- (.227,.4205);
              \draw [->,blue](.227,.4205) -- (.32,.4205);
                   \draw [->,blue](.32,.4205) -- (.4205,.4205);
                       \draw [->,blue](.4205,.4205) -- (.4205,.22);
                   \draw [->,blue](.4205,.22) -- (.4205,.03);
                          \draw [->,red] (1,.216) -- (.6,.216);
            \draw [->,red] (.6,.216) -- (.216,.216);
     \draw [->,red] (.216,.216) -- (.216,.33);
      \draw [->,red] (.216,.33) -- (.216,.44655);
     \draw [->,red] (.216,.44655) -- (.33,.44655);
      \draw [->,red] (.33,.44655) -- (.44655,.44655);
        \draw [->,red] (.44655,.44655) -- (.44655,.7);
      \draw [->,red] (.44655,.7) -- (.44655,.984);
         \node at (.72,1.05) {$\Box_1$};
                          \node at (.22,1.05) {$\Box_2$};
                            \node at (.01,1.05) {$\Box_3$};
                       \draw[domain=.4374:1,smooth,variable=\x] plot ({\x},{7/(\x+1.54)-2.54});
          \draw[domain=.002:.4374,smooth,variable=\x] plot ({\x},{7/(\x+1.54)-3.54});
                 \draw[domain=0:.002,smooth,variable=\x] plot ({\x},{7/(\x+1.54)-4.54});
                                  \node at (.5,-.06) {$N=7, \alpha=1.54$}; 
                      \node at (.5,-.16) {$|\T'(\alpha+1)|\approx1.09$}; 
                           \node at (.5,-.26) {$|\T'(\alpha)|\approx2.95$};  
    \end{tikzpicture}
 $$
\endminipage
\caption{Arrangements with one very small cylinder} \label{fig: optimal pd 6}
\end{figure}

The first case is illustrated with two arrangements in Figure \ref{fig: optimal pd 6}, one where $N=23$ and one where $N=11$. In both arrangements one outer cylinder is very small while the other one is full or almost full. In the arrangement where $N=23$, we see that $L$ is a very narrow strip between $p_5$ and $\T(\alpha)$, $\T^2(L)$ is not so narrow anymore, and $\T^3(L)$ is definitely wide enough to make clear that avoiding fixed points $f_4$ and $f_3$ is not possible. The middle arrangement, where $N=11$, is an example of the case where $\Delta_{d-2}$ is small and $\Delta_d$ is actually full. Here we have that $R:=\T(\Delta_1)\setminus\Delta_1$ is a very narrow strip between $\T(\alpha+1)$ and $p_2$ and that $\T^2(R)$ is only slightly larger than $\T(\Delta_1)$, whence eventually there will be an $n \in \N$ such that $f_3 \in \T^n(R)$ or $f_2 \in \T^n(R)$.\smallskip

The rightmost  arrangement in Figure \ref{fig: optimal pd 6} is an is an illustration of the second plausible case for the existence of gaps: here $\Delta_3$ is small, while $\Delta_1$ is incomplete but \emph{not} small. This arrangement illustrates the role $p_d$ might play in avoiding fixed points: in this case, taking $L:=\T(\Delta_3)\setminus\Delta_3$, we have $\T^3(L)=M_1\cup M_2$, with $M_1=[\T^4(\alpha),p_2]$ and $M_2=(p_2,\T^2(\alpha+1)]$. Since $\T^3(L)$ contains a discontinuity point, the expansion under $\T$ is interrupted. If $\T(M_1)$ would be a subset of $\T(\Delta_3)$ and $\T(M_2)$ would be a subset of $\T(L)$, the expansion would be finished and we would have three gaps: $(\T(\alpha),\T(\alpha+1))$, $(\T^3(\alpha),\T^4(\alpha))$ and $(\T^2(\alpha+1),\T^2(\alpha))$ -- but this is not the case, as we will shortly prove.  

Of course arrangements exist such that one of the outer cylinders is small, fixed points are avoided (in the sense we used above) for a long time and it takes more of $\T$ working on $L$ or $R$ before one of the discontinuity points is captured. But in these cases the interruption of the expansion is even weaker than in the cases above. We will first show that arrangements such as the rightmost one of Figure \ref{fig: optimal pd 6} exclude the existence of gaps (cf.\ Lemma \ref{wild case}) and will then consider cases such as the first two arrangements of Figure \ref{fig: optimal pd 6} (cf.\ Lemma \ref{gaplessness three cylinders hard}).\smallskip

\begin{Lemma}\label{wild case}
Let $I_{\alpha}=\Delta_d\cup\Delta_{d-1}\cup\Delta_{d-2}$.Then $\I$ is gapless if 
$$
f_{d-1}\leq \T(\alpha) \leq \T(\alpha+1)\,\,{\text{or}}\,\,\T(\alpha)\leq \T(\alpha+1) \leq f_{d-1}.
$$
\end{Lemma}

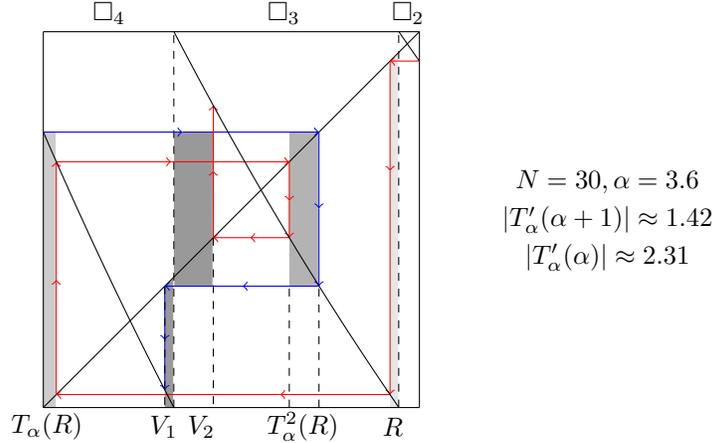
\begin{figure}[!htb]
$$
\begin{tikzpicture}[scale =5]
\draw[draw=white,fill=gray!20!white] 
plot[smooth,samples=100,domain=.922:.945] (\x,0) --
plot[smooth,samples=100,domain=.945:.922] (\x,\x);
\draw[draw=white,fill=gray!40!white] 
plot[smooth,samples=100,domain=0:.0346] (\x,0) --
plot[smooth,samples=100,domain=.0346:0] (\x,.733);
 \draw[draw=white,fill=gray!60!white] 
plot[smooth,samples=100,domain=.654:.733] (\x,.323) --
plot[smooth,samples=100,domain=.733:.654] (\x,.733);
 \draw[draw=white,fill=gray!80!white] 
plot[smooth,samples=100,domain=.347:.452] (\x,.323) --
plot[smooth,samples=100,domain=.452:.347] (\x,.733);
 \draw[draw=white,fill=gray!80!white] 
plot[smooth,samples=100,domain=.323:.347] (\x,0) --
plot[smooth,samples=100,domain=.347:.323] (\x,\x);
  \draw (0,0) --   (1,0);
      \draw [dashed] (.945,0) -- (.945,1);
             \draw [dashed] (.347,0) -- (.347,1);
               \draw [dashed] (.654,0) -- (.654,.323);
             \draw [dashed] (.733,0) -- (.733,.323);
       \draw (0,0) -- (1,1);
     \draw (0,0) -- (0,1);
     \draw (0,1) -- (1,1);
     \draw (1,1) -- (1,0);
          \draw [->,blue](0,.733) -- (.37,.733);
        \draw [->,blue](.37,.733) -- (.733,.733);
        \draw [->,blue](.733,.733) -- (.733,.53);
     \draw [->,blue](.733,.53) -- (.733,.323);
        \draw [->,blue](.733,.323) -- (.53,.323);
                 \draw [->,blue](.53,.323) -- (.323,.323);
                 \draw [->,blue](.323,.323) -- (.323,.18);
                 \draw [->,blue](.323,.18) -- (.323,.047);
                   \draw [->,red] (1,.922) -- (.922,.922);
      \draw [->,red] (.922,.922) -- (.922,.63);
     \draw [->,red] (.922,.63) -- (.922,.0346);
      \draw [->,red] (.922,.0346) -- (.63,.0346);
     \draw [->,red] (.63,.03466) -- (.0346,.0346);
      \draw [->,red] (.0346,.0346) -- (.0346,.34);
          \draw [->,red] (.0346,.34) -- (.0346,.654);
      \draw [->,red] (.0346,.654) -- (.34,.654);
     \draw [->,red] (.34,.654) -- (.654,.654);
      \draw [->,red] (.654,.654) -- (.654,.55);
     \draw [->,red] (.654,.55) -- (.654,.452);
       \draw [->,red] (.654,.452) -- (.55,.452);
     \draw [->,red] (.55,.452) -- (.452,.452);
       \draw [->,red] (.452,.452) -- (.452,.63);
     \draw [->,red] (.452,.63) -- (.452,.804);
     \draw [dashed] (.452,.452) -- (.452,-.01);
       \draw [dashed] (.323,.323) -- (.323,-.01);
                 \node at (.97,1.05) {$\Box_2$};
      \node at (.64,1.05) {$\Box_3$};
       \node at (.175,1.05) {$\Box_4$};
             \draw[domain=0.945:1,smooth,variable=\x] plot ({\x},{30/(\x+3.6)-5.6});
         \draw[domain=.347:.945,smooth,variable=\x] plot ({\x},{30/(\x+3.6)-6.6});
                      \draw[domain=0:.347,smooth,variable=\x] plot ({\x},{30/(\x+3.6)-7.6});
         \node at (1.5,.6) {$N=30, \alpha=3.6$}; 
                \node at (1.5,.5) {$|\T'(\alpha+1)|\approx1.42$}; 
                        \node at (.42,-.05) {$V_2$};  
                               \node at (.32,-.05) {$V_1$}; 
                                      \node at (.93,-.05) {$R$};  
                               \node at (.01,-.05) {$\T(R)$}; 
                                         \node at (.69,-.05) {$\T^2(R)$}; 
                        \node at (1.5,.4) {$|\T'(\alpha)|\approx2.31$};  
           \end{tikzpicture}
 $$
\caption{Arrangement illustrating Lemma \ref{wild case}} \label{fig: optimal pd 7}
\end{figure}

Proof of Lemma \ref{wild case}: We will confine ourselves to the first case of this lemma, that is when $f_{d-1}\leq \T(\alpha) \leq \T(\alpha+1)$, since the second one is proved similarly (in fact, this case is slightly harder due to the smaller size of the absolute value of the derivatives). Regarding our observations above, we may assume that the small cylinder set (which in this case is $\Delta_{d-2}$) is gapless (cf.~the remarks after Definition \ref{small set}). We will show that this implies the gaplessness of the other cylinder sets as well. We define $R:=\T(\Delta_{d-2})\setminus\Delta_{d-2}$ and try to determine $\alpha$ such that $p_d \in \T^3(R)$ (see the remark immediately preceding this lemma). Necessary conditions for this are $\T^2(\alpha)<p_d<\T^4(\alpha+1)$ (assuming that $f_d\not\in T_{\alpha}(R)$ and $f_{d-1}\not\in T_{\alpha}^2(R)$, since in either case we would be done). If these conditions are satisfied, we write $\T^3(R)=V_1\cup V_2$, with $V_1=[\T^2(\alpha),p_d]$ and $V_2=(p_d,\T^4(\alpha+1)]$. We will show that we cannot have both $\T(V_1)\subset \T(R)$ and $\T(V_2)\subset \T(\Delta_{d-2})$, which is necessary for limiting the expansion of $R$ under $\T$ and so not eventually capturing $f_d$ and $f_{d-1}$; see Figure \ref{fig: optimal pd 7}. 

We take an approach that is similar to the proof of Theorem \ref{gaplessness two cylinders}, for which several lemmas and a proposition where used, partially concerning a relation between $N$ and $d$ in the arrangements involved, partially concerning the slope in $\alpha+1$. In this proof we will not explicitly formulate similar statements as lemmas or propositions, nor do we prove them, since they require similar basic but very intricate calculations that we prefer to omit.\smallskip

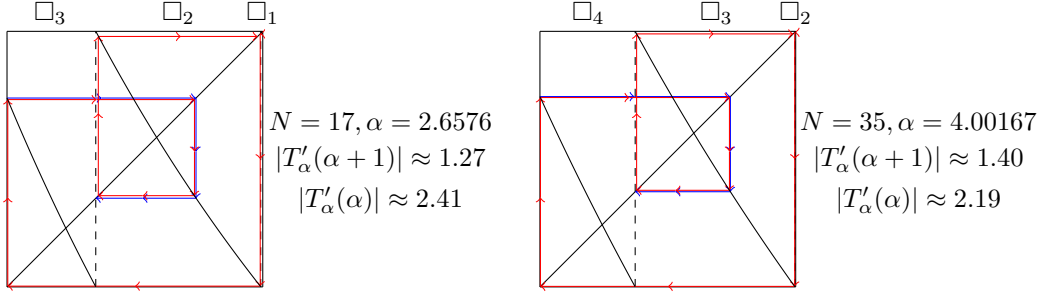
\begin{figure}[!htb]
\minipage{0.48\textwidth}
$$
\begin{tikzpicture}[scale =3.4]
  \draw (0,0) --   (1,0);
        \draw [dashed] (.347,0) -- (.347,1);
        \draw [dashed] (.992,0) -- (.992,1);
       \draw (0,0) -- (1,1);
     \draw (0,0) -- (0,1);
     \draw (0,1) -- (1,1);
     \draw (1,1) -- (1,0);
                 \draw [->,blue](0,.739) -- (.37,.739);
        \draw [->,blue](.37,.739) -- (.739,.739);
                \draw [->,blue](.739,.739) -- (.739,.53);
        \draw [->,blue](.739,.53) -- (.739,.347);
              \draw [->,blue](.739,.347) -- (.53,.347);
        \draw [->,blue](.53,.347) -- (.347,.347);
                   \draw [->,red] (1,.99) -- (.99,.99);
      \draw [->,red] (.99,.99) -- (.99,.5);
     \draw [->,red] (.99,.5) -- (.99,.003);
      \draw [->,red] (.99,.003) -- (.5,.003);
     \draw [->,red] (.5,.003) -- (.003,.003);
       \draw [->,red] (.003,.003) -- (.003,.35);
     \draw [->,red] (.003,.35) -- (.003,.733);
      \draw [->,red] (.003,.733) -- (.35,.733);
     \draw [->,red] (.35,.733) -- (.733,.733);
  \draw [->,red] (.733,.733) -- (.733,.53);
     \draw [->,red] (.733,.53) -- (.733,.357);
      \draw [->,red] (.733,.357) -- (.53,.357);
     \draw [->,red] (.53,.357) -- (.357,.357);
     \draw [->,red] (.357,.357) -- (.357,.68);
     \draw [->,red] (.357,.68) -- (.357,.98);
      \draw [->,red] (.357,.98) -- (.68,.98);
     \draw [->,red] (.68,.98) -- (.98,.98);
              \node at (1,1.07) {$\Box_1$};
                          \node at (.67,1.07) {$\Box_2$};
                            \node at (.17,1.07) {$\Box_3$};
                       \draw[domain=.992:1,smooth,variable=\x] plot ({\x},{17/(\x+2.6576)-3.6576});
          \draw[domain=.347:.992,smooth,variable=\x] plot ({\x},{17/(\x+2.6576)-4.6576});
                  \draw[domain=0:.347,smooth,variable=\x] plot ({\x},{17/(\x+2.6576)-5.6576});
              \node at (1.46,.64) {$N=17, \alpha=2.6576$}; 
                      \node at (1.46,.5) {$|\T'(\alpha+1)|\approx1.27$}; 
                           \node at (1.46,.34) {$|\T'(\alpha)|\approx2.41$};  
    \end{tikzpicture}
 $$
\endminipage
\minipage{0.48\textwidth}
$$
\begin{tikzpicture}[scale =3.4]
  \draw (0,0) --   (1,0);
        \draw [dashed] (.3724,0) -- (.3724,1);
        \draw [dashed] (.997,0) -- (.997,1);
       \draw (0,0) -- (1,1);
     \draw (0,0) -- (0,1);
     \draw (0,1) -- (1,1);
     \draw (1,1) -- (1,0);
                 \draw [->,blue](0,.7447) -- (.37,.7447);
        \draw [->,blue](.37,.7447) -- (.7447,.7447);
                \draw [->,blue](.7447,.7447) -- (.7447,.53);
        \draw [->,blue](.7447,.53) -- (.7447,.3724);
              \draw [->,blue](.7447,.3724) -- (.53,.3724);
        \draw [->,blue](.53,.3724) -- (.3724,.3724);
                   \draw [->,red] (1,.996) -- (.996,.996);
      \draw [->,red] (.996,.996) -- (.996,.5);
     \draw [->,red] (.996,.5) -- (.996,.002);
      \draw [->,red] (.996,.002) -- (.5,.002);
     \draw [->,red] (.5,.002) -- (.002,.002);
       \draw [->,red] (.002,.002) -- (.002,.35);
     \draw [->,red] (.002,.35) -- (.002,.741);
      \draw [->,red] (.002,.741) -- (.35,.741);
     \draw [->,red] (.35,.741) -- (.741,.741);
  \draw [->,red] (.741,.741) -- (.741,.53);
     \draw [->,red] (.741,.53) -- (.741,.378);
      \draw [->,red] (.741,.378) -- (.53,.378);
     \draw [->,red] (.53,.378) -- (.378,.378);
     \draw [->,red] (.378,.378) -- (.378,.68);
     \draw [->,red] (.378,.68) -- (.378,.99);
      \draw [->,red] (.378,.99) -- (.68,.99);
     \draw [->,red] (.68,.99) -- (.99,.99);
              \node at (1,1.07) {$\Box_2$};
                          \node at (.69,1.07) {$\Box_3$};
                            \node at (.19,1.07) {$\Box_4$};
                       \draw[domain=.997:1,smooth,variable=\x] plot ({\x},{35/(\x+4.00167)-6.00167});
          \draw[domain=.3724:.997,smooth,variable=\x] plot ({\x},{35/(\x+4.00167)-7.00167});
                  \draw[domain=0:.3724,smooth,variable=\x] plot ({\x},{35/(\x+4.00167)-8.00167});
              \node at (1.48,.64) {$N=35, \alpha=4.00167$}; 
                      \node at (1.48,.5) {$|\T'(\alpha+1)|\approx1.40$}; 
                   \node at (1.48,.34) {$|\T'(\alpha)|\approx2.19$};  
    \end{tikzpicture}
   $$
\endminipage
 \caption{Two arrangements in which almost $\T^2(\alpha)<p_d<\T^4(\alpha+1)$\\
(in fact, in both cases $p_d = \T^2(\alpha )$).} \label{fig: optimal pd 10}
\end{figure}

\begin{figure}[!htb]
$$
\begin{tikzpicture}[scale =5] 
 \draw [dashed] (.43,0) -- (.43,1);
 \draw [dashed] (0,.75488) -- (.75488,.75488);
  \draw [dashed] (.75488,.75488) -- (.75488,.43);
  \draw [dashed] (.75488,.43) -- (.43,.43);
         \draw (0,0) -- (1,1);
     \draw (0,0) -- (0,1);
     \draw (0,1) -- (1,1);
     \draw (1,1) -- (1,0);
        \draw (0,0) -- (1,0);
    \node at (.86,.754) {$(a,a)$};
      \node at (.94,.435) {$(a,1-a^2)$};
      \node at (1.22,.293) {$y=-(a+1)x+a+1$};
            \node at (-.16,.5) {$y=-(a+1)x+a$};
  \draw[domain=.43:1,smooth,variable=\x] plot ({\x},{1.75488-1.75488*\x});
    \draw[domain=0:.43,smooth,variable=\x] plot ({\x},{.75488-1.75488*\x});
   \node at (-.03,.75488) {$a$};    
      \node at (.43,1.06) {$(\frac a{a+1},1)$};   
     \node at (-.01,-.04) {$0$};   
         \node at (1.01,-.04) {$1$};  
            \end{tikzpicture}
  $$
\caption{The `limit graph' of $\T$, translated over $(-\alpha,-\alpha)$, under the conditions $\I=\Delta_d\cup\Delta_{d-1}$ and $N/(N/\alpha-d)-(d-1)=p_d$ for $N\to \infty$ (and $\alpha, d \to \infty$ accordingly). This `arrangement' can be seen as one with three cylinders, where $\Delta_{d-2}\Mod \alpha$, the one on the right, is infinitely small; see also the arrangements in Figure \ref{fig: optimal pd 10}.} {\label{fig: limit arrangement 2}}
\end{figure}
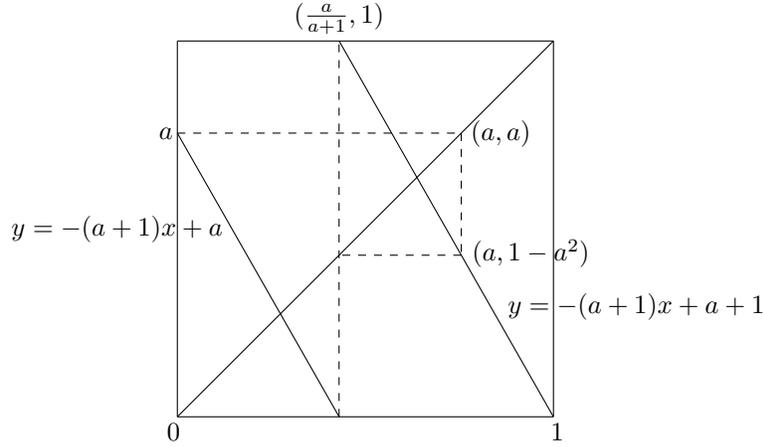

In order to find the relationship between $N$ and $d$ for arrangements with the conditions $\T^2(\alpha)<p_d$ and $\T^4(\alpha+1)>p_d$ mentioned above, we refer to some more relevant arrangements, as shown in Figure \ref{fig: optimal pd 10}. In both cases in Figure \ref{fig: optimal pd 10}, $\alpha$ is such that $\T^2(\alpha)=p_d$, which is a value of $\alpha$ that is only a little larger than the values for which $\T^2(\alpha)<p_d$ and $\T^4(\alpha+1)>p_d$. A `limit arrangement' (where the third, rightmost cylinder is infinitely small), similar to the `limit arrangement' used in the proof of Theorem \ref{gaplessness two cylinders}, is shown in Figure \ref{fig: limit arrangement 2}. The assumptions yield $a^3+a^2-1=0$, with real root $a=0.75487\cdots=:\gamma$.

Similar to the proof of Theorem \ref{gaplessness two cylinders} we then find that for arrangements as in Figure \ref{fig: optimal pd 10} we have 
$$
N\approx \frac{(d-1)(d-1+\gamma)(1+\gamma)}{\gamma^2}.
$$

Using this relationship, we can take a similar approach as in the proof of Proposition \ref{min slope}. We leave out the tedious steps and confine ourselves to observing that the slope of the line segments in Figure \ref{fig: limit arrangement 2} is $-(\gamma+1)=-1.75487\cdots$ and that in arrangements where $\T^2(\alpha)<p_d$ and $\T^4(\alpha+1)>p_d$, we will see that the slope $\T'(\alpha+1)$ approaches $-(\gamma+1)$ as $N$ tends to infinity. However, for our proof the inequality $|\T'(\alpha+1)|>1/2(\sqrt{5}+1)=1.61803\cdots=:G$ suffices, which turns out to hold for $N\in \N_{\geq 273}$. We will use this to show that for $N\in \N_{\geq 273}$ we have $|\T^3(R)| > |\T(\Delta_{d-2})|+|\T(R)|$. From this it immediately follows that we cannot have that both $\T(V_1)\subset \T(R)$ and $\T(V_2)\subset \T(\Delta_{d-2})$, and we are done with the proof of Lemma 11.\smallskip

Since $|\T'(x)|$ is a decreasing function on $\I$, and writing $\beta:=|\Delta_{d-2}|$, we have  
$$
 |\T(\Delta_{d-2})| >|\T'(\alpha+1)|\cdot\beta,\,\,{\text{so}}\,\, |R|>(|\T'(\alpha+1)|-1)\beta.
$$
It follows that 
$$
|\T(R)|>(|\T'(\alpha+1)|-1)\cdot|\T'(p_{d-1})|\beta,
$$
that
$$
|\T^2(R)|>(|\T'(\alpha+1)|-1)\cdot|\T'(p_{d-1})|\cdot|\T'(f_d)|\beta,
$$
and finally that
$$
|\T^3(R)|>(|\T'(\alpha+1)|-1)\cdot|\T'(p_{d-1})|^2\cdot|\T'(f_d)|\beta.
$$
We also have $|\T(\Delta_{d-2})|<|\T'(p_{d-1})|\beta$, so
$$
|R|<(|\T'(p_{d-1})|-1)\beta\,\,{\text{and}}\,\, |\T(R)|<|\T'(f_{d-1})|\cdot(|\T'(p_{d-1})|-1)\beta.
$$
It follows that 
\begin{align*}
|\T(\Delta_{d-2})|+|\T(R)|&<(|\T'(p_{d-1})|+|\T'(f_{d-1})|\cdot(|\T'(p_{d-1})|-1))\beta\\
&=(|\T'(p_{d-1})|-|\T'(f_{d-1})|+|\T'(f_{d-1})|\cdot|\T'(p_{d-1})|)\beta\\
&<|\T'(f_{d-1})|\cdot|\T'(p_{d-1})|\beta.
\end{align*}
So, although crudely, we certainly have that $|\T^3(R)| > |\T(\Delta_{d-2})|+|\T(R)|$ if 
$$
|\T'(f_{d-1})|\cdot|\T'(p_{d-1})|<(|\T'(\alpha+1)|-1)\cdot|\T'(p_{d-1})|^2\cdot|\T'(f_d)|,
$$
that is, if 
\begin{equation}\label{xx}
1<(|\T'(\alpha+1)|-1)\cdot|\T'(p_{d-1})|\cdot\frac{|\T'(f_d)|}{|\T'(f_{d-1})|}.
\end{equation}
Since 
$$
(|\T'(\alpha+1)|-1)\cdot|\T'(p_{d-1})|\cdot\frac{|\T'(f_d)|}{|\T'(f_{d-1})|}>(|\T'(\alpha+1)|-1)\cdot|\T'(p_{d-1})|>(|\T'(\alpha+1)|-1)\cdot|\T'(\alpha+1)|,
$$
we know that (\ref{xx}) holds for $|\T'(\alpha+1)|>G$, which in turn holds for all $N\in \N_{\geq 273}$. We remark that this value is quite a wide upper bound, since we did a rough approximation. Still, checking that we cannot have both $\T(V_1)\subset \T(R)$ and $\T(V_2)\subset \T(\Delta_{d-2})$ for smaller $N$ is not that hard and is left to the reader. This finishes the proof of Lemma \ref{wild case} (cf. case 2 in (\ref{conditions})).\hfill $\Box$

Lemma \ref{wild case} implies that in case $I_{\alpha}=\Delta_d\cup\Delta_{d-1}\cup\Delta_{d-2}$ and $f_{d-1}\leq \T(\alpha) \leq \T(\alpha+1)\,\,{\text{or}}\,\,\T(\alpha)\leq \T(\alpha+1) \leq f_{d-1}$ the division of an interval containing $p_d$ in two smaller ones cannot prevent an overall expansion that excludes any gaps. The other plausible case with three cylinder sets in which gaps might exist is when one outer cylinder set is very small, while the other one is full or nearly full, such that either $\T^3(\alpha+1)\geq\T(\alpha+1)$ (when $\Delta_{d-2}$ is the small cylinder set) or $\T^3(\alpha)\leq\T(\alpha)$ (when $\Delta_d$ is the small cylinder set). We will show that this is not possible either:

\begin{Lemma}\label{gaplessness three cylinders hard}
Let $I_{\alpha}=\Delta_d\cup\Delta_{d-1}\cup\Delta_{d-2}$. Then $\I$ is gapless if 
$$
f_{d-1}\leq \T(\alpha+1)\leq \T(\alpha)\,\,{\text{or}}\,\,\T(\alpha+1)\leq \T(\alpha)\leq f_{d-1}.
$$
\end{Lemma}

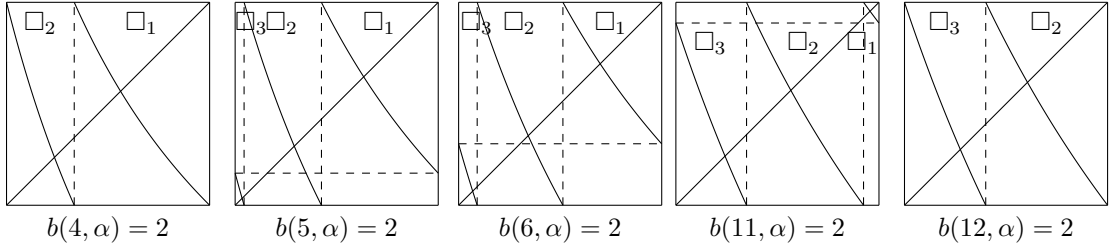
\begin{figure}[!htb]
\minipage{0.2\textwidth}
$$
\begin{tikzpicture}[scale =2.7] 
           \draw [dashed] (.33333,0) -- (.33333,1);
             \draw (0,0) -- (1,1);
     \draw (0,0) -- (0,1);
        \draw (0,0) -- (1,0);
     \draw (0,1) -- (1,1);
     \draw (1,1) -- (1,0);
      \draw[domain=.33333:1,smooth,variable=\x] plot ({\x},{4/(\x+1)-2});
   \draw[domain=0:.33333,smooth,variable=\x] plot ({\x},{4/(\x+1)-3});    
     \node at (.666,.9) {$\Box_1$};
      \node at (.167,.9) {$\Box_2$};
    \node at (.5,-.12) {$b(4,\alpha)=2$};   
    \end{tikzpicture}
 $$
\endminipage
\minipage{0.2\textwidth}
$$
\begin{tikzpicture}[scale =2.7] 
    \draw (0,0) -- (1,0);
         \draw [dashed] (.044,0) -- (.044,1);
        \draw [dashed] (.425,0) -- (.425,1);
          \draw [dashed] (0,.1583) -- (1,.1583);
             \draw (0,0) -- (1,1);
     \draw (0,0) -- (0,1);
     \draw (0,1) -- (1,1);
     \draw (1,1) -- (1,0);
     \draw[domain=0.425:1,smooth,variable=\x] plot ({\x},{5/(\x+1.1583)-2.1583});
           \draw[domain=.044:.425,smooth,variable=\x] plot ({\x},{5/(\x+1.1583)-3.1583});
         \draw[domain=0:.044,smooth,variable=\x] plot ({\x},{5/(\x+1.1583)-4.1583});
           \node at (.71,.9) {$\Box_1$};
      \node at (.23,.9) {$\Box_2$};
  \node at (.08,.9) {$\Box_3$};
         \node at (.5,-.12) {$b(5,\alpha)=2$};   
     \end{tikzpicture}
 $$
\endminipage
\minipage{0.2\textwidth}
$$
\begin{tikzpicture}[scale =2.7] 
     \draw (0,0) -- (1,0);
         \draw [dashed] (.0914,0) -- (.0914,1);
        \draw [dashed] (.5135,0) -- (.5135,1);
          \draw [dashed] (0,.303) -- (1,.303);
             \draw (0,0) -- (1,1);
     \draw (0,0) -- (0,1);
     \draw (0,1) -- (1,1);
     \draw (1,1) -- (1,0);
     \draw[domain=0.5135:1,smooth,variable=\x] plot ({\x},{6/(\x+1.303)-2.303});
           \draw[domain=.0914:.5135,smooth,variable=\x] plot ({\x},{6/(\x+1.303)-3.303});
         \draw[domain=0:.0914,smooth,variable=\x] plot ({\x},{6/(\x+1.303)-4.303});
            \node at (.75,.9) {$\Box_1$};
      \node at (.3,.9) {$\Box_2$};
  \node at (.09,.9) {$\Box_3$};
    \node at (.5,-.12) {$b(6,\alpha)=2$};   
        \end{tikzpicture}
 $$
\endminipage
\minipage{0.2\textwidth}
$$
\begin{tikzpicture}[scale =2.7] 
  \draw (0,0) -- (1,0);
         \draw [dashed] (.348,0) -- (.348,1);
        \draw [dashed] (.924,0) -- (.924,1);
          \draw [dashed] (0,.898) -- (1,.898);
             \draw (0,0) -- (1,1);
     \draw (0,0) -- (0,1);
     \draw (0,1) -- (1,1);
     \draw (1,1) -- (1,0);
     \draw[domain=0.924:1,smooth,variable=\x] plot ({\x},{11/(\x+1.898)-2.898});
           \draw[domain=.348:.924,smooth,variable=\x] plot ({\x},{11/(\x+1.898)-3.898});
         \draw[domain=0:.348,smooth,variable=\x] plot ({\x},{11/(\x+1.898)-4.898});
         \node at (.92,.8) {$\Box_1$};
      \node at (.63,.8) {$\Box_2$};
  \node at (.17,.8) {$\Box_3$};
  \node at (.5,-.12) {$b(11,\alpha)=2$};   
    \end{tikzpicture}
 $$
\endminipage
\minipage{0.2\textwidth}
$$
\begin{tikzpicture}[scale =2.7] 
 \draw (0,0) -- (1,0);
         \draw [dashed] (.4,0) -- (.4,1);
             \draw (0,0) -- (1,1);
     \draw (0,0) -- (0,1);
     \draw (0,1) -- (1,1);
     \draw (1,1) -- (1,0);
     \draw[domain=0.4:1,smooth,variable=\x] plot ({\x},{12/(\x+2)-4});
               \draw[domain=0:.4,smooth,variable=\x] plot ({\x},{12/(\x+2)-5});
  \node at (.7,.9) {$\Box_2$};
      \node at (.2,.9) {$\Box_3$};
        \node at (.5,-.12) {$b(12,\alpha)=2$};   
    \end{tikzpicture}
 $$
\endminipage
\caption{Arrangements of $\Upsilon_{N,\alpha}$ with $b=2$; in each case $\alpha = \frac{\sqrt{2N+1}-1}{2}$.} \label{fig: b=2}
\end{figure}

Proof of Lemma \ref{gaplessness three cylinders hard}: Taking into account our observations immediately following Definition \ref{small set} and the arrangements of Figure \ref{fig: optimal pd 6} for $N=23$ and $N=11$, we only have to prove that there are \emph{no} $\alpha$ such that $\T^3(\alpha)<\T(\alpha)$ is possible when $\T(\alpha+1)\leq \T(\alpha)\leq f_{d-1}$ (in case $\Delta_d$ is small) or such that $\T^3(\alpha+1)>\T(\alpha+1)$ is possible when $f_{d-1}\leq \T(\alpha+1)\leq \T(\alpha)$ (in case $\Delta_{d-2}$ is small). Note that the conditions $\T^3(\alpha)<\T(\alpha)$ and $\T^3(\alpha+1)>\T(\alpha+1)$ imply that the branch number is slightly larger than $2$. Now remember that $I_{\alpha}$ consists of $m$ full cylinder sets if and only if $\alpha=k$, $N=mk(k+1)$ and $d=(m-1)(k+1)$ for some $k \in \N$, cf.\ Theorem \ref{when I consists of m full cylinder sets}. Figure \ref{fig: b=2} shows for increasing values of $N$ a sequence of arrangements where the branch number $b$ is $2$, from one full arrangement (here for $N=4$) with two cylinders to the next one (here for $N=12$). Since $|\T'(\alpha)|>|\T'(\alpha+1)|$, the arrangements suggest that in case $b$ is slightly larger than 2, the most favourable arrangement for $\T^3(\alpha+1)= \T(\alpha+1)$ to have real roots is when $N=2k^2+2k-1$, where $k\in \N_{\geq 2}$, while for $\T^3(\alpha)= \T(\alpha)$ to have real roots is when $N=2k^2+2k+1$, where $k\in \N$. We will confine ourselves to investigating only the possibility of $\T^3(\alpha+1)= \T(\alpha+1)$; the calculations for the other case are similar. 

So we will try and find out if for $N=2k^2+2k-1$, $d=k+1$, with $k\in \N_{\geq 2}$, the positive root of $\T^3(\alpha+1)= \T(\alpha+1)$ lies in $\I$. To do this, we solve
$$
\cfrac{2k^2+2k-1}{\cfrac{\displaystyle
2k^2+2k-1}{\displaystyle \cfrac{2k^2+2k-1}{\alpha+1}-(k-1)}-k}-(k+1)=\frac{2k^2+2k-1}{\alpha+1}-(k-1),
$$
which is reducible to
$$
(2k^3+6k^2-k-1)\alpha^2+(2k^4+5k^2+k-2)\alpha-(4k^5+6k^4+2k^3-3k^2-k+1)=0,
$$
yielding
\begin{equation}\label{a for gap}
\alpha=\frac{\sqrt{36k^8+144k^7+164k^6-12k^5-95k^4-2k^3+21k^2-4k}-(2k^4+5k^2+k-2)}{2(2k^3+6k^2-k-1)}.
\end{equation}
A straightforward computation shows that this last expression is smaller than $f_{k+2}$, meaning that the root (\ref{a for gap}) lies outside $I_{\alpha}$ when $I_{\alpha}=\Delta_{k+1}\cup\Delta_k\cup\Delta_{k-1}$. Since $N=2k^2+2k-1$ was the most favourable option for investigation, this finishes our proof (cf.~case 3 in (\ref{conditions})).\hfill $\Box$

\begin{Remark}
The arrangement for $N=11$ in Figure \ref{fig: optimal pd 6} illustrates that the difference between $\T^3(\alpha+1)$ and $\T(\alpha+1)$ may be very small.
\end{Remark}
 
\subsection{A sufficient condition for gaplessness in case $\I$ consists of four cylinder sets}\label{four}\smallskip

In the previous subsection we proved Theorem \ref{gaplessness three cylinders or more} for $m=2$, by proving Lemmas \ref{gaplessness three cylinders easy}, \ref{wild case} and \ref{gaplessness three cylinders hard}. In this subsection we will consider $m=3$ and go into the analogons of Lemmas \ref{gaplessness three cylinders easy}, \ref{wild case} and \ref{gaplessness three cylinders hard}.\smallskip

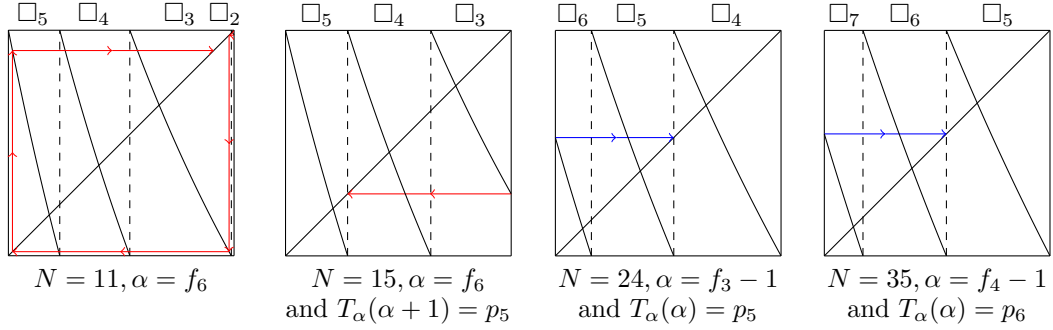
\begin{figure}[!htb]
\minipage{0.24\textwidth}
$$
\begin{tikzpicture}[scale =3] 
    \draw (0,0) -- (1,0);
    \draw [dashed] (.2276,0) -- (.2276,1);
      \draw [dashed] (.538,0) -- (.538,1);
        \draw [dashed] (.988,0) -- (.988,1);
               \draw (0,0) -- (1,1);
     \draw (0,0) -- (0,1);
     \draw (0,1) -- (1,1);
     \draw (1,1) -- (1,0);
     \draw[domain=0.988:1,smooth,variable=\x] plot ({\x},{11/(\x+1.472)-3.472});
      \draw[domain=0.538:.988,smooth,variable=\x] plot ({\x},{11/(\x+1.472)-4.472});
       \draw[domain=.2276:.538,smooth,variable=\x] plot ({\x},{11/(\x+1.472)-5.472});
         \draw[domain=0:.2276,smooth,variable=\x] plot ({\x},{11/(\x+1.472)-6.472});
                   \draw [->,red](1,.978) -- (.978,.978);
                     \draw [->,red](.978,.978) -- (.978,.495);
        \draw [->,red](.978,.495) -- (.978,.018);
   \draw [->,red](.978,.018) -- (.495,.018);
        \draw [->,red](.495,.018) -- (.018,.018);
           \draw [->,red](.018,.018) -- (.018,.46);
   \draw [->,red](.018,.46) -- (.018,.91);
   \draw [->,red](.018,.91) -- (.46,.91);
   \draw [->,red](.46,.91) -- (.91,.91);
   \node at (.96,1.07) {$\Box_2$};
      \node at (.76,1.07) {$\Box_3$};
  \node at (.38,1.07) {$\Box_4$};
   \node at (.11,1.07) {$\Box_5$};
 \node at (.5,-.11) {$N=11, \alpha=f_6$};   
   \node at (.5,-.3) {}; 
    \end{tikzpicture}
 $$
\endminipage
\minipage{0.24\textwidth}
$$
\begin{tikzpicture}[scale =3]
  \draw (0,0) --   (1,0);
            \draw [dashed] (.275,0) -- (.275,1);
        \draw [dashed] (.6438,0) -- (.6438,1);
       \draw (0,0) -- (1,1);
     \draw (0,0) -- (0,1);
     \draw (0,1) -- (1,1);
     \draw (1,1) -- (1,0);
            \draw [->,red] (1,.275) -- (.64,.275);
      \draw [->,red] (.64,.275) -- (.275,.275);
          \node at (.82,1.07) {$\Box_3$};
                   \node at (.46,1.07) {$\Box_4$};
           \node at (.14,1.07) {$\Box_5$};
       \draw[domain=.6438:1,smooth,variable=\x] plot ({\x},{15/(\x+1.899)-4.899});
          \draw[domain=.275:.6438,smooth,variable=\x] plot ({\x},{15/(\x+1.899)-5.899});
                  \draw[domain=0:.275,smooth,variable=\x] plot ({\x},{15/(\x+1.899)-6.899});
                \node at (.5,-.11) {$N=15,\alpha=f_6$}; 
                   \node at (.5,-.24) {and $\T(\alpha+1)=p_5$}; 
              \end{tikzpicture}
   $$
\endminipage
\minipage{0.24\textwidth}
$$
\begin{tikzpicture}[scale =3]
  \draw (0,0) --   (1,0);
         \draw [dashed] (.1596,0) -- (.1596,1);
             \draw [dashed] (.5247,0) -- (.5247,1);
       \draw (0,0) -- (1,1);
     \draw (0,0) -- (0,1);
     \draw (0,1) -- (1,1);
     \draw (1,1) -- (1,0);
          \draw [->,blue](0,.5247) -- (.27,.5247);
         \draw [->,blue](.27,.5247) -- (.5247,.5247);
                    \node at (.76,1.07) {$\Box_4$};
      \node at (.34,1.07) {$\Box_5$};
       \node at (.08,1.07) {$\Box_6$};
               \draw[domain=.5247:1,smooth,variable=\x] plot ({\x},{24/(\x+2.623475)-6.623475});
                      \draw[domain=.1596:.5247,smooth,variable=\x] plot ({\x},{24/(\x+2.623475)-7.623475});
                                 \draw[domain=0:.1596,smooth,variable=\x] plot ({\x},{24/(\x+2.623475)-8.623475});
                                    \node at (.5,-.11) {$N=24, \alpha=f_3-1$}; 
                                            \node at (.5,-.24) {and $\T(\alpha)=p_5$}; 
                           \end{tikzpicture}
$$
\endminipage
\minipage{0.24\textwidth}
$$
\begin{tikzpicture}[scale =3]
  \draw (0,0) --   (1,0);
         \draw [dashed] (.1713,0) -- (.1713,1);
             \draw [dashed] (.541,0) -- (.541,1);
       \draw (0,0) -- (1,1);
     \draw (0,0) -- (0,1);
     \draw (0,1) -- (1,1);
     \draw (1,1) -- (1,0);
         \draw [->,blue](0,.541) -- (.27,.541);
         \draw [->,blue](.27,.541) -- (.541,.541);
                 \node at (.77,1.07) {$\Box_5$};
      \node at (.36,1.07) {$\Box_6$};
       \node at (.09,1.07) {$\Box_7$};
                      \draw[domain=.541:1,smooth,variable=\x] plot ({\x},{35/(\x+3.245)-8.245});
                                 \draw[domain=.1713:.541,smooth,variable=\x] plot ({\x},{35/(\x+3.245)-9.245});
                                        \draw[domain=0:.1713,smooth,variable=\x] plot ({\x},{35/(\x+3.245)-10.245});
                                    \node at (.5,-.11) {$N=35, \alpha=f_4-1$}; 
                                            \node at (.5,-.24) {and $\T(\alpha)=p_6$}; 
                           \end{tikzpicture}
   $$
\endminipage
\caption{Four arrangements with two full cylinders} \label{fig: two two}
\end{figure}

When $\I$ consists of four cylinder sets, the analogon of Lemma \ref{gaplessness three cylinders easy} is that arrangements $I_{\alpha}$ are gapless when $I_{\alpha}=\Delta_d\cup\Delta_{d-1}\cup \Delta_{d-2}\cup \Delta_{d-3}$ while $\T(\alpha)\geq f_{d-1}$ {\bf{and}} $\T(\alpha+1)\leq f_{d-2}$. The analogon of Lemma \ref{gaplessness three cylinders hard} is that arrangements $I_{\alpha}$ are gapless when $I_{\alpha}=\Delta_d\cup\Delta_{d-1}\cup\Delta_{d-2}\cup\Delta_{d-3}$ while $f_{d-2}\leq \T(\alpha+1)\leq \T(\alpha)\,\,{\text{{\bf{or}}}}\,\,\T(\alpha+1)\leq \T(\alpha)\leq f_{d-1}$. In both cases branch numbers larger than $3$ are involved, in which case $|\T'(\alpha+1)|>2$ when $N\in \N_{\geq 18}$ (and Theorem \ref{Wilkinson} yields the desired result). The cases $2\leq N\leq 17$ can be checked manually and are left to the reader; in Figure \ref{fig: two two} the arrangement for $N=11$, associated with Lemma \ref{gaplessness three cylinders hard}, illustrates that gaps are out of the question. 

The analogon of Lemma \ref{wild case} is that arrangements $I_{\alpha}$ are gapless when $I_{\alpha}=\Delta_d\cup\Delta_{d-1}\cup \Delta_{d-2}\cup \Delta_{d-3}$ while $f_{d-1}\leq \T(\alpha) \leq \T(\alpha+1)\,\,{\text{{\bf{or}}}}\,\,\T(\alpha)\leq \T(\alpha+1) \leq f_{d-2}$. The arrangements for $N=15$, $N=24$ and $N=35$ in Figure \ref{fig: two two} are interesting illustrations of the analogon of Lemma \ref{wild case} in the case of two full cylinder sets instead of one. We will confine ourselves to the arrangement for $N=15$; the other ones have similar properties. 

\begin{figure}[!htb]
$$
\begin{tikzpicture}[scale =4.5] 
 \draw [dashed] (.382,0) -- (.382,1);
 \draw [dashed] (1,.382) -- (.382,.382);
  \draw [dashed] (.764,0) -- (.764,1);
         \draw (0,0) -- (1,1);
     \draw (0,0) -- (0,1);
     \draw (0,1) -- (1,1);
     \draw (1,1) -- (1,0);
        \draw (0,0) -- (1,0);
         \node at (1.1,.382) {$(1,a)$};
      \node at (1.1,.7) {$y=-\frac1ax+3$};
            \node at (.91,.2) {$y=-\frac1ax+2$};
             \node at (-.03,.5) {$y=-\frac1ax+1$};
  \draw[domain=.764:1,smooth,variable=\x] plot ({\x},{3-2.618034*\x});
    \draw[domain=.382:.764,smooth,variable=\x] plot ({\x},{2-2.618034*\x});
      \draw[domain=0:.382,smooth,variable=\x] plot ({\x},{1-2.618034*\x});
      \node at (.382,1.06) {$(a,1)$};   
     \node at (-.01,-.04) {$0$};   
         \node at (1.01,-.04) {$1$};  
            \end{tikzpicture}
  $$
\caption{The `limit graph' of $\T$, translated over $(-\alpha,-\alpha)$, under the conditions $\I=\Delta_d\cup\Delta_{d-1}\Delta_{d-2}$ and $N/(\alpha+1)-(d-2)=p_d$ for $N\to \infty$ (and $\alpha, d \to \infty$ accordingly)} {\label{fig: limit arrangement 3}}
\end{figure}
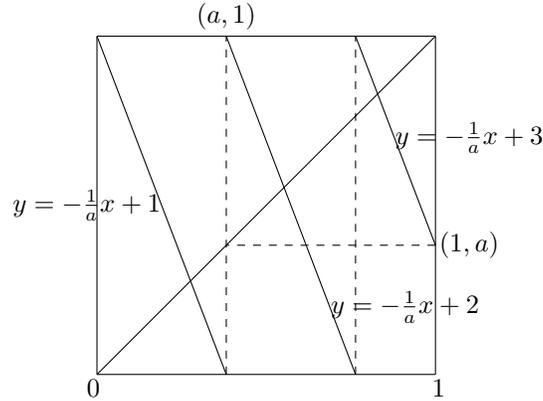

The arrangement for $N=15$ is the boundary case for the situation where we have four cylinders, the left one of which (that would be $\Delta_6$ in this example) is extremely small and the right one is such that almost $p_5\in\T^2(\T(\Delta_6)\setminus\Delta_6)$. The interesting thing is that this option would imply a quick interruption of the expansion of $\T(\Delta_6)\setminus\Delta_6$, involving two large gaps. But it is not really an option: the arrangement for $N=15$ in Figure \ref{fig: two two} is exceptional among relatively small $N$ (as well as the arrangements for $N=24$ and $N=35$ are), while for $N>36$ we have $|\T'(\alpha+1)|>2$ when $\I=\Delta_d\cup\Delta_{d-1}\Delta_{d-2}$ and $N/(\alpha+1)-(d-2)=p_d$ or $N/\alpha-d=p_d$. We derived this in a similar way as in the proof of Lemma \ref{wild case} (see Figure \ref{fig: limit arrangement 2}) or the preparations for Theorem \ref{gaplessness two cylinders} (see Figure \ref{fig: limit arrangement}). Figure \ref{fig: limit arrangement 3} shows the associated `limit graph', from which it is easily found that $a=1/2(3-\sqrt{5})$, yielding branch number $2+g$, with $g=1/G$ the small golden section.

With this, we conclude the proof of Theorem \ref{gaplessness three cylinders or more}. \hfill $\Box$\smallskip

In the next section we will prove that if $\I=\Delta_d\cup\Delta_{d-1}\cup\Delta_{d-2}\cup\Delta_{d-3}$ gaps exist only in very rare cases and if they do, that they are very large. After that, we will finish the proof of Theorem \ref{no gaps five cylinder case}, stating that all arrangements with five cylinders are gapless.

\section{Gaplessness in case $\I$ contains two full cylinder sets}\label{four and five}

When an arrangement contains three cylinders, and two of them are full, the arrangement is gapless according to Theorem  \ref{gaplessness three cylinders or more}. In this section we will proof that arrangements of four cylinders generally do not contain a gap either, save for special values of $N$. The core of this proof rests on values of $\alpha$ satisfying one of the equations
$$
\T(\alpha)=\T^3(\alpha)\,\,{\text{(with root $\alpha_{\ell}$)}}\,\, {\text{and}}\,\, \T(\alpha+1)=\T^3(\alpha+1)\,\,{\text{(with root $\alpha_u$)}}.
$$
We will show that for $N$ such that $\alpha_{\ell}<\alpha_u$ very large gaps exist for $\alpha\in[\alpha_{\ell},\alpha_u]$.\smallskip  

The central theorem of this section is the following:
\begin{Theorem}\label{gaps four cylinder case}
Let $N \in \N_{\geq 2}$ and $\I=\Delta_d\cup\Delta_{d-1}\cup\Delta_{d-2}\cup\Delta_{d-3}$. Then there is a gap in $\I$ if and only if $N=2k^2+2k-i$, with $k > 1$ and $i \in \{1,2,3\}$. Moreover, if there is a gap in $\I$, the gap contains $f_{d-1}$ and $f_{d-2}$, while $\Delta_d$ and $\Delta_{d-3}$ are gapless.
\end{Theorem}

Proof of Theorem \ref{gaps four cylinder case}: Suppose that there is a gap containing $f_{d-1}$ and $f_{d-2}$ in $\I$ and that $\Delta_d$ and $\Delta_{d-3}$ are gapless. Then, as a sub-interval of a gap, the interval $(f_{d-1},f_{d-2})$ is a gap. Since $f_{d-1}<f_{d-2}$, $N/(f_{d-1}+d-1)=f_{d-1}$ and $N/(f_{d-2}+d-2)=f_{d-2}$, we know that 
$$
(f_{d-1},f_{d-2}) \subsetneq \left (\frac N{f_{d-2}+d-1},\frac N{f_{d-1}+d-2}\right ),
$$  
where the larger open interval is a gap as well. What is more, the infinite sequence of intervals
$$
(f_{d-1},f_{d-2}) \subsetneq \left (\frac N{f_{d-2}+d-1},\frac N{f_{d-1}+d-2}\right )\subsetneq  \left (\frac N{\frac N {f_{d-1}+d-2}+d-1} ,\frac N{\frac N{f_{d-2}+d-1}+d-2}\right )\subsetneq \ldots
$$
consists of the union of $(f_{d-1},f_{d-2})$ with pre-images of $(f_{d-1},f_{d-2})$ in $\Delta_{d-1}$ and $\Delta_{d-2}$ respectively and therefore of gaps containing $f_{d-1}$ and $f_{d-2}$. It is contained in the closed interval
$[q,r]$, with
$$
q=[\overline{d-1,d-2}]_{N,\alpha}\in \Delta_{d-1} \quad {\text{and}} \quad r=[\overline{d-2,d-1}]_{N,\alpha}\in \Delta_{d-2},
$$
yielding 
\begin{equation}\label{conditions T^2}
T_{\alpha}^2(q)=q, \,\, T_{\alpha}(q)=r,\,\, T_{\alpha}(r)=q\,\,{\text{and}}\,\, T_{\alpha}^2(r) =r.
\end{equation}

Since $\Delta_d$ and $\Delta_{d-3}$ are gapless, $\T(\alpha)$ and $\T(\alpha+1)$ lie outside the interval $(q,r)$, which is to say 
$$
p_d<\T(\alpha)\leq q\,\,{\text{and}}\,\, r\leq \T(\alpha+1)<p_{d-2}.
$$

For the images of $\alpha$ under $\T$ this means that either $\T^2(\alpha)\in \Delta_{d-3}$ or $\T^2(\alpha)\in\Delta_{d-2}$, in the latter case of which we have, due to the expansiveness of $\T$ and the equalities of (\ref{conditions T^2}),
$$
|\T(\alpha)-q| \leq |\T^2(\alpha) - r|\leq |\T^3(\alpha) - q|,
$$
with equalities only in the case $\T(\alpha)=q$. From this we derive that
\begin{equation}\label{necessity left}
{\text{either}}\,\, \T^2(\alpha)\in \Delta_{d-3}\,\,{\text{or}}\,\, \T^2(\alpha)\in\Delta_{d-2}\wedge \T^3(\alpha) \leq \T(\alpha)
\end{equation}
and, similarly, that
\begin{equation}\label{necessity right}
{\text{either}}\,\, \T^2(\alpha+1)\in \Delta_d\,\,{\text{or}}\,\,  \T^2(\alpha+1)\in \Delta_{d-1}\wedge \T^3(\alpha+1) \geq \T(\alpha+1).
\end{equation}
\smallskip
In the following we will write $\alpha_u(N,m)$ ($u$ for `upper') for the positive root of the equation $\T^3(\alpha+1) = \T(\alpha+1)$ (so $\T(\alpha+1)=r$) and $\alpha_{\ell}(N,m)$\label{lower} ($l$ for `lower') for the positive root of the equation $\T^3(\alpha) = \T(\alpha)$ (so $\T(\alpha)=q$), with $m$ the number of full cylinder sets; in the current case we have $m=2$. Recall that $I_{\alpha}$ consists of $m$ full cylinder sets if and only if $\alpha=k$, $N=mk(k+1)$ and $d=(m-1)(k+1)$ for some $k \in \N$, cf.\ Theorem \ref{when I consists of m full cylinder sets}, so when $m=2$, we have arrangements consisting of two full arrangements only for $\alpha=k$, $N=2k(k+1)$ and $d=k+1$. If $N$ is $2k(k+1)-n$, with $n \in \N$, and $\alpha=k-x$, with $x\in \R$, we have 
$$
b(N,\alpha)=\frac{2k^2+2k-n}{(k-x)(k+1-x)}=2+\frac{4xk-n+2x-2x^2}{(k-1)(k+1-x)}>2+\frac{4xk-n}{k^2+x^2},
$$
which is a little bit larger than $2$ provided $x$ and $n$ are relatively small. For these arrangements we have
$$
d(\alpha)=\left \lfloor \frac{2k^2+2k-n}{k-x}-(k-x)\right \rfloor=\left \lfloor k+2+3x+\frac{2x^2+2x-n}{k-x}\right \rfloor=k+2
$$
and
$$
d(\alpha+1)=\left \lfloor \frac{2k^2+2k-n}{k+1-x}-(k+1-x)\right \rfloor=\left \lfloor k-1+3x+\frac{2x^2+2x-n}{k+1-x}\right \rfloor=k-1.
$$
So, for $x$ and $n$ relatively small, the arrangements consist of four cylinders, while the branch number is only a little bit larger than $2$. We will now use this to finish the forward implication of Theorem \ref{gaps four cylinder case}.\smallskip

Since $\Delta_{d-3}$ decreases and $\Delta_d$ increases as $\alpha$ decreases, we see that the assumption that there is a gap containing $f_{d-1}$ and $f_{d-2}$ in $\I$ implies $\alpha_u(N,2)\geq \alpha_{\ell}(N,2)$. We will shortly show that the only values of $N$ for which $\alpha_u(N,2)\geq \alpha_{\ell}(N,2)$ are $N=2k^2+2k-i$, with $k > 1$ and $i \in \{1,2,3\}$; in all cases $d=k+2$. Although we could keep $i$ as a variable in our calculations, we can limit ourselves to the case $i=3$, since $i=3$ is the least favourable value of $i$ allowing for a gap, as is suggested\footnote{Note that in Figure \ref{fig: b=2} we have $d=k+1$, while there is no small cylinder set $\Delta_{k+2}$} in Figures \ref{fig: b=2} and \ref{fig: N=11, alfa=1.8719} through \ref{fig: N=9, alfa=1.5946}. We will show that for $i=3$ indeed $\alpha_u(N,2)\geq \alpha_{\ell}(N,2)$. Subsequently we will show that for $4\leq i\leq 4k$ no gaps exist; the upper bound is $4k$, since $2k^2+2k-4k=2(k-1)^2+2(k-1)$, so as to confine the calculations to the group of arrangements where $d=k+2$.

\begin{figure}[!htb]
\minipage{0.44\textwidth}
$$
\begin{tikzpicture}[scale =7.2] 
\draw[draw=white,fill=gray!20!white] 
plot[smooth,samples=100,domain=0:.0013] (\x,0) --
plot[smooth,samples=100,domain=.0013:0] (\x,.004);
\draw[draw=white,fill=gray!20!white] 
plot[smooth,samples=100,domain=.0013:.0148] (\x,0) --
plot[smooth,samples=100,domain=.0148:.0013] (\x,1);
 \draw[draw=white,fill=gray!20!white] 
plot[smooth,samples=100,domain=.013:1] (\x,.958) --
plot[smooth,samples=100,domain=1:.013] (\x,1);
 \draw[draw=white,fill=gray!20!white] 
plot[smooth,samples=100,domain=.958:.969] (\x,0) --
plot[smooth,samples=100,domain=.969:.958] (\x,.958);
 \draw[draw=white,fill=gray!20!white] 
plot[smooth,samples=100,domain=.0148:.958] (\x,0) --
plot[smooth,samples=100,domain=.958:.0148] (\x,.0148);
 \draw[black,fill=black] (.0009,.0009) circle (.04ex);
\draw[black,fill=black] (.268,.268) circle (.04ex);
\draw[black,fill=black] (.592,.592) circle (.04ex);
\draw[black,fill=black] (.982,.982) circle (.04ex);
\draw[black,fill=black] (0,0) circle (.04ex);
\draw[black,fill=black] (.0148,0) circle (.04ex);
\draw[black,fill=black] (.958,0) circle (.04ex);
\draw[black,fill=black] (1,0) circle (.04ex);
    \draw (0,0) -- (.0148,0);
           \draw (.958,0) -- (1,0);
               \draw [dashed] (.0013,0) -- (.0013,1);
                        \draw [dashed] (.3858,0) -- (.3858,1);
         \draw [dashed] (.969,0) -- (.969,1);
        \draw (0,0) -- (1,1);
     \draw (0,0) -- (0,1);
     \draw (0,1) -- (1,1);
     \draw (1,1) -- (1,0);
         \draw [->,blue](0,.004) -- (.002,.004);
        \draw [->,blue](.002,.004) -- (.004,.004);
        \draw [->,blue](.004,.004) -- (.004,.4977);
       \draw [->,blue](.004,.4977) -- (.004,.991);
           \draw [->,blue](.004,.991) -- (.4977,.991);
              \draw [->,blue](.4977,.991) -- (.991,.991);
         \draw [->,red] (1,.958) -- (.979,.958);
             \draw [->,red] (.979,.958) -- (.958,.958);
    \draw [->,red] (.958,.958) -- (.958,.4864);
    \draw [->,red] (.958,.4864) -- (.958,.0148);
    \draw [->,red] (.958,.0148) -- (.4864,.0148);
    \draw [->,red] (.4864,.0148) -- (.0148,.0148);
     \draw [->,red] (.0148,.0148) -- (.0148,.486);
    \draw [->,red] (.0148,.486) -- (.0148,.958);
    \draw [->,red] (.0148,.958) -- (.486,.958);
    \draw [->,red] (.486,.958) -- (.958,.958);
                   \draw[domain=.969:1,smooth,variable=\x] plot ({\x},{11/(\x+1.872)-2.872});
       \draw[domain=.3858:.969,smooth,variable=\x] plot ({\x},{11/(\x+1.872)-3.872});
             \draw[domain=.0013:.3858,smooth,variable=\x] plot ({\x},{11/(\x+1.872)-4.872});
       \draw[domain=0:.0013,smooth,variable=\x] plot ({\x},{11/(\x+1.872)-5.872});
          \node at (.5,-.05) {$N=11, \alpha=1.8719\cdots$\,\,($k=2$, $i=1$)}; 
                        \node at (.5,-.12) {$\T^3(\alpha+1)=\T(\alpha+1)$}; 
                           \draw [->] (.04,.065) -- (.025,.03);
          \draw [->] (.93,.96) -- (.98,.98);
           \node at (.05,.09) {$F_4$};
      \node at (.31,.2715) {$F_3$};
       \node at (.635,.5955) {$F_2$};
      \node at (.91,.95) {$F_1$};
         \node at (1,.8) {$\Box_1$};
      \node at (.68,.8) {$\Box_2$};
            \node at (.2,.8) {$\Box_3$};
      \node at (0.02,.8) {$\Box_4$};
    \end{tikzpicture}
 $$
 \vspace*{-8mm} \caption{\label{fig: N=11, alfa=1.8719}}
\endminipage
\hfill
\minipage{0.44\textwidth}
$$
\begin{tikzpicture}[scale =7.2] 
\draw[draw=white,fill=gray!20!white] 
plot[smooth,samples=100,domain=0:.0058] (\x,0) --
plot[smooth,samples=100,domain=.0058:0] (\x,.018);
\draw[draw=white,fill=gray!20!white] 
plot[smooth,samples=100,domain=.0058:.018] (\x,0) --
plot[smooth,samples=100,domain=.018:.0058] (\x,1);
 \draw[draw=white,fill=gray!20!white] 
plot[smooth,samples=100,domain=.018:.975] (\x,.9615) --
plot[smooth,samples=100,domain=.975:.018] (\x,1);
 \draw[draw=white,fill=gray!20!white] 
plot[smooth,samples=100,domain=.975:1] (\x,.966) --
plot[smooth,samples=100,domain=1:.975] (\x,1);
 \draw[draw=white,fill=gray!20!white] 
plot[smooth,samples=100,domain=.9615:.975] (\x,0) --
plot[smooth,samples=100,domain=.975:.9615] (\x,.9615);
 \draw[draw=white,fill=gray!20!white] 
plot[smooth,samples=100,domain=.0058:.9615] (\x,0) --
plot[smooth,samples=100,domain=.9615:.0058] (\x,.018);
 \draw[black,fill=black] (.0044,.0044) circle (.04ex);
\draw[black,fill=black] (.2715,.2715) circle (.04ex);
\draw[black,fill=black] (.5955,.5955) circle (.04ex);
\draw[black,fill=black] (.9855,.9855) circle (.04ex);
\draw[black,fill=black] (0,0) circle (.04ex);
\draw[black,fill=black] (.0165,0) circle (.04ex);
\draw[black,fill=black] (.965,0) circle (.04ex);
\draw[black,fill=black] (1,0) circle (.04ex);
    \draw (0,0) -- (.0165,0);
          \draw (.965,0) -- (1,0);
               \draw [dashed] (.0058,0) -- (.0058,1);
                        \draw [dashed] (.391,0) -- (.391,1);
         \draw [dashed] (.975,0) -- (.975,1);
        \draw (0,0) -- (1,1);
     \draw (0,0) -- (0,1);
     \draw (0,1) -- (1,1);
     \draw (1,1) -- (1,0);
         \draw [->,blue](0,.018) -- (.009,.018);
        \draw [->,blue](.009,.018) -- (.018,.018);
        \draw [->,blue](.018,.018) -- (.018,.49);
       \draw [->,blue](.018,.49) -- (.018,.9615);
           \draw [->,blue](.018,.9615) -- (.49,.9615);
              \draw [->,blue](.49,.9615) -- (.9615,.9615);
                  \draw [->,blue](.9615,.9615) -- (.9615,.49);
       \draw [->,blue](.9615,.49) -- (.9615,.018);
           \draw [->,blue](.9615,.018) -- (.49,.018);
              \draw [->,blue](.49,.018) -- (.018,.018);
         \draw [->,red] (1,.966) -- (.983,.966);
          \draw [->,red] (.983,.966) -- (.966,.966);
    \draw [->,red] (.966,.966) -- (.966,.489);
    \draw [->,red] (.966,.489) -- (.966,.012);
    \draw [->,red] (.966,.012) -- (.489,.012);
    \draw [->,red] (.489,.012) -- (.012,.012);
\draw[domain=.975:1,smooth,variable=\x] plot ({\x},{11/(\x+1.8686)-2.8686});
       \draw[domain=.391:.975,smooth,variable=\x] plot ({\x},{11/(\x+1.8686)-3.8686});
             \draw[domain=.0058:.391,smooth,variable=\x] plot ({\x},{11/(\x+1.8686)-4.8686});
       \draw[domain=0:.0058,smooth,variable=\x] plot ({\x},{11/(\x+1.8686)-5.8686});
       \draw [->] (.04,.065) -- (.025,.03);
          \draw [->] (.93,.96) -- (.98,.98);
           \node at (.05,.09) {$F_4$};
      \node at (.31,.2715) {$F_3$};
       \node at (.635,.5955) {$F_2$};
      \node at (.91,.95) {$F_1$};
         \node at (1,.8) {$\Box_1$};
      \node at (.68,.8) {$\Box_2$};
            \node at (.2,.8) {$\Box_3$};
      \node at (0.02,.8) {$\Box_4$};
          \node at (.5,-.05) {$N=11, \alpha=1.8687\cdots \,\,(k=2, i=1)$}; 
                        \node at (.5,-.12) {$\T^3(\alpha)=\T(\alpha)$};       
    \end{tikzpicture}
 $$
\vspace*{-8mm} \caption{\label{fig: N=11, alfa=1.8687}}
\endminipage
\end{figure}

So, let $N=2k^2+2k-3$ and $d=k+2$. Then $\alpha_{\ell}(N,2) = [k+2,\overline{k+1,k}]$ and $\alpha_u(N,2)+1= [k-1,\overline{k,k+1}]$\footnote{We omit the suffix `$N,\alpha$' behind these expansions not only for eligibility but also because $\alpha$ has yet to be determined as the root of $\T^3(\alpha)=\T(\alpha)$ or $\T^3(\alpha+1)=\T(\alpha+1)$.}. Omitting straightforward calculations, we find that
$$
\alpha_u(2k^2+2k-3,2)=\frac{(2k^2+2k-3)\sqrt{D}-(2k^4+3k^2+3k-6)}{4k^3+12k^2-6k-6}
$$
and 
$$
\alpha_{\ell}(2k^2+2k-3,2)=\frac{(2k^2+2k-3)(\sqrt{D}-(k^2+5k+4))}{4k^3-18k-8},
$$
with
$$
D=9k^4 + 18k^3 -3k^2 - 12k=(3k^2+3k-2)^2-4.
$$
Since we assume that there is a gap containing $f_{d-1}$ and $f_{d-2}$, we have $\alpha_u\geq\alpha_{\ell}$. Omitting the basic calculations, we find that this inequality holds for $k\in \N_{\geq 2}$.
In the case $k=3$ (and so $N=21$), we have indeed
$$
\alpha_u(2k^2+2k-3,2)=\frac{\sqrt{508032}-192}{192}=2.7123\cdots>2.7122\cdots=\frac{\sqrt{508032}-588}{46}=\alpha_{\ell}(2k^2+2k-3,2);
$$
see Figure \ref{fig: N=21, alfa=2.7123}.
\smallskip

Some more basic calculations show that the cases $N=2k^2+2k-1$ and $N=2k^2+2k-2$ allow for larger intervals $[\alpha_{\ell},\alpha_u]$ where large gaps exist; see the next examples.
\begin{align*}
&\alpha_u(11,2)=\frac{\sqrt{9075}-26}{37}=1.8719\cdots\quad {\text{and}} \quad  \alpha_{\ell}(11,2)=\frac{99-\sqrt{9075}}2=1.8686\cdots \\
&\alpha_u(10,2)=\frac{\sqrt{1725}-12}{17}=1.7372\cdots \quad {\text{and}} \quad \alpha_{\ell}(10,2)=\frac{45-\sqrt{1725}}2=1.7334\cdots\\
&\alpha_u(9,2)=\frac{\sqrt{5103}-22}{31}=1.5946\cdots \quad {\text{and}} \quad \alpha_{\ell}(9,2)=\frac{27-\sqrt{567}}2=1.5941\cdots\\
&\alpha_u(8,2)=\frac{\sqrt{228}-5}7=1.4428\cdots \quad {\text{and}} \quad \alpha_{\ell}(8,2)=9-\sqrt{57}=1.4501\cdots
\end{align*}
We see that the intervals $\alpha_u-\alpha_{\ell}$ decrease as $N$ decreases, until (for $N=8$) the `interval' would have negative length, hence does not exist.\smallskip

\begin{figure}[!htb]
\minipage{0.44\textwidth}
$$
\begin{tikzpicture}[scale =6] 
 \draw[black,fill=black] (.0047,.0047) circle (.04ex);
\draw[black,fill=black] (.263,.263) circle (.04ex);
\draw[black,fill=black] (.58,.58) circle (.04ex);
\draw[black,fill=black] (.965,.965) circle (.04ex);
\draw[black,fill=black] (0,0) circle (.04ex);
\draw[black,fill=black] (.03,0) circle (.04ex);
\draw[black,fill=black] (.917,0) circle (.04ex);
\draw[black,fill=black] (1,0) circle (.04ex);
    \draw (0,0) -- (.03,0);
           \draw (.917,0) -- (1,0);
               \draw [dashed] (.006,0) -- (.006,1);
                        \draw [dashed] (.374,0) -- (.374,1);
         \draw [dashed] (.939,0) -- (.939,1);
        \draw (0,0) -- (1,1);
     \draw (0,0) -- (0,1);
     \draw (0,1) -- (1,1);
     \draw (1,1) -- (1,0);
         \draw [->,blue](0,.02) -- (.01,.02);
        \draw [->,blue](.01,.02) -- (.02,.02);
        \draw [->,blue](.02,.02) -- (.02,.485);
       \draw [->,blue](.02,.485) -- (.02,.954);
           \draw [->,blue](.02,.954) -- (.485,.954);
              \draw [->,blue](.485,.954) -- (.954,.954);
         \draw [->,red] (1,.917) -- (.979,.917);
             \draw [->,red] (.979,.917) -- (.917,.917);
    \draw [->,red] (.917,.917) -- (.917,.47);
    \draw [->,red] (.917,.47) -- (.917,.03);
    \draw [->,red] (.917,.03) -- (.47,.03);
    \draw [->,red] (.47,.03) -- (.03,.03);
     \draw [->,red] (.03,.03) -- (.03,.47);
    \draw [->,red] (.03,.47) -- (.03,.918);
    \draw [->,red] (.03,.918) -- (.47,.918);
    \draw [->,red] (.47,.918) -- (.918,.918);
    \draw[domain=.939:1,smooth,variable=\x] plot ({\x},{10/(\x+1.737)-2.737});
       \draw[domain=.374:.939,smooth,variable=\x] plot ({\x},{10/(\x+1.737)-3.737});
             \draw[domain=.006:.374,smooth,variable=\x] plot ({\x},{10/(\x+1.737)-4.737});
       \draw[domain=0:.006,smooth,variable=\x] plot ({\x},{10/(\x+1.737)-5.737});
          \node at (.5,-.05) {$N=10, \alpha=1.7372\cdots$\,\,($k=2$, $i=2$)}; 
                        \node at (.5,-.12) {$\T^3(\alpha+1)=\T(\alpha+1)$}; 
                           \draw [->] (.04,.065) -- (.025,.03);
               \node at (.05,.09) {$F_4$};
      \node at (.31,.26) {$F_3$};
       \node at (.62,.57) {$F_2$};
      \node at (.93,.97) {$F_1$};
         \node at (.97,.8) {$\Box_1$};
      \node at (.68,.8) {$\Box_2$};
            \node at (.2,.8) {$\Box_3$};
      \node at (0.02,.8) {$\Box_4$};
    \end{tikzpicture}
 $$
 \vspace*{-8mm} \caption{\label{fig: N=10, alfa=1.7372}}
\endminipage
\hfill
\minipage{0.44\textwidth}
$$
\begin{tikzpicture}[scale =6] 
 \draw[black,fill=black] (.011,.011) circle (.04ex);
\draw[black,fill=black] (.2595,.2595) circle (.04ex);
\draw[black,fill=black] (.568,.568) circle (.04ex);
\draw[black,fill=black] (.947,.947) circle (.04ex);
\draw[black,fill=black] (0,0) circle (.04ex);
\draw[black,fill=black] (.05,0) circle (.04ex);
\draw[black,fill=black] (.874,0) circle (.04ex);
\draw[black,fill=black] (1,0) circle (.04ex);
    \draw (0,0) -- (.05,0);
          \draw (.874,0) -- (1,0);
               \draw [dashed] (.014,0) -- (.014,1);
                        \draw [dashed] (.364,0) -- (.364,1);
         \draw [dashed] (.909,0) -- (.909,1);
        \draw (0,0) -- (1,1);
     \draw (0,0) -- (0,1);
     \draw (0,1) -- (1,1);
     \draw (1,1) -- (1,0);
         \draw [->,blue](0,.05) -- (.025,.05);
        \draw [->,blue](.025,.05) -- (.05,.05);
        \draw [->,blue](.05,.05) -- (.05,.465);
       \draw [->,blue](.05,.465) -- (.05,.88);
           \draw [->,blue](.05,.88) -- (.465,.88);
              \draw [->,blue](.465,.88) -- (.88,.88);
                        \draw [->,red] (1,.874) -- (.94,.874);
          \draw [->,red] (.94,.874) -- (.874,.874);
    \draw [->,red] (.874,.874) -- (.874,.46);
    \draw [->,red] (.874,.46) -- (.874,.051);
    \draw [->,red] (.874,.051) -- (.46,.051);
    \draw [->,red] (.46,.051) -- (.051,.051);
     \draw [->,red] (.051,.051) -- (.051,.46);
    \draw [->,red] (.051,.46) -- (.051,.875);
    \draw [->,red] (.051,.875) -- (.46,.875);
    \draw [->,red] (.46,.875) -- (.875,.875);
\draw[domain=.909:1,smooth,variable=\x] plot ({\x},{9/(\x+1.5946)-2.5946});
       \draw[domain=.364:.909,smooth,variable=\x] plot ({\x},{9/(\x+1.5946)-3.5946});
             \draw[domain=.014:.364,smooth,variable=\x] plot ({\x},{9/(\x+1.5946)-4.5946});
       \draw[domain=0:.014,smooth,variable=\x] plot ({\x},{9/(\x+1.5946)-5.5946});
               \node at (.05,.03) {$F_4$};
      \node at (.31,.26) {$F_3$};
       \node at (.62,.56) {$F_2$};
      \node at (.9,.95) {$F_1$};
         \node at (.96,.8) {$\Box_1$};
      \node at (.68,.8) {$\Box_2$};
            \node at (.2,.8) {$\Box_3$};
      \node at (0.02,.8) {$\Box_4$};
          \node at (.5,-.05) {$N=9, \alpha=1.5946\cdots$\,\,($k=2$, $i=3$)}; 
                        \node at (.5,-.12) {$\T^3(\alpha+1)=\T(\alpha+1)$};       
    \end{tikzpicture}
 $$
\vspace*{-8mm} \caption{\label{fig: N=9, alfa=1.5946}}
\endminipage
\end{figure}

Now suppose $N=2k^2+2k-4$ (note that voor $k=2$ we have $N=8$). Then 
$$
\alpha_u(2k^2+2k-4,2)=\frac{(k+2)\sqrt{D}-(k^3+k^2+2k+4)}{2(k^2+4k+2)}
$$
and
$$
\alpha_{\ell}(2k^2+2k-4,2)=\frac{(k-1)\sqrt{D}-(k^3+4k^2-k-4)}{2(k^2-2k-1)},
$$
with
$$
D=9k^4+18k^3-7k^2-16k.
$$
There are no gaps provided $\alpha_{\ell}(2k^2+2k-4,2)-\alpha_u(2k^2+2k-4,2)>0$. Once more we omit the calculations, finding that this inequality holds for $k\in \N_{\geq 2}$,
so we conclude that there are no gaps in case $N=2k^2+2k-4$. When we replace the number $4$ in $N=2k^2+2k-4$ by larger integers (if possible), there will not be any gaps either: the length of the `interval' $[\alpha_{\ell},\alpha_u]$ would only become more negative. This concludes the proof that if $\I=\Delta_d\cup\Delta_{d-1}\cup\Delta_{d-2}\cup\Delta_{d-3}$ and there is a gap containing $f_{d-1}$ and $f_{d-2}$ in $\I$, then $N=2k^2+2k-i$, with $k > 1$ and $i \in \{1,2,3\}$.\smallskip

For the \emph{converse statement}, we assume that $N = 2k^2 + 2k - i$, with $k \in\N$ and $i \in \{1, 2, 3\}$. If also $d=k+2$, then earlier
in this proof we showed that only then $\alpha_{\ell}(N,2) \leq \alpha_u(N,2)$. We will show that for $\alpha$ such that
$\alpha_{\ell}(N,2) \leq \alpha \leq \alpha_u(N,2)$ there is a gap in $I_{\alpha}=\Delta_d\cup \Delta_{d-1}\cup \Delta_{d-2}\cup\Delta_{d-3}$ containing both $f_{d-1}$ and $f_{d-2}$.\medskip\

As earlier in this proof, we set 
$$
q = [\,\overline{d-1,\, d-2}\, ]_{N,\alpha} \quad \text{and}\quad r = [\,\overline{d-2,\, d-1}\, ]_{N,\alpha}.
$$
Set $G=(q,r)$, then clearly both $f_{d-1}\in G$ and $f_{d-2}\in G$. Furthermore, by definition of $\alpha_{\ell}(N,2)$ and $\alpha_u(N,2)$ we have that for every $\alpha \in [\alpha_{\ell}(N,2), \, \alpha_u(N,2)]$ that
$$
T_{\alpha}(\alpha )\leq q\quad (\text{and therefore $T_{\alpha}^2(\alpha )\geq r$})
$$
and that
$$
T_{\alpha}(\alpha +1)\geq r\quad (\text{and therefore $T_{\alpha}^2(\alpha +1)\leq q$}).
$$
Note that $T_{\alpha}((p_d,q) = (r, \alpha +1)$ and that $T_{\alpha}((r,p_{d-2}) = (\alpha , q)$. Then we have that $T_{\alpha}(G^c) = G^c$, where $G^c$ is the complement of $G$ in $I_{\alpha}$. We are left to show that $G=(q,r)$ is a gap; i.e.\ that for almost all $x\in G$ there exists an $n=n(x)$ such that $T_{\alpha}^n(x)\in G^c$.\medskip

To show this, consider the map $T:I_{\alpha}\to I_{\alpha}$, defined by
\begin{equation}\label{eq:NewExpansion}
T(x) = \begin{cases}
\frac{-x}{p_d-\alpha} + \frac{(\alpha +1)p_d - \alpha^2}{p_d-\alpha}, & \text{if $x\in \Delta_{d}$}; \\
\frac{N}{x} - (d-1), & \text{if $x\in \Delta_{d-1}$}; \\
\frac{N}{x} - (d-2), & \text{if $x\in \Delta_{d-2}$}; \\
\frac{-x}{\alpha + 1 - p_{d-2}} + \frac{(\alpha +1)^2 - \alpha p_{d-2}}{\alpha +1 - p_{d-2}}, & \text{if $x\in \Delta_{d-3}$}.
\end{cases}
\end{equation}
So on $\Delta_d$ and on $\Delta_{d-3}$ we have that $T$ is a straight line segment with negative slope, through $(\alpha,\alpha +1)$ and $(p_d,\alpha)$ on $\Delta_d$, resp.\ through $(p_{d-2},\alpha +1)$ and $(\alpha +1, \alpha )$ on $\Delta_{d-3}$. For $x\in \Delta_{d-1}\cup \Delta_{d-2}$ we have that $T(x) = T_{\alpha}(x)$. To show that $G$ is a gap, it is enough to show the ergodicity of $T$. Then the maximality 
of $G$ follows from the fact that the support of the absolutely continuous invariant measure is $G^c$, since $T_{\alpha}(G^c) = G^c$.  The proof of the existence of the absolutely continuous invariant measure for $T$ and its ergodicity is similar to the proof of Theorem \ref{ergodic}.  Here all branches are complete and the proof is rather simpler. Once we have the ergodicity of $T$, it is obvious that for a.e.~$x \in G$ there exists $n_{0}=n_{0}(x)$ such that $z=T^{n_0}(x) \in G^c$. Then $z$ never returns in $G$ under iterations of $T_{\alpha}$. This finishes the proof of Theorem \ref{gaps four cylinder case}.\hfill $\Box$

We stress that the in case of $N=2k^2+2k-3$ the intervals $[\alpha_{\ell},\alpha_u]$ on which gaps exist may be very small; see Figure \ref{fig: N=21, alfa=2.7123}. On the other hand, in case $N=2k^2+2k-4$, the gaplessness may be a very close call; see Figure \ref{fig: N=20, alfa=2.6124}. Table \ref{table 2} illustrates how fast these differences between $\alpha_{\ell}$ and $\alpha_u$ decrease as $N$ increases:

\begin{table}[h!]
\centering
 \begin{tabular}{ l | c | r }
  {} & $\alpha_{\ell}(N,2)$\phantom{xxx\,} &  $\alpha_u(N,2)$\phantom{xxx} \\
  \hline
  $N=9$ & $1.594119\cdots$ & $1.594686\cdots$\\
  $N=21$ & $2.712252\cdots$ & $2.712310\cdots$ \\
   $N=37$ & $3.776839\cdots$ & $3.776851\cdots$ \\
     $N=57$ & $4.817672\cdots$ & $4.817675\cdots$ \\
  \hline
  $N=8$ & $1.450165\cdots$ & $1.442809\cdots$\\
  $N=20$ & $2.613247\cdots$ & $2.611575\cdots$ \\
   $N=36$ & $3.700989\cdots$ & $3.700407\cdots$ \\
     $N=56$ & $4.756087\cdots$ & $4.755832\cdots$ \\
\end{tabular}
\vspace{3mm}
\caption{The thin thread between having a gap or not}
\label{table 2}
\end{table}

\begin{Remark} 
While a fixed point $f_i$ is repellent for points within $\Delta_i$, the fixed points in two adjacent cylinder sets behave mutually contracting for all other points in these cylinder sets. As a consequence, it may take quite some time before the orbit of points in the full cylinders of gap arrangements with four cylinders leave these full cylinders for the first time. As an example we take the gap arrangement for $k=50$ (according to the notations used above). Then $N=2\cdot50^2+2\cdot50-3=5097$, $d=52$ and $\alpha\approx\alpha_u\approx\alpha_{\ell}\approx49.98019737$. Table \ref{table 3} shows for ten values of $x$ between $\alpha$ and $\alpha+1$ the smallest $n$ such that $\T^n(x)\not\in\Delta_{51}\cup\Delta_{50}$. What is more, there are \emph{uncountably many} $x$ in the gap $(a,b)$ that contains $f_{d-1}$ and $f_{d-2}$ such that $\T^n(x)\in (a,b)$ for all $n\in \N\cup\{0\}$. Indeed, for any sequence $(d_1,d_2,\ldots,d_n,\ldots)$ such that $d_n\in\{d-1,d-2\}$, with $n\in \N$, we have that $x=[d_1,d_2,d_3,\ldots ]_{N, \alpha}\in(a,b)$.
\end{Remark}

\begin{table}[h!]
\centering
 \begin{tabular}{ c | c c c c c c c c c c}
$x$ & $50$ &  $50.1$  & $50.2$ &  $50.3$ & $50.4$ &  $50.5$  & $50.6$ &  $50.7$ & $50.8$ &  $50.9$\\
\hline
$n$ & $5417$ &  $2090$  & $3568$ &  $1123$ & $4776$ &  $185$  & $5816$ &  $16231$ & $5646$ &  $7604$\\
 \end{tabular}
\vspace{3mm}
\caption{The difficulty of leaving the gap: with $N=5097$, $\alpha=49.98019737$, for each of ten values of $x \in [\alpha,\alpha+1]$ the smallest $n$ is given such that $\T^n(x)\not\in\Delta_{51}\cup\Delta_{50}$.}
\label{table 3}
\end{table}

\begin{figure}[!htb]
\minipage{0.44\textwidth}
$$
\begin{tikzpicture}[scale =7] 
\draw[draw=white,fill=gray!20!white] 
plot[smooth,samples=100,domain=0:.0149] (\x,0) --
plot[smooth,samples=100,domain=.0149:0] (\x,.0434);
\draw[draw=white,fill=gray!20!white] 
plot[smooth,samples=100,domain=.0149:.951] (\x,0) --
plot[smooth,samples=100,domain=.951:.0149] (\x,1);
\draw[draw=white,fill=gray!20!white] 
plot[smooth,samples=100,domain=.951:1] (\x,.924) --
plot[smooth,samples=100,domain=1:.951] (\x,1);
\draw[black,fill=black] (.011,.011) circle (.04ex);
\draw[black,fill=black] (.2866,.2866) circle (.04ex);
\draw[black,fill=black] (.605,.605) circle (.04ex);
\draw[black,fill=black] (.97,.97) circle (.04ex);
\draw[black,fill=black] (0,0) circle (.04ex);
    \draw[black,fill=black] (1,0) circle (.04ex);
         \draw (0,0) --  (1,0);
        \draw [dashed] (.0149,0) -- (.0149,1);
                            \draw [dashed] (.412,0) -- (.412,1);
         \draw [dashed] (.951,0) -- (.951,1);
        \draw (0,0) -- (1,1);
     \draw (0,0) -- (0,1);
     \draw (0,1) -- (1,1);
     \draw (1,1) -- (1,0);
        \draw [dotted, ->, blue, very thick] (0,.0434) -- (.0434,.0434);
             \draw [->,blue](.0434,.0434) -- (.0434,.48);
     \draw [->,blue](.0434,.48) -- (.0434,.918);
     \draw [->,blue](.0434,.918) -- (.48,.918);
     \draw [->,blue](.48,.918) -- (.918,.918);
          \draw [->,blue](.918,.918) -- (.918,.47);
     \draw [->,blue](.918,.47) -- (.918,.052);
     \draw [->,blue](.918,.052) -- (.47,.052);
     \draw [->,blue](.47,.052) -- (.052,.052);
      \draw [->,blue](.052,.052) -- (.052,.47);
     \draw [->,blue](.052,.47) -- (.052,.893);
     \draw [->,blue](.052,.893) -- (.47,.893);
     \draw [->,blue](.47,.893) -- (.893,.893);
       \draw [->,blue](.893,.893) -- (.893,.49);
     \draw [->,blue](.893,.49) -- (.893,.092);
     \draw [->,blue](.893,.092) -- (.49,.092);
     \draw [->,blue](.49,.092) -- (.092,.092);
         \draw [->,blue](.092,.092) -- (.092,.43);
     \draw [->,blue](.092,.43) -- (.092,.782);
     \draw [->,blue](.092,.782) -- (.43,.782);
         \draw [->,blue](.43,.782) -- (.782,.782);
         \draw [->,red] (1,.924) -- (.96,.924);
             \draw [->,red] (.96,.924) -- (.924,.924);
    \draw [->,red] (.924,.924) -- (.924,.46);
    \draw [->,red] (.924,.46) -- (.924,.043);
    \draw [->,red] (.924,.043) -- (.48,.043);
    \draw [->,red] (.48,.043) -- (.043,.043);
     \draw [->,red] (.043,.043) -- (.043,.48);
    \draw [->,red] (.043,.48) -- (.043,.92);
    \draw [->,red] (.043,.92) -- (.48,.92);
    \draw [->,red] (.48,.92) -- (.92,.92);
       \draw [->,red] (.92,.92) -- (.92,.48);
    \draw [->,red] (.92,.48) -- (.92,.05);
    \draw [->,red] (.92,.05) -- (.48,.05);
    \draw [->,red] (.48,.05) -- (.05,.05);
   \draw [->,red] (.05,.05) -- (.05,.48);
    \draw [->,red] (.05,.48) -- (.05,.9);
    \draw [->,red] (.05,.9) -- (.48,.9);
    \draw [->,red] (.48,.9) -- (.9,.9);
     \draw [->,red] (.9,.9) -- (.9,.49);
    \draw [->,red] (.9,.49) -- (.9,.08);
    \draw [->,red] (.9,.08) -- (.49,.08);
    \draw [->,red] (.49,.08) -- (.08,.08);
       \draw [->,red] (.08,.08) -- (.08,.45);
    \draw [->,red] (.08,.45) -- (.08,.808);
    \draw [->,red] (.08,.808) -- (.45,.808);
    \draw [->,red] (.45,.808) -- (.808,.808);
   \draw[domain=.951:1,smooth,variable=\x] plot ({\x},{20/(\x+2.6124)-4.6124});
     \draw[domain=.412:.951,smooth,variable=\x] plot ({\x},{20/(\x+2.6124)-5.6124});
        \draw[domain=.0149:.412,smooth,variable=\x] plot ({\x},{20/(\x+2.6124)-6.6124});
             \draw[domain=0:.0149,smooth,variable=\x] plot ({\x},{20/(\x+2.6124)-7.6124});
     \node at (.05,.025) {$F_5$};
      \node at (.33,.287) {$F_4$};
       \node at (.65,.605) {$F_3$};
      \node at (.92,.97) {$F_2$};
         \node at (.975,.8) {$\Box_2$};
      \node at (.68,.8) {$\Box_3$};
            \node at (.23,.8) {$\Box_4$};
      \node at (0.0075,.8) {$\Box_5$};
                 \node at (.5,-.05) {$N=20, \alpha=2.6124$}; 
    \node at (.5,-.12) {$\alpha_{\ell}=2.6132\cdots$\,\, and\,\, $\alpha_u=2.6115\cdots$}; 
         \end{tikzpicture}
 $$
\vspace*{-8mm} \caption{\label{fig: N=20, alfa=2.6124}}
\endminipage
\hfill
\minipage{0.44\textwidth}
$$
\begin{tikzpicture}[scale =7] 
\draw[draw=white,fill=gray!20!white] 
plot[smooth,samples=100,domain=0:.0106] (\x,0) --
plot[smooth,samples=100,domain=.0106:0] (\x,.03);
\draw[draw=white,fill=gray!20!white] 
plot[smooth,samples=100,domain=.0106:.03] (\x,0) --
plot[smooth,samples=100,domain=.03:.0106] (\x,1);
 \draw[draw=white,fill=gray!20!white] 
plot[smooth,samples=100,domain=.03:1] (\x,.945) --
plot[smooth,samples=100,domain=1:.03] (\x,1);
 \draw[draw=white,fill=gray!20!white] 
plot[smooth,samples=100,domain=.945:.964] (\x,0) --
plot[smooth,samples=100,domain=.964:.945] (\x,.945);
 \draw[draw=white,fill=gray!20!white] 
plot[smooth,samples=100,domain=.03:.945] (\x,0) --
plot[smooth,samples=100,domain=.945:.03] (\x,.03);
\draw[black,fill=black] (.008,.008) circle (.04ex);

\draw[black,fill=black] (.287,.287) circle (.04ex);
\draw[black,fill=black] (.61,.61) circle (.04ex);
\draw[black,fill=black] (.978,.978) circle (.04ex);
\draw[black,fill=black] (0,0) circle (.04ex);
           \draw[black,fill=black] (.03,0) circle (.04ex);
\draw[black,fill=black] (.945,0) circle (.04ex);
\draw[black,fill=black] (1,0) circle (.04ex);
         \draw (0,0) -- (.03,0);
          \draw (.945,0) -- (1,0);
        \draw [dashed] (.0106,0) -- (.0106,1);
                            \draw [dashed] (.416,0) -- (.416,1);
         \draw [dashed] (.964,0) -- (.964,1);
        \draw (0,0) -- (1,1);
     \draw (0,0) -- (0,1);
     \draw (0,1) -- (1,1);
     \draw (1,1) -- (1,0);
        \draw [->,blue] (0,.03) -- (.03,.03);
             \draw [->,blue] (.03,.03) -- (.03,.48);
     \draw [->,blue] (.03,.48) -- (.03,.945);
     \draw [->,blue] (.03,.945) -- (.48,.945);
     \draw [->,blue] (.48,.945) -- (.945,.945);
          \draw [->,blue] (.945,.945) -- (.945,.47);
     \draw [->,blue] (.945,.47) -- (.945,.03);
     \draw [->,blue] (.945,.03) -- (.47,.03);
     \draw [->,blue] (.47,.03) -- (.03,.03);
                   \draw [->,red] (1,.945) -- (.97,.945);
             \draw [->,red] (.97,.945) -- (.945,.945);
    \draw [->,red] (.945,.945) -- (.945,.46);
    \draw [->,red] (.945,.46) -- (.945,.03);
    \draw [->,red] (.945,.03) -- (.48,.03);
    \draw [->,red] (.48,.03) -- (.03,.03);
     \draw [->,red] (.03,.03) -- (.03,.48);
    \draw [->,red] (.03,.48) -- (.03,.945);
    \draw [->,red] (.03,.945) -- (.48,.945);
    \draw [->,red] (.48,.945) -- (.945,.945);
          \draw[domain=.964:1,smooth,variable=\x] plot ({\x},{21/(\x+2.7123)-4.7123});
     \draw[domain=.416:.964,smooth,variable=\x] plot ({\x},{21/(\x+2.7123)-5.7123});
        \draw[domain=.0106:.416,smooth,variable=\x] plot ({\x},{21/(\x+2.7123)-6.7123});
             \draw[domain=0:.0106,smooth,variable=\x] plot ({\x},{21/(\x+2.7123)-7.7123});
     \node at (.05,.025) {$F_5$};
      \node at (.33,.287) {$F_4$};
       \node at (.65,.605) {$F_3$};
      \node at (.93,.975) {$F_2$};
         \node at (.975,.8) {$\Box_2$};
      \node at (.68,.8) {$\Box_3$};
            \node at (.23,.8) {$\Box_4$};
      \node at (0.0075,.8) {$\Box_5$};
                 \node at (.5,-.05) {$N=21, \alpha=2.7123$}; 
                   \node at (.5,-.11) {$\alpha_{\ell}=2.7122\cdots$\,\, and\,\, $\alpha_u=2.7123\cdots$}; 
         \end{tikzpicture}
 $$
\vspace*{-8mm} \caption{\label{fig: N=21, alfa=2.7123}}
\endminipage
\end{figure}

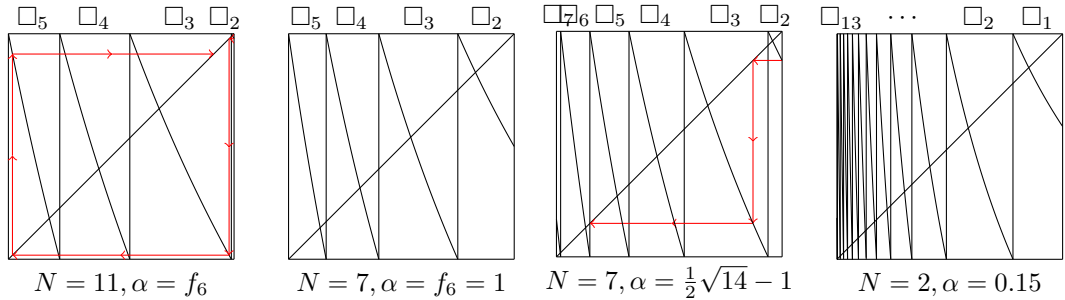
\begin{figure}[!htb]
\minipage{0.24\textwidth}
$$
\begin{tikzpicture}[scale =3]
    \draw (0,0) -- (1,0);
    \draw (.2276,0) -- (.2276,1);
      \draw  (.538,0) -- (.538,1);
        \draw (.988,0) -- (.988,1);
               \draw (0,0) -- (1,1);
     \draw (0,0) -- (0,1);
     \draw (0,1) -- (1,1);
     \draw (1,1) -- (1,0);
     \draw[domain=0.988:1,smooth,variable=\x] plot ({\x},{11/(\x+1.472)-3.472});
      \draw[domain=0.538:.988,smooth,variable=\x] plot ({\x},{11/(\x+1.472)-4.472});
       \draw[domain=.2276:.538,smooth,variable=\x] plot ({\x},{11/(\x+1.472)-5.472});
         \draw[domain=0:.2276,smooth,variable=\x] plot ({\x},{11/(\x+1.472)-6.472});
                   \draw [->,red](1,.978) -- (.978,.978);
                     \draw [->,red](.978,.978) -- (.978,.495);
        \draw [->,red](.978,.495) -- (.978,.018);
   \draw [->,red](.978,.018) -- (.495,.018);
        \draw [->,red](.495,.018) -- (.018,.018);
           \draw [->,red](.018,.018) -- (.018,.46);
   \draw [->,red](.018,.46) -- (.018,.91);
   \draw [->,red](.018,.91) -- (.46,.91);
   \draw [->,red](.46,.91) -- (.91,.91);
    \node at (.96,1.07) {$\Box_2$};
      \node at (.76,1.07) {$\Box_3$};
  \node at (.38,1.07) {$\Box_4$};
   \node at (.11,1.07) {$\Box_5$};
 \node at (.5,-.11) {$N=11, \alpha=f_6$};   
    \end{tikzpicture}
    $$
\endminipage
\minipage{0.24\textwidth}
$$
\begin{tikzpicture}[scale =3]
    \draw (0,0) -- (1,0);
    \draw (.1667,0) -- (.1667,1);
      \draw  (.4,0) -- (.4,1);
        \draw (.75,0) -- (.75,1);
     \draw (0,0) -- (1,1);
     \draw (0,0) -- (0,1);
     \draw (0,1) -- (1,1);
     \draw (1,1) -- (1,0);
     \draw[domain=0.75:1,smooth,variable=\x] plot ({\x},{7/(\x+1)-3});
      \draw[domain=0.4:.75,smooth,variable=\x] plot ({\x},{7/(\x+1)-4});
       \draw[domain=.1667:.4,smooth,variable=\x] plot ({\x},{7/(\x+1)-5});
         \draw[domain=0:.1667,smooth,variable=\x] plot ({\x},{7/(\x+1)-6});
            \node at (.88,1.07) {$\Box_2$};
     \node at (.58,1.07) {$\Box_3$};
  \node at (.28,1.07) {$\Box_4$};
   \node at (.08,1.07) {$\Box_5$};
 \node at (.5,-.11) {$N=7, \alpha=f_6=1$};   
    \end{tikzpicture}
   $$
\endminipage
\minipage{0.24\textwidth}
$$
\begin{tikzpicture}[scale =3] 
      \draw (.01853,0) -- (.01853,1);
       \draw (.148,0) -- (.148,1);
        \draw (.3215,0) -- (.3215,1);
            \draw (.5663,0) -- (.5663,1);
               \draw (.9376,0) -- (.9376,1);
     \draw (0,0) -- (1,0);
      \draw (0,0) -- (1,1);
     \draw (0,0) -- (0,1);
     \draw (0,1) -- (1,1);
     \draw (1,1) -- (1,0);
            \draw [->,red](1,.871) -- (.871,.871);
                     \draw [->,red](.871,.871) -- (.871,.51);
        \draw [->,red](.871,.51) -- (.871,.1484);
   \draw [->,red](.871,.1484) -- (.51,.1484);
        \draw [->,red](.51,.1484) -- (.1484,.1484);
                 \node at (.97,1.07) {$\Box_2$};
          \node at (.75,1.07) {$\Box_3$};
             \node at (.44,1.07) {$\Box_4$};
        \node at (.24,1.07) {$\Box_5$};
      \node at (.08,1.07) {$\Box_6$};
    \node at (.01,1.07) {$\Box_7$};
  \draw[domain=.9376:1,smooth,variable=\x] plot ({\x},{7/(\x+0.8708)-2.8708});
      \draw[domain=.5663:.9376,smooth,variable=\x] plot ({\x},{7/(\x+0.8708)-3.8708});
      \draw[domain=.3215:.5663,smooth,variable=\x] plot ({\x},{7/(\x+0.8708)-4.8708});
       \draw[domain=.148:.3215,smooth,variable=\x] plot ({\x},{7/(\x+0.8708)-5.8708});
        \draw[domain=0.01853:.148,smooth,variable=\x] plot ({\x},{7/(\x+0.8708)-6.8708});
         \draw[domain=0:.01853,smooth,variable=\x] plot ({\x},{7/(\x+0.8708)-7.8708});
          \node at (.5,-.11) {$N=7, \alpha=\tfrac12\sqrt{14}-1$};
    \end{tikzpicture}
$$
\endminipage
\minipage{0.24\textwidth}
$$
\begin{tikzpicture}[scale =3] 
      \draw (.78,0) -- (.78,1);
      \draw (.485,0) -- (.485,1);
     \draw (.332,0) -- (.332,1);
     \draw (.23835,0) -- (.23835,1);
        \draw (.1752,0) -- (.1752,1);
      \draw (.1297,0) -- (.1297,1);
     \draw (.0954,0) -- (.0954,1);
     \draw (.0686,0) -- (.0686,1);
   \draw (.047,0) -- (.047,1);
      \draw (.0294,0) -- (.0294,1);
     \draw (.0146,0) -- (.0146,1);
     \draw (.0021,0) -- (.0021,1);
             \draw (0,0) -- (1,0);
             \draw (0,0) -- (1,1);
     \draw (0,0) -- (0,1);
     \draw (0,1) -- (1,1);
     \draw (1,1) -- (1,0);
   \node at (.02,1.07) {$\Box_{13}$};
              \node at (.3,1.07) {$\cdots$};
  \node at (.63,1.07) {$\Box_2$};
              \node at (.89,1.07) {$\Box_1$};
                  \draw[domain=.78:1,smooth,variable=\x] plot ({\x},{2/(\x+.15)-1.15});
                   \draw[domain=.485:.78,smooth,variable=\x] plot ({\x},{2/(\x+.15)-2.15});
 \draw[domain=.332:.485,smooth,variable=\x] plot ({\x},{2/(\x+.15)-3.15});
 \draw[domain=.23835:.332,smooth,variable=\x] plot ({\x},{2/(\x+.15)-4.15});
 \draw[domain=.1752:.23835,smooth,variable=\x] plot ({\x},{2/(\x+.15)-5.15});
 \draw[domain=.1297:.1752,smooth,variable=\x] plot ({\x},{2/(\x+.15)-6.15});
 \draw[domain=.0954:.1297,smooth,variable=\x] plot ({\x},{2/(\x+.15)-7.15});
 \draw[domain=.0686:.0954,smooth,variable=\x] plot ({\x},{2/(\x+.15)-8.15});
 \draw[domain=.047:.0686,smooth,variable=\x] plot ({\x},{2/(\x+.15)-9.15});
 \draw[domain=.0294:.047,smooth,variable=\x] plot ({\x},{2/(\x+.15)-10.15});
 \draw[domain=.0146:.0294,smooth,variable=\x] plot ({\x},{2/(\x+.15)-11.15});
 \draw[domain=.0021:.0146,smooth,variable=\x] plot ({\x},{2/(\x+.15)-12.15});
 \draw[domain=0:.0021,smooth,variable=\x] plot ({\x},{2/(\x+.15)-13.15});
    \node at (.5,-.11) {$N=2, \alpha=0.15$};
       \end{tikzpicture}
   $$
\endminipage
\caption{Borderline cases for part II of the proof of Theorem \ref{no gaps five cylinder case}} \label{fig: five cylinders}
\end{figure}

For the final, second part of the proof of Theorem \ref{no gaps five cylinder case}, we will consider one by one all cases left, that is $N \in \{2,\ldots,11\}$. When $N=11$ and $\alpha\geq f_7$, $\I$ consists of five cylinder sets if and only if $\alpha\in (f_2-1,f_6)$; see the left arrangement of Figure \ref{fig: five cylinders}, which we already saw in Figure \ref{fig: two two}. Since $\alpha_{\ell}(11,3)>\alpha_u(11,3)$ (cf.~page \pageref{lower}), we conclude on similar grounds as in the proof of Theorem \ref{gaps four cylinder case}, that the arrangement is gapless. When $\alpha \in [f_7,f_2-1]$, the interval $\I$ consists of four cylinder sets, implying gaplessness because of Theorem \ref{gaps four cylinder case}. Since $|\T'(f_7+1)|=2.04\cdots$, gaps are also excluded for all $\alpha\leq f_7$. A similar approach works for $N=10$ (with $|\T'(f_7+1)|=2.03\cdots$), $N=9$ (with $|\T'(f_7+1)|=2.02\cdots$) and even $N=8$, in which case $f_7=1$, $|\T'(f_7+1)|=2$, and the arrangement with four cylinders is full.  \smallskip

For $N \in \{3,\ldots,7\}$ we take a different approach, confining ourselves to the case $N=7$; the cases $N \in \{3,\ldots,6\}$ are done similarly. We will omit most calculations, which are generally quite tedious and do hardly elucidate anything. So let $N=7$. Then $\I$ consists of at least five cylinder sets if and only if $\alpha <f_6$; see the second arrangement of Figure \ref{fig: five cylinders}. We have $|\T'(\alpha+1)|=2$ for $\alpha=\frac12\sqrt{14}-1$, in which case $d_{\min}=2$; see the third arrangement of Figure \ref{fig: five cylinders}. Now suppose $\frac12\sqrt{14}-1\leq \alpha <f_6=1$. We have $|\T'(f_6+1)|=\frac74$ and $|\T'(\alpha)|>7$. Regarding these relatively large values, it is not hard to understand that $\Delta_2$ is gapless. The part of the orbit of $\alpha+1$ under $\T$ in the third arrangement of Figure \ref{fig: five cylinders} illustrates that even in the case of $\alpha=\frac12\sqrt{14}$, the expansion of $[\T(\alpha+1),p_3]$ under $\T$ clearly excludes the existence of gaps.\smallskip

Finally, let $N=2$. We have $|\T'(f_1)|=2$, indicating the rapid increase of $|\T'|$ on $\I$ when $\alpha$ decreases. The large expansiveness of $\T$ left of $f_1$ assures the gaplessness of $\Delta_1$. We will show that for any $\alpha \in (0,\sqrt{2}-1]$ the image of $[\T(\alpha+1),p_2])$ contains most of the fixed points, implying the gaplessness of $\I$; see the last arrangement of Figure \ref{fig: five cylinders} for an illustration of this. When $\T(\alpha+1)\leq f_2$ this is quite obvious, so we assume $\T(\alpha+1)> f_2$. Suppose that $\T^2(\alpha+1)=f_s$, for some $s \in \N_{\geq 2}$. Then, omitting some basic calculations, we have $\alpha=(s+1-\sqrt{s^2+8})/(2s-7)$, whence 
$$
d=d(\alpha)=\left \lfloor \frac{4s^2-11s+(4s-13)\sqrt{s^2+8}-15}{2s-7}\right \rfloor\geq \frac{4s^2-13s+(4s-13)\sqrt{s^2+8}-8}{2s-7},
$$
from which we derive that $d\geq4s$.\smallskip

This finishes the proof of Theorem \ref{no gaps five cylinder case}.\hfill $\Box$

\section*{Acknowledgements}
The first author of this paper is very thankful to the Dutch organisation for scientific research Nederlandse Organisatie voor Wetenschappelijk Onderzoek, which funded his research for this paper with a so-called {\it{Promotiebeurs voor Leraren}}, with grant number 023.003.036. We also want to thank Julian Lyczak, then working at Leiden University, for his help with some programming in the preparatory phase of Section \ref{four and five}. Finally, we thank Niels Langeveld for providing us with Figure \ref{N=50}. The third author's research was partially supported by the Japan Society for the Promotion of Science with Grant-in-Aid for Scientific Research (C) 20K03661.

\end{document}